\documentclass[a4paper,leqno,12pt]{article}
\usepackage{a4wide}
\usepackage{hyperref} 
\usepackage{amsmath,amsfonts,amssymb}
\usepackage{mathrsfs}

\makeindex

\newcommand{\bbb}{\begin{eqnarray}}
\newcommand{\eee}{\end{eqnarray}}

\DeclareMathOperator{\erf}{{\rm erf}}
\DeclareMathOperator{\erfc}{{\rm erfc}}

\DeclareMathOperator{\Ss}{\mathcal S}

\newcommand{\ta}{{}^t\!}

\newcommand{\GL}{{\rm GL}}

\newcommand{\SL}{{\rm SL}}

\newcommand{\Lie}{{\rm Lie\,}}

\newcommand{\R}{{\mathbb R}}
\newcommand{\C}{{\bf C}}
\newcommand{\Z}{{\bf Z}}
\newcommand{\N}{{\bf N}}

\newcommand{\Or}{{\rm O}}

\newcommand{\SO}{{\rm SO}}

\newcommand{\ot}{{\mathfrak o}}

\newcommand{\p}{{\mathfrak p}}

\newcommand{\ka}{{\mathfrak k}}

\newcommand{\g}{{\mathfrak g}}

\newcommand{\Ima}{\, {\rm Im}\,\,}
\newcommand{\Rea}{\, {\rm Re}\,\,}
\newcommand{\sig}{\,\, {\rm sig}\,\,}

\newcounter{EQNR}[NNN]
\newcommand{\nne}{\refstepcounter{EQNR}  \tag{\theNNN .\theEQNR} }
\newcommand{\nn}{\vskip 2mm 
\refstepcounter{NNN}\noindent {\bf \theNNN . }}

\nopagebreak%

 \DeclareMathOperator{\Ei}{Ei}
\DeclareMathOperator{\mXi}{\widetilde{\Xi}_\rho}

\newcommand{\ZZ}{{\mathbb Z}}
 
\newcommand{\HH}{{\mathbb H}}
 \newcommand{\kzxz}[4]{\left(\begin{smallmatrix} #1 & #2 \\ #3 & #4\end{smallmatrix}\right) }
 
\def\?{\ ???\ \immediate\write16{}%
\immediate\write16{Warning: There was still a question mark . . . }%
\immediate\write16{}}

\begin{document}
\title{On Kudla's Green function for signature (2,2), part I}
\author{by Rolf Berndt and Ulf K\"uhn}
\date{}
\maketitle
 
\begin{abstract} Around 2000 Kudla presented conjectures about deep relations between
arithmetic intersection theory, Eisenstein series and their derivatives, and special values of Rankin $L-$series.
The aim of this text is to work out the details of an
old unpublished draft on the second  authors attempt to prove these conjectures for the
case of the product of two modular curves. 

The mayor difficulties in our situation are of analytical nature, therefore this text assembles some material concerning Kudla's Green function
associated to this situation. We present a mild modification of this
Green function  that satisfies the
requirements of being a Green function in the sense of Arakelov
theory   on the natural compactification  in addition. Only this allows us to define arithmetic special cycles. We then prove
that the generating series of those modified arithmetic special cycles is as predicted by Kudla's conjectures 
a modular form with values in the first arithmetic Chow group. 
 
In order to make our work accessible for non experts in the
theory of orthogonal groups, we give detailed description of various
not well-documented facts needed in our calculation in several
appendices.  
\end{abstract}
%
 
\tableofcontents

\newpage

 \addcontentsline{toc}{section}{\bf Introduction}

\section*{Introduction}
\setcounter{section}{0}

This is the first part of our attempt to prove a particular case of Kudla's conjecture on the relation of arithmetic intersection numbers on Shimura varieties of type $\Or(n,2)$, Eisenstein series and their derivatives, and Rankin $L$-series.
In \cite{Ku3} this general picture is split into
several sub-conjectures each of independent interest.

The first  of it is to prove that a generating series for particular
arithmetic special cycles is a modular form 
with values in an arithmetic Chow group. The second 
of these problems is then to identify the derivative of an Eisenstein
series with the arithmetic intersection product of that particular
arithmetic generating series with some power of the first arithmetic
Chern class of the line bundle of modular forms.

   The product of two modular curves $X(1)$ is a particular case of a Shimura
variety associated to the orthogonal group $\Or(n,2)$. Allthough its geometry 
and arithmetic is comparably simple, due to the co-rank 2 situation the 
analytic problems are expected to be more involved then in the general case. 
  
A major open aspect of Kudla's conjectures is to deal with the behaviour at the compactification divisor of the Green's function he proposed, since on any natural compactification these Green functions in general fail to be Green's function in the sense of Arakelov theory of Gillet-Soul\'e or its extension by Burgos-Kramer-K\"uhn.   One has at least two ways to overcome these difficulties: either to develop a meaningful theory of push-forwards in arithmetic intersection theory with compact support, or to find a consistent way to cut off the severe singularities at the boundary. Inspired by recent work of Funke-Millson we present here in the case of the product of the modular curve with itself a modification of Kudla's Green functions depending on a choice of a certain partition of unity. The generating series of differential
forms associated with those modified Green's functions remains modular.  

Using the general machinery of arithmetic intersection theory now allows us to 
prove the main result of this article.
It settles the first sub-conjecture, which in our case can be seen as
  an arithmetic version of the Hirzebruch-Zagier theorem.\\

{\bf Main Theorem.} \emph{For any choice of $\rho$   the generating function of modified special arithmetic cycles $\widehat{Z}_\rho(v,z,m)$ 
\begin{align*} 
\widehat \phi_{K,\rho} := \sum_{m \in \ZZ}  \widehat{Z}_\rho(v,z,m) q^m  
\end{align*} 
is an $\widehat{CH}^1(X) $-form valued  weight $2$ modular form for the group $\SL_2(\ZZ)$, here  $q = e^{2\pi i\tau }$,  $v = \Ima \tau >0$ 
and more details for the definition of $\widehat{Z}_\rho(v,z,m)=(T(m),\mXi(v,z,m)) \in \widehat{CH}^1(X)$ is given below.}\\

We intend to show how the arithmetic intersection
theory with pre-log-log forms developed in
\cite{BKK} can now be applied to compute various
arithmetic intersection numbers involving the line bundle of modular
forms equipped with the natural invariant metric on the product of two
modular curves.  Those arithmetic applications will be discussed in a following paper.

We point to the fact that the principles and methods of proof are
easily generalized to other lattices and groups of type $\Or(2,n)$ and in particular a similar result 
hold also for Hilbert modular surfaces.\\

{\bf Discussion of the results and leitfaden of the paper}\\

In section 1 we recall the basic results on Hecke corresponences  $T(m)$ in the product of modular curves $X= X(1)\times X(1)$. Here we have as usual $X(1) = \SL_2(\ZZ) \setminus \HH \cup \infty $.\\

In section 2 we introduce Kudla's Green function
$$  
\Xi (v,z,m) = \Xi (v,z_1,z_2,m) = (1/2)\sum_{M\in L_m} - \operatorname{Ei}(- \pi (v/y_1y_2)|a-bz_2-cz_1+dz_1z_2|^2)
$$
 associated with $T(m)$ where 
$L=\operatorname{M}_2(\ZZ)$ is the lattice of integral matrices and
$L_m= \{M = (\begin{smallmatrix}a&b\\c&d \end{smallmatrix})\in L:\, \det(M)=m\}$ be the subset of matrices of
determinant $m \in\ZZ$. Here we used the notation $z=(z_1, z_2) \in \HH^2$ and $y_1 = \Ima z_1$ and $y_2 = \Ima z_2$ 
and $v$ is the imaginary part of a symplectic variable $\tau$ coming in when we discuss generating and modular functions.

Its behaviour at the boundary divisor
\begin{align*}
D = \infty \times X(1) + X(1)\times \infty = \{q_1 q_2 =0 \} \subset X,
\end{align*}
where as  usual  $q_1= e^{2 \pi i z_1}, q_2= e^{2 \pi i z_2}$,
is given by\\

{\bf Proposition.} \emph{Near the compactification divisor $D$ up to a smooth function vanishing along $D$ we have}
\begin{align*}
\Xi(v,z_1,z_2,m) \sim \operatorname{L}(q_1,q_2,m)
+ (t/\sqrt v) ( \operatorname{B}(v,s,m) - \operatorname{I}(v,s,m)) ,
\end{align*}
where $t=\sqrt{y_1y_2}$ and $s= \sqrt{y_1/y_2}$,
\begin{align*}
\operatorname{L}(q_1,q_2,m)
&=  \begin{cases} (1/2)\sum_{-bc=m} 
- \log|q_1^{|c|}-q_2^{|b|}|^2, &\quad bc>0\\
0 ,& \quad bc\leq 0 
\end{cases} \\
 \operatorname{B}(v,s,m)&=  (1/2) \sum_{-bc=m}  \frac{2 \sqrt{\pi}}{\sqrt{v}} \operatorname{erfc}_1(|b/s -cs | \sqrt{v})\\
\intertext{and}
\operatorname{I}(v,s,m) &=  (1/2)  \begin{cases} \sum_{-bc=m} 
4 \pi \sqrt v \min(|b/s|, |cs|) , &\quad -bc>0\\
0 , &\quad -bc\leq 0 .
\end{cases}
\end{align*}

The proof of this statement will be a long calculation whose main steps are
a careful analysis of the following sub sum in $\Xi(v,z,m)$ 
\begin{align*}
 \sum_{-bc=m}\sum_{a \in \Z} - \operatorname{Ei}(-2 \pi v |a -bz_2 - cz_1|^2/ (2 y_1 y_2))
 \end{align*}
  by means of Poisson summation w.r.t. $a$. Then the detection and separation
of the singular terms is very analogous to  computations in the context of  local Borcherds products.\\

The  singularities along $D$ of the function 
$ (t/\sqrt v) ( \operatorname{B}(v,s,m) - \operatorname{I}(v,s,m))$ seem to prevent $\Xi _0(v,z,m)$ to be a Green function in any reasonable arithmetic Chow theory. This is not astonishing since for the Kudla-Millson forms 
it is well known that
$$
[-\varphi_{KM}(m)] =[ \operatorname{dd}^c \Xi(m)] \neq [T(m)] \in H^2(X,\Z).$$

However we were able to circumvent these problems with the following method:\\

Let $\rho : X \to \R$ be a partition of unity w.r.t. $D$, i.e. $\rho$ is a smooth function and there are open sets
$U_0 \supsetneq U_1 \supset D$ 
 with  $\rho \equiv 0$ on $X \setminus U_0$ and $\rho\equiv1$ on $U_1$.\\

{\bf Theorem A.} \emph{Let $m \neq 0$. For any choice of $\rho$ the modified Green functions 
$$ \mXi(v,z,m):=  \Xi(v,z,m)  - \rho(z) (t/\sqrt v) (\operatorname{B} (v,s,m)- \operatorname{I}(v,s,m) )
$$
are Green functions for $T(m) \subset X$.} \\

Since the zeroth Green function is a little bit more complicated, it gets devoted the section 3.
It is natural to define the zero term of Kudla's Green function on $X$ by
\begin{align*}
\Xi(v,z,0) &:= (1/2)\sum_{M\in L_0\setminus 0}  - \operatorname{Ei}(-2\pi v \,R(z,M))  + (1/2)\log \|\Delta (z_1)\Delta (z_2)\| ^{1/6}.
 \end{align*}
In analogy to the cases $m\neq 0$  this function will never be a 
Green's function for the divisor $T(v,0)$ appropriate to this case, since it has no logarithmic growth along the 
divisor $T(v,0):= (-1/24 + 1/(8 \pi v) ) \cdot D$ and moreover it has some additional severe singularities along $D$. But we can correct this failure.\\

{\bf Theorem B.} \emph{For any choice of $\rho$ the modified  functions 
$$ \mXi (v,z,0):=  \Xi(v,z,0)  - \rho(z) (t/\sqrt v) \left(\operatorname{B}(v,s,0) - \frac{(1/2)}{\sqrt v} \left( s +\frac{1}{s} \right) \right)
$$
are Green functions for $T(v,0) \subset X$.}\\

In section 4 we give a glance at why our choices made in section 2 and 3 are the right one. 
Indeed, 
most striking is our next result on our modification of Kudla's Green functions, which shows that our modifications are for all choices of the partition of unity $\rho$  in its definition still 
modular in a certain sense. For this we set 
\begin{align*}
\varphi_\rho (v,z,m):= - dd^c \mXi (v,z,m)
\end{align*}

{\bf Theorem C.} \emph{For any choice of $\rho$ and with $q = e^{2\pi i\tau }$ the generating function 
\begin{align*}
\sum_{m \in \ZZ} \varphi_\rho(v,z,m) q^m  
\end{align*} 
is an $A^{1,1}(X) $-form valued  weight $2$ modular form for the group $\SL_2(\ZZ)$.}\\

The proof of the main theorem fills section 5. Our Shimura variety $X$ and the cycles $T(m)$ have obvious models over the integers. 
As indicated before theorems A and B are needed to define the arithmetic 
special cycles $\widehat{Z}_\rho(v,z,m)=(T(m),\mXi(v,z,m))$ as elements  of the  arithmetic Chow group $\widehat{CH}^1(X)$ introduced by Gillet and Soule. 
A comparison of the generating series $\phi_{K,\rho}$ with another 
arithmetic generating series that is modular by trivial reasons \cite{BKK1} 
plus theorem C reduces the proof of modularity to a simple height calculation.\\

In an appendix we give a dictionary to some language from the orthogonal world, 
 detailed proof of facts that have been only sketched in the literature, and an independent proof of modularity of some Kudla-Millson theta series.
 We also relate our considerations to recent work of Funke-Millson 
on a new proof of the Hirzebruch-Zagier theorem.\\

The idea to consider this exceptional example of arose from discussion with Steve Kudla going back to the years 1999-2002. Beginning 2002 the second author wrote an unpublished draft, where a "good behavior" of Kudla's Green function with respect to an reasonabe arithmetic Chow theory was part of the assumptions. This approach
is totally different to that of this paper and it seems unclear whether it could work at all.  Our understanding of the analysis on orthogonal groups  and some proofs in this text largely took profit by conversations with Jan Bruinier and Jens Funke.   We thank them all for their patience with our questions.\\

\section{Products of modular curves and Hecke Correspondences}
Here we follow the language of the classical description as in Hirzebruch's Mannheim lecture \cite{Hi} and report on some analytic theory 
(See also the Argos Seminar \cite{GR}).\\

\nn {\bf Definition.}  Let $\mathbb{H}$ denote the upper
half plane, and $ X(1)={\rm SL}_{2}(\mathbb{Z})\backslash\mathbb{H}
\cup \{{\infty}\} $ the modular curve with the cusp ${\infty}$.
We denote by $X$ the modular surface $X(1) \times X(1)$. 
The divisor
$D= X(1) \times {\infty} + {\infty} \times X(1)$ is called the boundary 
of $X$. \\

Another description of $X$ may be got by considering a homgeneous space of the 
orthogonal group $\Or(2,2).$ For details we refer to our volominous appendix. 
In the text below we will assume at several places knowlwdge of these two `languages'.\\

 Let $N$ be a
positive integer and let $L_{N}$ be the set of integral $(2\times 2)$-matrices of determinant $N\neq 0$. 
Then the Hirzebruch-Zagier divisors, also referred to as Heegner divisors or Hecke correspondences, 
are given by
\begin{align}
\label{T(N)}
T(N)=\left\{(z_{1},z_{2})\in\mathbb{H}\times\mathbb{H}\,\big|
\,\exists\,\gamma\in L_{N}:\,z_{1}=\gamma z_{2}\right\} \subset \HH \times \HH.\nne
\end{align}

These  divisors descent to curves in $X$.
If $N<0$ then $T(N) = \emptyset$ and if $N>0$ the normalisation of $T(N)$
  is a
finite sum of modular curves $X_{0}(m)$ associated to the congruence
subgroups $\Gamma_{0}(m)$ for $m$ dividing $N$. These curves are also
referred to as the graphs of the Hecke correspondences.

Observe if $\gamma =\left(
 \begin{smallmatrix} a & b\\ c & d\end{smallmatrix}\right) \in 
L_N$ then $z_1 = \gamma z_2$ is equivalent to
\begin{align}\label{eq:TN}
a z_2 - d z_1 - c z_1 z_2 +b = 0.\nne
\end{align}
Looking at the real part of \eqref{eq:TN} it is easy to see that only
for $c=0$ such an equation can hold for large $y_1 y_2$ and bounded
$x_1$, $x_2$, where $x_i= \Rea z_i$ and $y_i=\Ima z_i$.  Therefore the
local equation for $T(N)$ near the cusp $\infty\times \infty$ is
given by
\begin{align}
\left\{ q_1^a -q_2^d =0 \,| \qquad ad=N \right\},\nne
\end{align}
where $q_j= e^{2\pi i z_j}$.\\

Another description of $T(N)$ is by means of the modular polynomial. Recall that
the group ${\rm SL}_{2}(\mathbb{Z})$ acts form the left on the set
$L_{N}$ and that a complete set of representatives for this action is
given by the set
\begin{displaymath}
R_{N}=\left\{\gamma=\begin{pmatrix}a&b\\0&d\end{pmatrix}\,\bigg|\,
a,b,d\in\mathbb{Z};\,ad=N;\,d>0;\,0\le b<d\right\}.
\end{displaymath}
and the cardinality of $R_{N}$ equals $\sigma_1(N)= \sum_{d|N} d$.
Denoting by $\Delta(z)$ the discriminant and by $j(z)$ the $j$-function, then the modular form 
\begin{align}\label{def:PsiN}
\Psi_{N}(z_{1},z_{2})=\left( \Delta(z_{1})\Delta(z_{2})\right)^{\sigma_1(N)}
\prod_{\gamma\in R_{N}}\left(j( z_{1})-j(\gamma z_{2})\right)\nne
\end{align} 
for the group ${\rm SL}_{2}(\mathbb{Z})\times {\rm SL}_{2}(\mathbb{Z})$ has divisor 
$T(N) \subset X$. Observe the boundary components of the divisor determined by 
the modular polynomial are compensated by the powers of the discrimant.

Irreducible components of $T(N)$ correspond to decomposition of $R_N$ in square-free $R_n$

We write 
\begin{align*}
[T(N)] \in H^2(X,\ZZ)
\end{align*}
for its class in the middle cohomology group of the surface $X$.\\

Let $W$ be a $\mathbb R$-vector space, e.g. $H^2(H,\mathbb R)$,  and
\begin{align*}
f(q)=\sum_{n\in \Z} w(n,\tau) q^n \in W((q))
\end{align*} 
be a formal $q$-series with $W$-coefficients.  If there exists a
decomposition $ f = \sum_{j \in J} w_j \otimes_\mathbb R f_j $ where each $w_j
\in W$ and each $f_j$ is a modular form with  $q$-coefficients $a(n,z)
\in \mathbb R$, then $f$ is is called a modular form with values in $W$.

\nn {\bf Theorem:} (Hurwitz, Hirzebruch, et.~al.) \label{thm:hurwitz-hirzebruch}
\emph{If we define 
\begin{align*}
T(0) := T(v,0) := \left( \frac{-1}{24} + \frac{1}{8 \pi v} \right) D, 
\end{align*}
then the
generating series for the cohomology classes
$[T(N)] \in H^2(H,\mathbb R)$ of the Hecke correspondences}
\begin{align*}
[\phi_{can}]= [T(v,0)] + \sum_{N>0} [T(N)] q^N  \qquad ( q = e^{2 \pi i \tau} , v = \Ima \tau)  
\end{align*}
\emph{is a weight $2$ modular form for $\SL_2(\ZZ)$ with values in $H^2(X,\mathbb R)$.}\\

{\bf Proof.} One checks that 
 $T(N)$ is rationally equivalent to $\sigma_1(N) T(1)$, where $\sigma_1(N)=\sum_{d|N} d$ and 
$T(0)$ is rationally equivalent to $ \left( \frac{-1}{24} + \frac{1}{8 \pi v} \right) T(1)$, thus
 we then derive
\begin{align*}
[\phi_{can}]= [T(1)] \otimes E_2(\tau,1).
\end{align*}
Here 
\begin{align}\label{def:E2-hecke}
E_2(\tau,1) := -\frac{1}{24} +\frac{1}{8 \pi v} +
\sum_{m>0} \sigma_1(N)\, q^N \nne
\end{align}
is Hecke's non-holomorphic Eisenstein series for $\SL_2(\ZZ)$
of weight $2$.\hfill $\square$\\

{\bf Observe:}  Any linear functional  on 
$H^2(X,\mathbb R)$ applied to $[\phi_{can}]$
yields a modular form (namely a multiple of the Eisenstein series).

\newpage
\section{Some results on Kudla's Green function}

Here we use different coordinates than those in the previous section. 
The translation and some kind of reason for this as well as proofs for standard facts will be given in the appendix,\\ 

{\bf Definition of Kudla's Green function}\\

Let $M =\left(
 \begin{smallmatrix} a & b\\ c & d\end{smallmatrix}\right) \in 
\operatorname{M}_2(\R)$ and $z=(z_1,z_2) \in  \HH^2$, then we write
\begin{align}
\label{R(z,M)}
R(z,M) = R(z_1,z_2,M)= 
\frac{|a   - b z_2 - c z_1  +d z_1 z_2 |^2}{2 y_1y_2} .\nne
\end{align}

\nn {\bf Lemma.} \label{lem:R-trafo} 
(i)\emph{ We have $R(z,M)=0$, if and only if 
$z_1 = \frac{-b z_2 +a }{-d z_2 +c}$.}

(ii)\emph{ Let $\gamma \in \operatorname{SL}_2(\R)^2$, 
then we have $R(\gamma z,M) = R(z, M^\gamma)$, where
$\gamma z =(\gamma_1 z_1,\gamma_2 z_2)$ is given by the usual action
of $\operatorname{SL}_2(\R)^2$ on $\HH^2$ and $M^\gamma= \gamma _1M{}^t\!\gamma _2$.}\\

{\bf Proof.} See Appendix (\ref{inv}). \hfill $\Box$\\

Also in the Appendix we will more systematically treat some here relevant $\Or(2,2)-$theory. For the moment 
we only indicate the following notions which we shall need: 
For $M$ as above, one has the quadratic form
\begin{align}
\label{(M,M)}
(M,M) := 2\det M = 2(ad-bc)  .\nne
\end{align}
All the time, we shall identify
\begin{align}
\label{(M,a)}
\operatorname{M}_2(\R) \ni M =\left( \begin{smallmatrix} a & b\\ c & d\end{smallmatrix}\right) \longleftrightarrow 
{\bf a} = \ta(a,b,c,d) \in \R^4 .\nne
\end{align}
With the usual $SL_2(\R)-$ matrices $g_z = (\begin{smallmatrix} y^{1/2} & xy^{-1/2}\\ &y^{-1/2} \end{smallmatrix})$ such that $g_z(i) = z \in \HH $ 
 a 4 by 4-matrix $A(z)$ is defined by
\begin{align}
\label{(M,z)}
 M^z:= g_{z_1}M \ta g_{z_2} \longleftrightarrow {\bf a}(z) =: A(z)^{-1}{\bf a} .\nne
\end{align}
 (see (\ref{Az})), For ${\bf a}$ we put
 \begin{align}
\label{(Norma)}
\| {\bf a} \| := \sqrt {a^2+b^2+c^2+d^2}.\nne
\end{align}
Hence, one can propose
 
\nn {\bf Lemma.\label{(LemNorma)}} \emph{For the majorant $(M,M)_z$ of the quadratic form given by (\ref{(M,M)}) in $z\in \HH^2,$ 
we have the identities}
$$
(M,M)_z = 2R(z,M) + (M,M) =  {}^t(A(z)^{-1}{\bf a})(A(z)^{-1}{\bf a}) = \| {\bf a}(z)\|^2.
$$
{\bf Proof.} See (\ref{maj}) and (\ref{Key}) in the Appendix . \hfill $\Box$

\nn {\bf Definition.} \label{DefKugr} Let
$L=\operatorname{M}_2(\ZZ)$ be the lattice of integral matrices and
$L_m= \{M \in L:\, \det(M)=m\}$ be the subset of matrices of
determinant $m \in\ZZ$.  For $m \in \ZZ, m\not=0,$ define  {\it Kudla's Green function}
 $\Xi(m)$ (c.f. Kudla \cite{Ku1} p.330 and the footnote \footnote{Observe, this is not exactly what Kudla did. Motivated by the theory of theta functions, which will be recalled later, Kudla 
prefers to consider the Green function weighted by an exponential factor
\begin{align}
\label{xi}\nne
\Xi (z,m) := - (1/2)\sum _{M \in L_m}\operatorname{Ei}\left( - 2 \pi R(z_1,z_2,M)  \right) \exp(-\pi(M,M)).
\end{align}
Here we left out the exponential factor, it will occur naturally in the context of generating series, see eg. corollary \ref{CorKUDLAGREEN}.} below)
  by 
 \begin{align}
\label{Xi0}
\Xi(v,z,m) := (1/2)\sum_{M \in L_m} - \operatorname{Ei}(-2\pi v \,R(z,M)) \nne
\end{align}
with
\begin{align*}
- \operatorname{Ei}(-2\pi vR(z,M)) = - \operatorname{Ei}(-2\pi vR(z,-M))
=\int_1^\infty e^{-2\pi vR(z,M)r} dr/r
\end{align*}
  and, moreover, 
  following the prescription in appendix (\ref{taux}), we have introduced an additional positive real parameter $v,$ the imaginary part of 
a symplectic variable $\tau \in \HH$, which we will  mentioned again in this section.\\

The definition also makes sense for $m=0$ if one excludes the matrix $M=0.$ 
This shall be discussed later separately, hence, up to further notion, we shall assume $m\not=0$.

\nn {\bf Proposition.} \label{PropKudla's Green}\emph{If $m\neq 0$, then Kudla's Green function $\Xi (m) := \Xi(v,z,m)$
is a Green function for $T(m)$
 on the open variety $X \setminus D$.}\\

{\bf Proof.} By the above Lemma \ref{lem:R-trafo} the function $\Xi$ is, once we showed the convergence of the series used to define it, invariant under the action of $\Gamma$ and, hence, in the region of convergence, a function on $X.$\\

For the convergence, we remind that one has  
 \begin{align}
\label{Eibound}
 - {\rm Ei}\,(-t) < \exp(-t)\nne
\end{align}
 
 for $t> 1 $ (c.f. 5.1.19 Abramovich-Stegun \cite{AS}).
Therefore we have
  \begin{align}\label{Xivzm}
2\cdot \Xi(v,z,m)& =  \sum_{M 
\in L_m} \operatorname{Ei}\left(-
 2\pi v R(z,M)  \right),\nonumber\\
&\le 
\sum_{\overset{M 
\in L_m}{2 \pi v R(z,M) \le 1}
} \operatorname{Ei}\left(
 -2 \pi v R(z,M)  \right)+
\sum_{M 
\in L_m} \exp\left(-
 2 \pi v R(z,M)  \right), \nne
\end{align}
where the last infinite sum may be majorized by a indefinite (Siegel-)theta series, and, hence, we have to show that the 
first sum stays finite  or has logarithmic singularities along $T(m)$. \\

If  we set
$$
n_z(\kappa ) := \sharp \{M\in L_m ; \,\,2R(z,M) \leq \kappa \},
$$
then for any compact $C\subset \mathbb D$ there is an $n = n(C)\in \mathbb N$ such that $n_z(1) \leq n$ for all $z\in C$. And this may be seen like this: 
We have 
\begin{align*}
2R((i,i),M)
&= (a-d)^2+(b+c)^2 = a^2+b^2+c^2+d^2-2(ad-bc),
\end{align*}
and obviously there are only finitely many $M$ with $\det(M) =m$ and 
$$
 a^2 +b^2 +c^2 +d^2  \le  \kappa + 2m,
$$ 
i.e. $n_{(i,i)}(1) =: n_i$ is finite (and roughly bounded by $8(1+2m)^4$).
By Lemma \ref{lem:R-trafo}
\begin{align*}
2R(z,M) &= 2R((i,i),M^z)\\
&= (a'-d')^2+(b'+c')^2 = a'^2+b'^2+c'^2+d'^2-2(a'd'-b'c'),
\end{align*}
 
where $M^z$ is the matrix to be identified with ${\bf a'} = (A(z)^{-1}{\bf a})$ for some $A(z) \in \GL_2(\R)$.
 Now we show that the linear transformation ${\bf a}\rightarrow A(z)^{-1}{\bf a}$ produces only 
a finite number of elements ${\bf a'}$ with (euclidean) norm $\| {\bf a'}\| \leq \kappa + 2 m$. 
We do this in a most crude manner. The linear map $M \rightarrow M^z$ resp. ${\bf a}\rightarrow A(z)^{-1}{\bf a}$ 
is bounded by, say, $N(z)$. Note $N(z)$ depends continuously on $z$ and, hence, one has bounds $S,T$ with $S\leq N(z)\leq T$  
for all $z \in C.$ We have
$$
n_z(1) := \sharp \{M^z\in L_m^z : \,\,2R((i,i),M^z) \leq 1 \} \leq N(z)^4 n_i \leq T^4n_i
$$
for all $z\in C$. Therefore the first sum in \eqref{Xivzm} has only finitely many summands\footnote{Another way to the finiteness of $n_z(\kappa ) $ could be as follows. Let $\lambda (z)$ be the smallest eigen value 
of $A(z)^{-1}.$ Again one has an $\ell$ with $\lambda (z) \leq \ell$ for all $z\in C$ and one has the estimate
$$
n_z(\kappa )\leq \lambda (z)^4n_i \leq \ell^4 n_i.
$$}.\\

For the first sum in \eqref{Xivzm} we recall that the exponential integral $\operatorname{Ei}$ has at $t=0$  the expansion
$$
-\operatorname{Ei}(-t) = - \log(t) -\gamma - \int_0^t \frac{e^r-1}{r}
dr.
$$
and for large $t$ it has exponential decay. 

Therefore, on open sets disjoint to $T(m)$  and  the boundary $D,$ $\Xi(m)$ is
 a well-defined $\Gamma$-invariant
function. It has logarithmic growth along $T(m)$, i.e.,~up to a smooth function near $T(m)$ one has
\begin{align*}
\Xi(v,z,m) \approx - \log|a  - b z_2 - c z_1 + d z_1z_2|^2 .
\end{align*}

In proposition \ref{PropKUDLA} in the appendix we shall (re-)prove and later we will use that one has the relation between Kudla's  Green function and 
Kudla-Millson's Schwartz form
\begin{align*}
dd^c \Xi (m)   + \delta_{T(m)} =  - \varphi_{KM}(m) 
\end{align*}
As $\varphi_{KM}(m)$ is known to be a smooth form the claim follows.
\hfill$\Box$
\nn {\bf Corollary.}\label{CorKUDLAGREEN} \emph{The proof of the convergence in the Proposition \ref{PropKudla's Green} also shows 
the convergence of the generating series}
\begin{align}\label{KUDLAGREEN}
\Xi (\tau ,z) := \sum_m \Xi(v,z,m)q^m, \,\,z\in \mathbb H^2 \setminus \cup _{m\in\mathbb N} T(m)\nne
\end{align}
(\emph {$\Xi (v,z,0)$ will be later defined in }Definition \ref{DefXivz0}). \\

{\bf Proof.} The estimation (\ref{Xivzm}) above shows that the series 
may be divided into two parts where the second term is majorized by a Siegel theta series 
and the part coming from the first term via $|q^m| = e^{-2\pi vm}$ is majorized by
$$
\sum_m 8T^4v(1+2m)^4 e^{-2\pi vm}.
$$\hfill$\Box$

This generating series will reappear in our section 4.\\ 
 \newpage
{\bf Kudla's Green function $\Xi $ on $\Gamma \setminus \mathbb D$ at the boundary}\\

We now study the
behaviour of $\Xi(v,z,m)$ along the boundary $D$, we are only interested in the terms that might become singular.\\
We  start with the decomposition
\begin{align}
\label{Xi0d}
\Xi(v,z,m) = \Xi^0(v,z,m) +  \sum_{d=1}^\infty \Xi^d(v,z,m),\nne
\end{align}
where 
$$
\Xi^d(v,z,m)= \sum_{ \kzxz{*}{*}{*}{d} \in L_m} -\operatorname{Ei}\left(- 2 \pi v
R\left(z_1,z_2,\kzxz{*}{*}{*}{d} \right)\right)
$$
and
$$
\Xi^0(v,z,m)= (1/2)\sum_{ \kzxz{*}{*}{*}{0} \in L_m} -\operatorname{Ei}\left(- 2 \pi v
R\left(z_1,z_2,\kzxz{*}{*}{*}{0} \right)\right).
$$

Observe each summand $\Xi ^d(v,z,m)$ for itself is not $\Gamma$-equivariant 
and hence $\Xi ^d(v,z,m)$ does not induce  a function on $X$. 
However it is invariant under the parabolic group $\Gamma_\infty \times \Gamma_\infty$. \\
We point to the fact that the following considerations are purely local.\\

\nn {\bf Lemma.} \label{lem:Xism} \emph{If $d \neq 0$, then $\Xi^d(v,z,m)$ defines a smooth function in a neighbourhood of $D.$
 Moreover $\sum_{d\in\mathbb N}\Xi^d(v,z,m)$ vanishes along $D.$} {\footnote{$\mathbb N$ is natural numbers, i.e., without 0}\\

{\bf Proof.} 
Indeed, clearly, if  $M \in L_m$ with $d \neq 0,$ one has
$$ -\operatorname{Ei}\left(-
 2 \pi v R(z_1,z_2,M)  \right) \to 0$$
  for $t^2= y_1 y_2 \to \infty$. But let us take a closer look.\\

For $t\rightarrow \infty $ and $d\not=0$ one has $R(z,M)\geq 1$ and hence by \eqref{Eibound} we may estimate
\begin{align*}
\sum_{d\in \mathbb N} \Xi ^d(v,z,m) = & \sum_{a,b,c\in \mathbb Z,d\in \mathbb N, ad-cb=m}\int_1^\infty e^{-2\pi vR(z,M)r}dr/r\\
\leq & \sum_{a,b,c\in \mathbb Z,d\in \mathbb N,  ad-cb=m}e^{-2\pi vR(z,M)}.
\end{align*}
With Lemma \ref{(LemNorma)} we get
\begin{align*}\label{eqDdd}
\sum_{d\in \mathbb N} \Xi ^d(v,z,m) e^{-2\pi vm} 
&\leq  \sum_{a,b,c\in \mathbb Z,d\in \mathbb N, ad-cb=m} e^{-2\pi vR(z,M)}e^{-\pi v(M,M)}\\\nonumber
&\leq  \sum_{a,b,c\in \mathbb Z,d\in \mathbb N,  ad-cb=m}e^{-\pi v\| \bf a(z) \|^2 }.\nne
\end{align*}

At first, we want to see what happens at the boundary and we take $z_1=iy_1,z_2=iy_2$ with $y_1$ or $y_2 > T>>0.$  Then  one has 
with $iy=(iy_1,iy_2)$
\begin{align*}\label{eqD0}
 \sum_{d\in \mathbb N} \Xi ^d(v,iy,m)e^{-2\pi vm} 
&\leq  \sum_{a,b,c\in \mathbb Z,d\in \mathbb N} e^{-\pi v(a^2/(y_1y_2)+b^2y_2/y_1+c^2y_1/y_2+d^2y_1y_2)}\\\nonumber
&= (1/2)\vartheta (iv/(y_1y_2))(\vartheta (ivy_1y_2)-1)\vartheta (ivy_1/y_2)\vartheta (ivy_2/y_1).\nne
\end{align*}
Here we use the standard facts concerning the theta function $\vartheta ,$ namely (see for instance \cite{Mu} p.40)  for $v\rightarrow \infty $
one has $\vartheta (iv) \rightarrow 1$ and $|\vartheta (iv)-1| < Ce^{-\pi v}.$ And with the functional equation 
$\vartheta (i/v) = \sqrt v \vartheta (iv) $ for $y_1\rightarrow \infty $ (and fixed $y_2>0$), we get 
\begin{align*}\label{eqD1}
 \sum_{d\in \mathbb N} \Xi ^d(v,iy,m)e^{-2\pi vm} &\leq   (1/2)(y_1/v)\vartheta (iy_1y_2/v)(\vartheta (ivy_1y_2)-1)\vartheta (ivy_1/y_2)\vartheta (iy_1/(y_2v))\\\nonumber
& \leq C(y_1/v)e^{-\pi vy_1y_2} \rightarrow 0.\nne
\end{align*}
Similarly, for $y_2\rightarrow \infty $ and fixed $y_1>0$ we get $ \sum_{d\in \mathbb N} \Xi ^d(v,iy,m)e^{-2\pi vm}  \rightarrow 0.$\\



For general $z$ we go back to \eqref{eqDdd} and write
\begin{align*}
A(z) = A(y) A(y)^{-1} A(z),
\end{align*}
then we easily see that we have a constant $K_T >0 $ with 
\begin{align*}
\| {\bf a}(z)\|^2 = {}^t(A(z)^{-1}{\bf a})(A(z)^{-1}{\bf a}) 
  \geq K_T \,\,{}^t(A(y)^{-1}{\bf a})(A(y)^{-1}{\bf a})
  =K_T  \|{\bf a}(y)\|^2
\end{align*}
for all $z\in U_T^1 = \{ z \in \mathbb H^2 : |x_j| <1/2,\,j=1,2,\, y_1>0, y_2 \geq T \}$. 
We estimate
\begin{align*}
\sum_{d\in \mathbb N} \Xi ^d(v,z,m)e^{-2\pi vm} =&  \sum_{a,b,c\in \mathbb Z,d\in \mathbb N}e^{-\pi v\parallel {\bf a}(z)\parallel ^2}\\
\leq & \sum_{a,b,c\in \mathbb Z,d\in \mathbb N}e^{-\pi vK_T\parallel {\bf a}(y)\parallel ^2}.
\end{align*}
As above in \eqref{eqD0} the last sum will have a finite value and similarly for $z \in U_T^2$ with one and two interchanged.
\hfill$\Box$\\

Thus for $t\gg0$ we can neglect these contributions.
So we are left with analyzing  the remaining term $\Xi^0(v,z,m)$.
First observe, if $d = 0$, then we have
\begin{align}
\label{Rd0}
R(z,\kzxz{a}{b}{c}{0} ) &= (v/(2t^2))((a-(bx_2+cx_1))^2 + (by_2+cy_1)^2)\nonumber\\
&= (v/(2t^2))((a-{\bf x})^2 + {\bf y} ^2) \nne
\end{align}
with
\begin{align}
\label{beta}
{\bf x} := bx_2+cx_1,\,\, {\bf y} := |by_2+cy_1| \nne
\end{align}

Now, if we write 
\begin{align}
\label{hbc}
\Xi^0(v,z,m) = (1/2)\sum _{-bc=m} \left(\sum _{a\in \mathbb Z} 
-\operatorname{Ei}\left(- 2 \pi v 
R\left(z_1,z_2,\kzxz{a}{b}{c}{0} \right)\right) \right) ,\nne
\end{align}
the inner sum  with respect to $a$ has period 1 as a function in ${\bf x}$. 
In order to describe the induced Fourier expansion
\begin{align}
\label{hbc1}
\sum _{a\in \mathbb Z} 
-\operatorname{Ei}\left(- 2 \pi y 
R\left(z_1,z_2,\kzxz{a}{b}{c}{0} \right)\right) 
= \sum_{n \in \Z} a_{bc}(n) e(n{\bf x}),\nne
\end{align}
 we recall the error function 
$$\erf(x) = \frac{2}{\sqrt{\pi}}  \int_0^x e^{-t^2} dt,$$
it satisfies 
$\erf(-x)=-\erf(x)$. We also use the complementary error function $\erfc(x)= 1- \erf(x)$. 
Then we finally we need another function.

\nn {\bf Definition.} We define
\begin{align}\label{Bvs}
\operatorname{B}(v,s;b,c)&=    
\int^\infty _1 e^{-\pi/(t^2){\bf y} ^2vu }u^{-3/2}du
= \int^\infty _1 e^{-\pi ( b/s +c s)^2 vu }u^{-3/2}du.
\nne
\end{align}
\nn {\bf Remark.} Here we adapt the notation (up to a factor the function $\beta (\alpha )$ appearing in [HZ] and in other sources)
\begin{align}
\label{B}
B(\alpha ) :=   \int_1^\infty e^{-\alpha u}u^{-3/2}du \nne
\end{align}
and with
$\alpha := \alpha (s)  = (\pi /t^2){\bf y} ^2 = \pi (bs^{-1}+cs)^2$.
We get $\operatorname{B}(v,s;b,c) = B(\alpha (s)).$
As often we shall use this later, we remark that one has
\begin{align}
\label{B0} \nne
B(0) = 2\,\,\, {\rm and}\,\,\lim_{\alpha \rightarrow \infty } B(\alpha ) = 0.
\end{align}

\nn {\bf Lemma}. \label{lem:fourierexp}
 \emph{ 
The coefficients in the
Fourier expansion \eqref{hbc1} are given for $n=0$ by
\begin{align*} 
a_{bc}(0) &=  \frac{t}{\sqrt v} \operatorname{B}(v,s;b,c) 
\intertext{and for $n \not= 0,$ one has}
\label{abcn}
a_{bc}(n) &= \frac{e^{-2\pi {\bf y} \mid n \mid }}{|n|}  
  -\frac{e^{2\pi {\bf y} \mid n \mid }}{2\mid n \mid }\erfc(\sqrt{\pi }t|n| /\sqrt v + \sqrt{\pi v}{\bf y} /t)  
\nonumber\\ 
& \quad - \frac{e^{-2\pi {\bf y} \mid n \mid }}{2|n|}\erfc(\sqrt{\pi }t|n| /\sqrt v - \sqrt{\pi v}{\bf y} /t ) \nne  
\end{align*}
}
{\bf Proof.}
By the definition of the exponential integral, namely
$ - \Ei(-r) = \int_1^\infty \exp(-tr) t^{-1} dt$, we get
\begin{align}
\label{abc}
a_{bc}(n) &=  \int^1_0(\sum _{a \in \Z} \int^\infty _1 e^{-(\pi/(t
^2))v((a-{\bf x})^2 + {\bf y} ^2)r} dr/r)e(-n{\bf x})d{\bf x}\nonumber\\[.3cm]
&=  \int^\infty _1 \int^\infty _{-\infty } e^{-(\pi/(t^2))vr({\bf x}^2+{\bf y} ^2) -2\pi in{\bf x}}d{\bf x}dr/r\nonumber\\[.3cm]
&=  (t/\sqrt v) \int^\infty _1 e^{-(\pi/(t^2))v{\bf y} ^2u - \pi t^2n^2/(vu)}u^{-3/2}du\nonumber\\[.3cm]
&=   (t/\sqrt v)\int^1_0 e^{-(\pi/(t^2))v{\bf y} ^2/w - (\pi t^2n^2/v)w}w^{-1/2}dw\nonumber\\[.3cm]
&=   2(t/\sqrt v)\int^1_0 e^{-(\pi/t^2)v{\bf y} ^2/x^2 + (\pi t^2n^2/v)x^2)}dx \nne
\end{align}
Here one can use formula 7.4.33 from \cite{AS}: For $a\not=0,$ we have 
\footnote{This comes out as follows: One has
$$
\erf z := (2/\sqrt \pi) \int^z_0 e^{-u^2}du,
$$
hence
$$
(\erf z)' = (2/\sqrt \pi)e^{-z^2}
$$
and
$$
(\erf(ax+b/x))' = (2/\sqrt \pi)e^{-(ax+b/x)^2}(a-b/x^2)
$$
$$
(\erf(ax-b/x))' = (2/\sqrt \pi)e^{-(ax-b/x)^2}(a+b/x^2).
$$
From here we easily get the formula above which, with $\erf (\pm \infty ) = \pm 1,$ leads to
$$
\int^1_0 e^{-a^2x^2-b^2/x^2}dx = \sqrt\pi/(4a)(e^{2ab}(\erf(a+b) - 1)+ e^{-2ab}(\erf(a-b) + 1)).
$$}
$$
\int e^{-a^2 x^2-b^2/x^2}dx = \sqrt\pi/(4a)(e^{2ab}\erf(ax+b/x) + e^{-2ab}\erf(ax-b/x))
$$ 
 Using this identity and $\erfc z = 1 - \erf z$ we get the claim. \hfill $\Box$

\nn{\bf Lemma.}  \label{lem:localexp} {\emph{For $n\neq 0$ 
we we modify the Fourier coefficients in \eqref{hbc1} by
\begin{align*} 
 \tilde a_{bc}(n):= a_{bc}(n) - \frac{\exp(-{\bf y} |n|)}{|n|}  . 
\end{align*}
Then the modified coefficients   satisfy 
\begin{align}\label{localexp}
 \tilde a_{bc}(n)   \le 2\sqrt v \frac{\exp(-\pi t^2|n|^2)}{\pi |n|t},\nne
\end{align}
moreover all those modified coefficients vanish along $D$.}
}\\

{\bf Proof.}
For $t \geq 0,$ from \cite{AS} 7.1.13 we take over the estimation
\begin{align}
\label{7.1.13}
\frac{1}{t+\sqrt{t^2+2}} < e^{t^2}\int_t^\infty e^{-u^2}du \leq \frac{1}{t+\sqrt{t^2+4/\pi}}.\nne
\end{align}
Therefore for $n \not= 0$ in Lemma \ref{lem:fourierexp} with $t_n := \sqrt{\pi t}|n|/\sqrt v + \sqrt{\pi v}{\bf y}/t$ one has
\begin{align*}
 \frac{e^{2\pi {\bf y} |n|}} {2|n|}\erfc(\sqrt{\pi t}|n|/\sqrt v + \sqrt{\pi v}{\bf y}/t ) 
&=  \frac{e^{2\pi {\bf y} |n|}} {2|n|}\erfc(t_n)\\
&=  \frac{e^{2\pi {\bf y} |n|}} {2|n|}\frac{2}{\sqrt \pi }\int_{t_n}^\infty e^{-t^2}dt\\\nonumber
&\leq  \frac{e^{2\pi {\bf y} |n|}} {2|n|}\frac{2}{\sqrt \pi }\frac{e^{-t_n^2}}{t_n+\sqrt{t_n^2+4/\pi}}\\\nonumber
&\leq  \sqrt v \frac{e^{-\pi t^2n^2/v}} {\pi |n|t}.\nne
\end{align*}

One has a similar estimation for the term
\begin{align*}
 \frac{e^{-2\pi {\bf y} |n|}} {2|n|}\erfc(\sqrt{\pi t}|n|/\sqrt v - \sqrt{\pi v}{\bf y}/t ) 
&\leq  \sqrt v \frac{e^{-\pi t^2n^2/v}} {\pi |n|t}\nne
\end{align*}
and hence the claim follows from Lemma \ref{lem:fourierexp}.\hfill$\Box$\\

Hence only the first summand in $a_{bc}(n)$ has to be considered. \\

\nn {\bf Proposition.} \label{PropSbc}\emph{With ${\bf x} = bx_2+cx_1$ and ${\bf y} = |by_2+cy_1|$ we set
\begin{align}\label{Sbc}
\mathscr{L}_{bc}({\bf x},{\bf y}) &:=  \sum _{n \in \Z, n\not=0} \frac{e^{-2 \pi {\bf y} |n|}}{\mid n \mid } \,e^{2\pi in{\bf x}} ,\nne
\intertext{For all ${\bf y}>0$ the function $\mathscr {L}_{bc}({\bf x})$ is a smooth function, in particular}
 \mathscr {L}_{bc}({\bf x},{\bf y}) &= 
- 2\log \left | 1 \,-\!\, e^{2\pi i({\bf x} +i{\bf y})}\right|.\nne
\intertext{If $-bc =m >0$, then $\mathscr {L}_{bc}({\bf x},{\bf y})$ has a logarithmic singularity 
along the irreducible component  of $T(-bc)$ given by
$$\{ {\bf x}={\bf y} =0\} =\{b z_1  + c z_2 =0\} = \{ q_1^{\mid c \mid }   - q_2^{\mid b \mid } =0 \}
 $$
 and an additional singularity in the point
$(\infty,\infty)\in X$. More precisely we have the expansion}   
 \mathscr{L}_{bc}({\bf x},{\bf y}) &= -\log \left | q_1^{\mid c \mid } - q_2^{\mid b \mid } \right |^2 -
2 \min\left(-\log\left|q_1^{\mid c \mid }\right| 
 ,-\log\left|q_2^{\mid b \mid } \right| \right) \nne
\end{align}
If $-bc =m <0,$ then $\mathscr{L}_{bc}({\bf x},{\bf y})$ takes finite values for all $y_1,y_2 > 0.$ }\\

{\bf Proof.}
For $\mid u \mid < 1$ one has the development
$$
- \log (1 - u) = u + u^2/2 + u^3/3 + \dots .
$$
Hence, we get for ${\bf y} = \,\mid by_2+cy_1 \mid \,> 0$
\begin{align}
\label{Sx}
\log (1 \,-\!\, e^{-2\pi({\bf y} -i{\bf x})}) + \log (1 - e^{-2\pi({\bf y} +i{\bf x})}) 
&= - \sum _{n \in \N} e^{-2\pi ({\bf y} +i{\bf x})n}/n - \sum _{n \in \N} e^{-2\pi ({\bf y} -i{\bf x})n}/n\nonumber\\[.3cm]
&= - \sum _{n \in \N} e^{-2\pi n{\bf y}}(e^{-2\pi ni{\bf x}} + e^{2\pi ni{\bf x}})/n\nonumber\\[.3cm]
&= - \sum _{n \in \Z, n\not=0} \frac{e^{-2\pi{\bf y} \mid n \mid }}{\mid n \mid } \,e^{2\pi in{\bf x}} \nonumber\\[.3cm]
&= \,\,-  \mathscr {L}_{bc}({\bf x},{\bf y}).\nne
\end{align}
We see that for large ${\bf y} $ this sum  goes to zero.
Moreover, we remark that for $bc > 0$ (because then ${\bf y} >0$) the sum 
$$
 \mathscr {L}_{bc}({\bf x},{\bf y}) =  - \log (1 \,-\! e^{-2\pi({\bf y} -i{\bf x})}) - \log (1 - e^{-2\pi({\bf y} +i{\bf x})})
$$
goes to a finite value for all $y_1,y_2.$ Hence we have to discuss $ \mathscr {L}_{bc}({\bf x},{\bf y})$ only for $m = -bc >0.$ 
One has $ \mathscr {L}_{bc}({\bf x},{\bf y}) = \mathscr {L}_{bc}(-{\bf x},{\bf y}),$ i.e., $ \mathscr {L}_{bc}({\bf x},{\bf y})$ is real, and, to further analyze this, one introduces the usual $q-$variables
$$
q_j = e^{2\pi i({\bf x}_j+iy_j)}, \,j =1,2.
$$
1. For $by_2+cy_1>0,$  we get 
$$
 \mathscr {L}_{bc}({\bf x},{\bf y}) = - 2\log \mid 1 - q_1^cq_2^b \mid .
$$
1.1 For $b > 0, c < 0$ we have $by_2 > cy_1,$ i.e. also $\mid by_2 \mid > \mid cy_1 \mid$ and 
$$
\begin{array}{rl} 
 \mathscr {L}_{bc}({\bf x},{\bf y}) = - 2\log \mid 1 - q_1^cq_2^b \mid \,\,=& -2\log \mid q_1^{-c} - q_2^b \mid - 2 \log \mid q_1^c \mid \\
 =& -2\log \mid q_1^{-c} - q_2^b \mid - (-4\pi cy_1) .
\end{array}
$$
1.2 For $b < 0, c > 0$ we have $by_2 < cy_1,$ i.e. also $\mid by_2 \mid < \mid cy_1 \mid$ and 
$$
\begin{array}{rl} 
 \mathscr {L}_{bc}({\bf x},{\bf y}) = - 2\log \mid 1 - q_1^cq_2^b \mid \,\,=& -2\log \mid q_2^{-b} - q_1^c \mid - 2 \log \mid q_2^b \mid \\
 =& -2\log \mid q_1^{-c} - q_2^b \mid - (-4\pi by_2) .
\end{array}
$$
2. For $by_2+cy_1< 0,$  we get 
$$
 \mathscr {L}_{bc}({\bf x},{\bf y}) = - 2\log \mid 1 - q_1^{-c}q_2^{-b} \mid .
$$
2.1 For $b > 0, c < 0$ we have $by_2 < -cy_1,$ i.e. also $\mid by_2 \mid < \mid cy_1 \mid$ and 
$$
\begin{array}{rl} 
 \mathscr {L}_{bc}({\bf x},{\bf y}) = - 2\log \mid 1 - q_1^{-c}q_2^{-b} \mid \,\,=& -2\log \mid q_2^{b} - q_1^{-c} \mid - 2 \log \mid q_2^{-b} \mid \\
 =& -2\log \mid q_1^{-c} - q_2^b \mid - (-4\pi (-b)y_2) .
\end{array}
$$
2.2 For $b < 0, c > 0$ we have $cy_1 < -by_2,$ i.e. also $\mid cy_1 \mid < \mid by_2 \mid$ and 
$$
\begin{array}{rl} 
 \mathscr L_{bc}({\bf x},{\bf y}) = - 2\log \mid 1 - q_1^{-c}q_2^{-b} \mid \,\,=& -2\log \mid q_1^{c} - q_2^{-b} \mid - 2 \log \mid q_1^{-c} \mid \\
 =& -2\log \mid q_1^{c} - q_2^{-b} \mid - (-4\pi (-c) y_1) .
\end{array}
$$
Hence for $b > 0, c < 0$ we have
$$ 
 \mathscr L_{bc}({\bf x},{\bf y}) =  -2\log \mid q_1^{-c} - q_2^{b} \mid - 4\pi {\rm min}\,(\mid c y_1\mid ,\mid by_2\mid ) 
$$
and for $b < 0, c > 0$
$$
 \mathscr L_{bc}({\bf x},{\bf y}) = -2\log \mid q_1^{c} - q_2^{-b} \mid - 4\pi {\rm min}\,(\mid c y_1\mid ,\mid by_2\mid ) ,
$$
i.e., finally,  for $-bc = m > 0$
\begin{align}
\label{Sxm}
 \mathscr L_{bc}({\bf x},{\bf y}) = -2\log \mid q_1^{\mid c \mid } - q_2^{\mid b \mid } \mid - 4\pi {\rm min}\,(\mid c y_1\mid ,\mid by_2\mid ) ,\nne
\end{align}

\hfill $\Box$

\nn{\bf Definition.}
 We define
 \begin{align}
 \label{Ivs} 
\operatorname{I}(v,s,b,c)
&:=  \begin{cases} 
(4\pi \sqrt v ) \min (|b/s|,|cs|)   & \textrm{if} \,\, \, bc<0 \\
 0 & \text{else}
 \end{cases}\nne
 \end{align}
 Apparently, for $bc<0,$ one also has 
$$
\operatorname{I}(v,s,b,c) = (\sqrt v/t) \min\left(-\log\left|q_1^{\mid c \mid }\right| ^2
 ,-\log\left|q_2^{\mid b \mid } \right| ^2\right) .
$$ 

We   recall
\begin{align*}
\operatorname{B}(v,s;b,c)  &= 
 \int^\infty _1 e^{-\pi ( b/s +c s)^2 vu }u^{-3/2}du
\end{align*}
and propose the following\\
\nn{\bf Notation.}\label{Nota} We write $ f(z) \sim g(z) $  if near $D$ the difference $f(z)-g(z)$ equals a smooth function vanishing along $D,$ 
and such functions are also refered to as being {\it harmless}.\\

This allows for the formulation of a central result.\\

\nn {\bf Theorem.}\label{ThnearD}
\emph{Near $D$ we have }

\begin{align*}
\Xi (v,z,m) \sim (1/2)\sum_{-bc=m} \left(  
-\log \left | q_1^{\mid c \mid } - q_2^{\mid b \mid } \right |^2
+ \frac{1}{\sqrt v} t
\left( \operatorname{B}(v,s,(b,c)) - \operatorname{I}(v,s,(b,c)) \right)\right)
\end{align*}

{\bf Proof.}  By means of Lemma \ref{lem:Xism}, \ref{lem:fourierexp}, \ref{lem:localexp}, and 
Proposition \ref{PropSbc} we have  
\begin{align}
\label{Xi0m}
2\cdot \Xi ^0(v,z,m) &\sim \sum_{-bc=m}\sum_{n\in\Z} a_{bc}(n)e(n{\bf x})\nonumber\\[.3cm]
&\sim \sum_{-bc=m}(\sum_{n\in\mathbb Z, n\not=0}  \frac{e^{-2\pi {\bf y} |n|}}{|n|}e(n{\bf x}) + (t/\sqrt v)B(v,s;b,c))\nonumber\\[.3cm]
&\sim \sum_{-bc=m}( - 2\log \mid q_1^{\mid c \mid }- q_2^{\mid b \mid } \mid \nne\\
&\,\,\,\quad +\, (t/\sqrt v)B(v,s;b,c) - 4\pi t \min ( \mid cs \mid , \mid bs^{-1} \mid ) ).\nonumber
\end{align}\hfill $\Box$

{\bf Boundary corrections}\\

We first observe the following behavior along $D$ near $T(-bc).$

\nn {\bf Lemma.} \label{lem:boundarysing} \emph{In a small neighborhood of} $\{(\infty ,\infty )\}$ (as intersection point of $D$ and $T(-bc)$),
 \emph{we have}
\begin{align*}
\label{old check function}
(t/\sqrt v)\,\, ( \operatorname{B}(v,s;b,c) - \operatorname{I}(v,s;b,c)) 
&= \,(2t/\sqrt v) - 2\pi (|tb/s| + |tcs|)+ O(\alpha )\\\nonumber
&= 2\sqrt {y_1y_2/v} - 2\pi (|b|y_2+|c|y_1) + O(\alpha ).\nne 
\end{align*}
\emph{Here one has 
\begin{align}\label{alphasmall}\nne 
\alpha :=  \pi v(by_2+cy_1)^2/(y_1y_2) = \pi v(b/s+cs)^2
\end{align} 
is assumed to be small.}\\

{\bf Proof.} Analogously to the reasoning around \eqref{eq:TN} resp. \eqref{eqTm}, $T(-bc)$ near $\{(\infty ,\infty )\}$ is given by
$q_2^b-q_1^c = 0$ i.e., $bx_2+cx_1=0, by_2+cy_1=0.$ Hence, here one has
\begin{align} 
\alpha =  \pi v(by_2+cy_1)^2/(y_1y_2) = \pi v(b/s+cs)^2
\end{align} 
small and for such small $\alpha $ one has
\begin{align}\label{Balpha}\nne 
B(\alpha ) &=  2\,(e^{-\alpha } - 2\sqrt \alpha \int_{\sqrt \alpha}^\infty e^{-s^2} ds)\\[.3cm]
&\simeq   2 - 2\pi \sqrt v \mid b/s+cs \mid + \,O(\alpha ) \nonumber.
\end{align} 

and hence as functions in $s$ and $t$ with $bc\not=0$
\begin{align*} \label{alpha small}
(t/\sqrt v)\,\, ( \operatorname{B}(v,s;b,c) &- \operatorname{I}(v,s;b,c)) \\\nonumber
&= (t/\sqrt v)(2 - 2\pi \sqrt v \mid b/s+cs \mid - 4\pi \sqrt v \min (|b/s|,|cs|) + \dots ) .\nne
\end{align*}


The smoothness of this function in $s,t$ is essentially determined by the following consideration 
(here take $-b/s = v, cs = u$):\\

We put 

\begin{align*} 
B(u,v) := \mid u-v \mid , \,\,
I(u,v) := \begin{cases} \min (\mid u \mid , \mid v \mid )\,\,&{\rm if}\,\,uv > 0\\
\,\,0\,\,&{\rm if}\,\,uv \leq  0.
\end{cases}
\end{align*}

Hence, one has

\begin{align*}
 B(u,v) = \begin{cases} u-v \,\,&{\rm if}\,\,u>v\\
 v-u \,\,&{\rm if}\,\,u\leq v\\
\end{cases},\,\,\,
I(u,v) = \begin{cases}v \,\,&{\rm if}\,\,u>v>0\\
 u \,\,&{\rm if}\,\,v>u>0\\
 -u \,\,&{\rm if}\,\,0\geq u\geq v\\
 -v \,\,&{\rm if}\,\,0\geq v\geq u\\ 
 \,0 \,\,&{\rm if}\,\,uv \leq 0.\\
 \end{cases}
\end{align*}
This adds up to the continuous function
 \begin{align*}
B(u,v) + 2I(u,v) &= \begin{cases}u-v+2v= u+v \,\,&{\rm if}\,\,u>v>0\\
 v-u+2u = u+v \,\,&{\rm if}\,\,v\geq u>0\\
 u-v-2u = -v-u \,\,&{\rm if}\,\,0\geq u>v\\
 v-u-2v = -u-v \,\,&{\rm if}\,\,0\geq v\geq u\\ 
 u-v \,\,&{\rm if}\,\,u \geq  0, v \leq  0\\
 v-u \,\,&{\rm if}\,\,u \leq  0, v \geq  0 
 \end{cases}\\
 &=\,\, \mid u \mid + \mid v \mid 
\end{align*}
which even is smooth for $uv \not= 0.$\hfill$\Box$\\

We observe that $\Xi(v,z,m)$ becomes singular along $D$ and moreover\\

\nn {\bf Proposition.} \label{PropnoGreen}\emph{ Kudla's Green function $\Xi(v,z,m)$ descends neither to a Green function in the sense of Gillet-Soul\'e nor to a Green function with log-log singularities 
in the sense of Burgos-Kramer-K\"uhn for the divisor $T(m)$ on the compact variety $X$.}\\

{\bf Proof.} 
First we recall that $T(m)$ intersects the boundary divisor
$$
D = \{q_1,q_2 : q_1q_2 = 0 \} = (\{\infty \, \}\times X(1)) \cup \, (X(1)\times \{\infty \,\})
$$
only in the point $\{\infty ,\infty \}.$ Thus, if $\Xi(v,z,m)$ were a Green function in the sense of Giller-Soul\'e, 
it must be smooth on $D\backslash \{\infty ,\infty \}.$ This is contradicted by the function 
$$
f(q_1,q_2) := t = (y_1y_2)^{1/2} = (1/2\pi )(\log |q_1|\log |q_2|)^{1/2}
$$
in the boundary behaviour of $\Xi(v,z,m).$ If it were a Green function with log-log singularities along $D$ in the sense of 
\cite{BBK}, then $f(q_1,q_2)$ would be of log-log growth as in Definition 1.2 on p.8 in \cite{BBK}, i.e., it would be bounded by 
$\prod_{j=1,2} \log(\log (1/r_j))^M $ for some $M\in \mathbb N.$ Obviously, this is not possible, thus the claim follows. 
\hfill$\Box$\\

\nn {\bf Remark.} The function $f(q_1,q_2)$ is of pre-log growth along $D$ in the sense of Definition 1.4 on p.9 in \cite{BBK}: 
One has $\partial \log |q_j| = -(2\pi )\partial y = (1/2)dq_j/q_j$ and hence
\begin{align*}
\partial f =& (1/2\pi )(1/2)(\log |q_1|\log |q_2|)^{-1/2}\partial (\log |q_1|\log |q_2|)\\
=& -(1/4)((y_2/t) dq_1/q_1 + (y_1/t) dq_2/q_2).
\end{align*}
Moreover, one has 
\begin{align*}
\bar{\partial } y^{1/2} =& \,\,(1/2)y^{-1/2}\bar{\partial }y = -(1/8\pi )y^{-1/2}d\bar{q}/\bar{q}\\
\bar{\partial } y^{-1/2} =& -(1/2)y^{-3/2}\bar{\partial }y = (1/8\pi )y^{-3/2}d\bar{q}/\bar{q}\\
\end{align*}
and hence
\begin{align*}
\bar{\partial } \partial f = \,\,&(1/16\pi t)(\bar{d}q_2\wedge dq_1/(\bar{q}_2 q_1) + \bar{d}q_1\wedge dq_2/(\bar{q}_1 q_2) \\
&+ (y_2/t^2)dq_1\wedge \bar{d}q_1/|q_1^2| + (y_1/t^2)dq_2\wedge \bar{d}q_2/|q_2^2| ).
\end{align*}
However one may check that $f(q_1,q_2) = t$ is not a Green function for $D.$ Indeed if it were a 
Green function, its pullback to any curve intersecting $D$ would be a Green function for its 
intersection point (respecting the intersection multiplicity $n_0$). If we take the curve $C = \{y_1 = c = \,{\rm const}\,\},$ 
we get locally on $C$ the function $f_c = \log |q_2|^{1/2}$ and 
specializing Proposition 1.13 on p.13 in \cite{BBK} one should have
$$
\lim_{\epsilon \rightarrow 0} \int_{\partial B_\epsilon (0)} d^cf_c = n_0
$$ 
with $B_\epsilon (0)$ an $\epsilon -$disc around zero.
But one has
$$
d^cf_c = (1/4\pi i)(\partial -\bar{\partial })(\log |q|)^{1/2} = (-1/16\pi i)y^{-1/2}(dq/q-d\bar{q}/\bar{q})
$$     
and the limes above comes out as zero.\\

Hence, in order to get a Green function for $T(m)$ on $\bar X$ we need to get rid of the singularities. 
The procedure we propose to solve this issue is the most naive one: 
\emph{cut off the singularities smoothly}.\\

\nn {\bf Definition.} \label{Defpartunity} Let $\rho : X \to \R$ be a partition of unity w.r.t. $D$, i.e. $\rho$ is a smooth function and there are open set
$D \subset U_1 \subsetneq U_0$ with $\rho\equiv1$ on $U_1$ and $\rho \equiv 0$ on $X \setminus U_0$.\\

\nn {\bf Theorem.} \label{ThemXim}\emph{Let $m \neq 0$. For any choice of $\rho$ the modified Green functions 
\begin{align}
\label{mXim}
 \mXi(v,z,m):=  \Xi(v,z,m)  - \rho(z) \frac{t}{\sqrt{v}}  (\operatorname{B} (v,s,m)-\operatorname{I}(v,s,m) )\nne
\end{align}
are Green functions in the sense of Arakelov theory \`{a} la Gillet-Soul\'e for $T(m) \subset X$.
In particular, if  $m<0$, then $\mXi(v,z,m)$ is a smooth function on $X$.
} \\

{\bf Proof:} The function $\widetilde \Xi _\rho (m)$  is away from $T(m)$ a smooth function and along 
$T(m)$ it has logarithmic singularities. It remains to show that $dd^c \widetilde \Xi _\rho (m)$ is a 
smooth differential form. But this follows from the
 properties of the Schwartz form of the Kudla-Millson theory (where again we cut off bad terms at the boundary), 
which we recall with more details in the appendix, together with the smoothness
of the cut off function $\rho $ (see Corollary \ref{Corphirho}). \hfill $\Box$\\

\nn {\bf Remark.} We do not lift this construction to 
 a $\Gamma$-invariant function on $\HH^2$. If we want to lift $ \mXi(v,z,m)$ to a  function on $D$ smooth away from $T(m)$, then 
the cut off function $\rho $ has to be chosen more carefully. 

\newpage

{\bf Boundary function}\\

 We want to analyze the content of Theorem \ref{ThnearD}in the context of Funke-Millson [FM] .
Hence we propose the following notation.\\

\nn {\bf Definition.} We put
\begin{align}
\label{Xicheck}
\check \Xi (v,z,m) := (1/2)\sum_{-bc=m} \check \xi (v,z;b,c)\nne
\end{align}
with
\begin{align}
\label{boundary function}
\check \xi (v,z,;b,c) &:=  \,(t/\sqrt v)\,\, ( \operatorname{B}(v,s;b,c) - \operatorname{I}(v,s;b,c))  \nonumber\\
&= \begin{cases} (t/\sqrt v)\, \int_1^\infty e^{-\pi v(b/s+cs)^2r}r^{-3/2}dr - 4 \pi t\,\min (\mid bs^{-1} \mid ,\mid cs \mid )) \,\, {\rm if}\,\,-bc > 0\\
  (t/\sqrt v)\, \int_1^\infty e^{-\pi v(b/s+cs)^2r}r^{-3/2}dr  \,\,\quad \quad {\rm if}\,\,-bc \leq  0.
  \end{cases}\nne
\end{align}
and (for $m\not=0$) denote it as a  {\it boundary function.} The boundary function $\check \xi$ is also the result of a procedure going to the boundary as in Funke-Millson [FM] where they use a mixed model of the Weil representation. This also will be exploited in an appendix later.\\

As we shall need this later in our section 4, we add the following:\\

\nn {\bf Remark.}\label{RemXichek0} The definition(\ref{Xicheck}) also makes sense for $m=0$ 
even if we get an infinite series
\begin{align}\label{Xichek0}
\check{\Xi}(v,z,0) := (1/2)(t/\sqrt v)( \sum_{b\not=0}B(v,s;b,0) + \sum_{c\not=0}B(v,s;0,c) + B(v,s;0,0)).\nne  
\end{align}
{\bf Proof.}  This follows imediately from (\ref {Balpha}) and (\ref{7.1.13}).\hfill$\Box$\\

With the same reasoning, hence, one has parallel to the Corollary \ref{CorKUDLAGREEN}\\

\nn{\bf Corollary.}\label{Corgenbound} \emph{The generating series
\begin{align}\label{genbound}
\check{\Xi }(\tau ,z) := \sum_m \check{\Xi }(v,z,m) q^m \nne
\end{align}
converges uniformly on $\mathbb H^2$.}\hfill$\Box$ \\

As we find a nice formula, for further reference we analyze a bit more the singularities of $\Xi(v,z,0)$. \\

\nn {\bf Remark.} For (\ref{Xichek0}) we get
 \begin{align}
2\cdot \check{\Xi}(v,z,0) = &(t/\sqrt v)( \sum_{b\not=0}B(v,s;b,0) + \sum_{c\not=0}B(v,s;0,c) + B(v,s;0,0)).\nonumber\\
= & -2t/\sqrt v + t(s+1/s)(1/v + (2/\pi )\zeta (2)) \nonumber\\
&-(2t/\pi)((1/s)\sum _{b\in \mathbb N}e^{-\pi s^2b^2/v}/b^2 + s\sum_{c\in \mathbb N}e^{-\pi c^2/(s^2v)}/c^2 ).\nne 
\end{align}

 {\bf Proof.}  For $\vartheta (\tau ) = \sum_n e^{\pi i \tau n^2}$ using the standard theta transformation formula 
$$
\vartheta (i/v) = \sqrt v \vartheta (iv)
$$
 we get

\begin{align}
\label {two}
\sum_{b\in\mathbb Z} t/\sqrt v B(v,s;b,0)&=  t/\sqrt v  \int_1^\infty \sum_{b\in\mathbb Z}e^{-\pi vr(b/s)^2}dr/r^{3/2}\nonumber\\
&=  t/\sqrt v  \int_1^\infty \vartheta (is^{-2}vr))dr/r^{3/2}\nonumber\\
&=  ts/v  \int_1^\infty \vartheta (is^2/(vr)))dr/r^2\nonumber\\
&=  ts/v  \int_1^\infty \sum_{b\in\mathbb Z}e^{-(\pi s^2/(vr)))b^2}dr/r^2\nonumber\\
&=  ts/v + t/s \sum_{b\in\mathbb Z,b\not=0}(1/b^2)\int_0^{(bs)^2/v} e^{-\pi u}du\nonumber\\
&=  ts/v + (2t/(s\pi ))(\zeta (2) - \sum_{b\in \mathbb N}(1/b^2) e^{-\pi (bs)^2/v})\nne
\end{align}

and
\begin{align}
\label {three}
\sum_{b\in\mathbb Z,b\not=0}   t/\sqrt v B(v,s;b,0)&=  ts/v - 2t/\sqrt v + (2t/(s\pi ))(\zeta (2) - \sum_{b\in \N}(1/b^2) e^{-\pi (sb)^2/v}).\nne
\end{align}
\hfill$\Box$

 \newpage
\section{The zero-term}
 

The expression
\begin{align}\label{green-}
\xi (v,z,M) = \int_1^\infty e^{-2\pi vR(z,M)r} dr/r \nne
\end{align}
appearing for $m = \det M$ in the definition (\ref{Xi0}) of Kudla's Green function on $\mathbb D$ also works for matrices $M$ with $\det M\not=0$ if $M\not=0$. Although for $M=0$ the function $\xi(v,z,M)$ doesn't make sense at all, 
the $(1,1)$-form determined by $dd^c \xi(z,M)$ is well defined for $M=0$.
Actually we have (this will be shown in Section 7 in the appendix)

\begin{align*}
 dd^c \xi(v,z,\kzxz{0}{0}{0}{0}) := dd^c \int_1^\infty e^{-2\pi vR(z,M)r}dr/r \big |_{M=\kzxz{0}{0}{0}{0}}  = - \varphi _{KM}(v,z,\kzxz{0}{0}{0}{0}).
  \end{align*}
 We want to have a replacement for $\xi(z,\kzxz{0}{0}{0}{0})$ that satisfies the same differential equation. 
  Consider the usual modular form of weight 12
$$
\Delta (z) = q\prod _n(1-q^n)^{24}
$$
and set
\begin{align}
\label{XI00=}
\Xi ^0(v,z,0) &:= (1/2)\log \|\Delta (z_1)\Delta (z_2)\| ^{1/6}\nne
\end{align} 
where the Petersson norm for a $\Gamma -$modular form $F(z_1,z_2)$ of weight $k$ is given by 
$$
\| F(z_1,z_2) \| ^2 = |F(z_1,z_2)|^2(16\pi ^2y_1y_2)^k.
$$
Then it is well known (and again in the appendix (see(\ref{dd^cl})), that 
\begin{align*}
dd^c \Xi ^0(v,z,0) =  - (1/2)\varphi _{KM}(v,z,\kzxz{0}{0}{0}{0})
\end{align*}

  With $L_0^\ast = \{M = \left( \begin{smallmatrix} a & b\\ c & d\end{smallmatrix}\right), \det M = 0, M\not=0\}$ we set
\begin{align}
\Xi ^\ast (v,z,0) := (1/2)\sum_{M \in L_0^\ast } \xi (v,z,M).\nne  
\end{align}

\nn {\bf Definition.} \label{DefXivz0} We define the zero term of Kudla's Green function to be
\begin{align}\label{Xivz0}
\Xi (v,z,0) := \Xi ^\ast (v,z,0) + \Xi ^0(v,z,0).\nne
\end{align}

\nn {\bf Theorem.}  \label{ThnearD0}
\emph{Near the compactification divisor $D$ up to a smooth function zero along $D$ we have
\begin{align} \label{nearD0}
2\cdot \Xi (v,z,0) \sim  \frac{1}{6} \log|q_1 q_2| +
(t/\sqrt v) (\sum_{b\not=0}B (v,s;b,0) + \sum_{c\not=0}  B (v,s;0,c) + 2) \nne.
\end{align}
Using the definition (\ref{Xichek0}) this is}
\begin{align} \label{nearD00}
 \Xi (v,z,0) \sim& \,\, \frac{1}{12} \log|q_1 q_2| + \check{\Xi}(v,z,0)\nonumber\\
 \sim& \,-(\pi /6)t(s+1/s) + \check{\Xi}(v,z,0)\nne.
\end{align}

{\bf Proof.} {\bf Case 1. \,$M\not=0:$}\,\,
As in (\ref{Xi0d}) in the previous section,  we decompose
\begin{align}
\label{Xi00}
\Xi ^\ast (v,z,0) &= (1/2)\sum_{M \in L_0^\ast } \xi (v,z,M) = {\Xi ^\ast}^0  (v,z,0) + \sum_{d\in\mathbb N}{\Xi ^\ast }^d (v,z,0) \nne
\end{align}
and neglect the summands for which $d \not= 0$ as those vanish as Lemma \ref{lem:Xism} also holds for $m=0.$\\

For $d=0$ one has $b=0$ or $c=0$ or both. 
Take $c=0,$ hence
$$
R(z,M) = (1/2t^2)((a-bx_2)^2 + (by_2)^2).
$$
and similarly in the other cases. We write
$$
2{\Xi ^\ast }^0(z) = I + II + III
$$
where
$$
\begin{array}{rl}
I =& \sum_{b\in\mathbb Z,b\not=0}\sum_{a\in\mathbb Z}  \int_1^\infty e^{-\pi /t^2((a-bx_2)^2+(by_2)^2)vr}dr/r\\[.3cm]
II =& \sum_{c\in\mathbb Z,c\not=0}\sum_{a\in\mathbb Z}  \int_1^\infty e^{-\pi /t^2((a-cx_1)^2+(cy_1)^2)vr}dr/r\\[.3cm]
III =& \sum_{a\in\mathbb Z,a\not=0}  \int_1^\infty e^{-(\pi /t^2)a^2vr}dr/r\\[.3cm]
\end{array}
$$
As before in (\ref{hbc1}) we have
\begin{align}\label {zero}
I &= \sum_{b\in\mathbb Z,b\not=0}\sum_{a\in\mathbb Z}  \int_1^\infty e^{-\pi /t^2((a-bx_2)^2+(by_2)^2)vr}dr/r\nonumber\\
&= \sum_{b\in\mathbb Z,b\not=0} \sum_{n\in\mathbb Z}a_b(n)e(n{\bf x})\nne
\end{align}
with $bx_2 =:{\bf x},\, by_2 =: {\bf y} ,\, y_2^2/t^2=s^2$ where
\begin{align}\label {one}
a_b(0) &=  \int_{-\infty }^\infty \int_1^\infty e^{-(\pi /
^2)({\bf x}^2+{\bf y} ^2)vr}dr/r d{\bf x}\nonumber\\
&=  t/\sqrt v  \int_1^\infty e^{-\pi s^{-2}vrb^2}dr/r^{3/2}\nne\\
&=  t/\sqrt v B (v,s,(b,0))\nonumber.
\end{align}
 For $n\not=0$ we get as in Lemma \ref{lem:fourierexp} and Lemma \ref{lem:localexp}
\begin{align}\label {abn}
a_b(n) &=  \int_{-\infty }^\infty \int_1^\infty e^{-(\pi /t^2)({\bf x}^2+{\bf y} ^2)vr}e^{-2\pi in{\bf x}}dr/r d{\bf x}\nonumber\\
&=  2t/\sqrt v  \int_0^1 e^{-(\pi /t^2v{\bf y} ^2/{\bf x}^2 + \pi t^2n^2/v{\bf x}^2)}d{\bf x}\nne\\
&=  \frac{e^{-2\pi {\bf y} |n|}}{|n|} + \,\,\Phi \nonumber
\end{align}
where $\Phi $ indicates a harmless finite term vanishing along $D.$ Hence, again here with ${\bf x}=bx_2, {\bf y} = |by_2|,$ 
we get up to harmless terms in the sense of our Notation \ref{Nota} 
$$
\begin{array}{rl}
\sum_{n\not=0} a_b(n)e(n{\bf x}) \sim &  \sum_{n\in\N}(1/n)( e^{-2\pi ({\bf y} +i{\bf x})n} +  e^{-2\pi ({\bf y} -i{\bf x})n})\\
= & - \log(1 -  e^{-2i ({\bf y} +i{\bf x})n}) - \log(1 -  e^{-2i ({\bf y} -i{\bf x})n})\\
= &   \mathscr L_{b0}({\bf x},{\bf y})
\end{array}
$$
With $q_2 = e(z_2)$ for $b<0$ one has
$$
 \mathscr L_{b0}({\bf x},{\bf y}) = - \log |1-q_2^b|^2
$$
and for $b<0$
$$
 \mathscr L_{b0}({\bf x},{\bf y}) = - \log |1-q_2^{-b}|^2.
$$

For
\begin{align}\label {zeroc}
II &= \sum_{c\in\mathbb Z,c\not=0}\sum_{a\in\mathbb Z} - \int_1^\infty e^{-\pi /t^2((a-cx_1)^2+(cy_1)^2)vr}dr/r\nonumber\\
&= \sum_{c\in\mathbb Z,c\not=0} \sum_{n\in\mathbb Z}a_c(n)e(n{\bf x})\nne
\end{align}
with ${\bf x} = c{\bf x}_1, {\bf y} =|cy_1|$ we get in the same way up to a harmless $\Phi $
\begin{align}
\label {twoc}
\sum_{c\in\mathbb Z,c\not=0}  \sum_{n\in\mathbb Z}a_c(n)e(n{\bf x}) 
&= \sum_{c\in\mathbb Z,c\not=0} ( - \log|1-q_1^{|c|}|^2 + t/\sqrt v B (v,s;0,c) + \Phi \nne.
\end{align}
The third term is more tedious. Fortunately, there is a calculation using a regularization procedure recently 
done by Funke (see also \cite{Mu} p.86). It goes like this: With $\varrho \in\C$ we take
$$
\begin{array}{rl}
III_\varrho  :=& \sum_{a\in\mathbb Z,a\not=0}  \int_1^\infty e^{-(\pi /t^2)a^2vr}r^\varrho dr/r =: A - B\\[.3cm]
 A :=& \sum_{a\in\mathbb Z,a\not=0}  \int_0^\infty e^{-(\pi /t^2)a^2vr}r^\varrho dr/r \\[.3cm]
 B :=& \sum_{a\in\mathbb Z,a\not=0}  \int_0^1 e^{-(\pi /t^2)a^2vr}r^\varrho dr/r \\[.3cm]
\end{array}
$$
and, with $Z(2\varrho ) = \pi ^{-\varrho }\Gamma (\varrho )\sum_{a\in\N}a^{-2\varrho },$ get
$$
\begin{array}{rl}
A =& \sum_{a\in\mathbb Z,a\not=0}  \int_0^\infty e^{-(\pi /t^2)a^2vr}r^{\varrho-1}dr \\[.3cm]
=& \sum_{a\in\mathbb Z,a\not=0}  t^{2\varrho}(\pi va)^{-\varrho}\int_0^\infty e^{-r}r^{\varrho-1}dr \\[.3cm]
=& \sum_{a\in\mathbb Z,a\not=0}  t^{2\varrho}(\pi va^2)^{-\varrho}\,\Gamma (\varrho) \\[.3cm]
=& 2t^{2\varrho} v^{-\varrho}\,Z(2\varrho) \\[.3cm]
=& 2(1 + (\log t^2/v)\varrho + O(\varrho^2))(-1/(2\varrho) + {\rm const} \,+ O(\varrho)) \\[.3cm]
=& - 1/\varrho - (\log t^2/v) + 2\,{\rm const}\, + O(\varrho)) \\[.3cm]
\end{array}
$$
where const $= -(1/2)\log \pi - \log 2 + (1/2)\gamma , \gamma $ the Euler constant,  
and, using the theta transformation formula,
$$
\begin{array}{rl}
B =&  \int_0^1(\vartheta (ivr /t^2) - 1)r^\varrho dr/r\\[.3cm]
=& \int_1^\infty (\vartheta (it^2r/v)t\sqrt {r/v} - 1)r^{-\varrho}dr/r\\[.3cm]
=& - \int_1^\infty r^{-\varrho-1}dr + (t/\sqrt v)\int_1^\infty \vartheta (it^2r/v))r^{-\varrho-1/2}dr\\[.3cm]
=& - \int_1^\infty r^{-\varrho-1}dr + (t/\sqrt v)\int_1^\infty r^{-\varrho-1/2}dr 
+ 2(t/\sqrt v)\sum_{a\in\mathbb N}\int_1^\infty  e^{-(\pi t^2)a^2r/v}r^{-\varrho-1/2}dr\\[.3cm]
=& - 1/\varrho - (t/ \sqrt v)(1/((1/2)-\varrho) + \Phi (t,\varrho)
\end{array}
$$
where $\varrho >1/2$ and $\Phi $ is something harmless: For $t\rightarrow \infty $ one may estimate
\begin{align*}
2(t/\sqrt v)\sum_{a\in\mathbb N}\int_1^\infty  e^{-(\pi t^2)a^2r/v}r^{-\varrho -1/2}dr 
&= 2(t/\sqrt v)\int_1^\infty  (1/2)(\vartheta (it^2r/v)-1)r^{-\varrho -1/2}dr\\
&\leq 2(t/\sqrt v)\int_1^\infty  (1/2)Ce^{-t^2r/v}r^{-\varrho-1/2}dr.
\end{align*} 
We add $A$ and $-B,$ take $\varrho  \rightarrow  0$ and get
\begin{align}
\label{III}
 III =  \lim _{\varrho \rightarrow 0}III_\varrho  =   2(t/\sqrt v) - \log (t^2/v) + \Phi \nne.
\end{align}
Finally, we add the contributions from the three terms $I,II,III$ where we only keep the contributions which don't 
go to zero as $t\longrightarrow \infty $ and obtain
\begin{align}
2{\Xi^\ast }^0(v,z,0) \sim  t/\sqrt v (\sum_{b\not=0}B (v,s;b,0) + \sum_{c\not=0}  B (v,s;0,c) + 2) - \log (t^2/v)\nne.
\end{align}

{\bf Case 2.\, $M=0:$}\,\,It is straightforward to see
\begin{align}
2\cdot {\Xi }^0(v,z,0)  
&=  \log t^2 - (1/3)\pi (y_1+y_2) + \log (16\pi ^2) + \sum _{m,n\in \N} \log |1-q_1^m|^4\log |1-q_2^n|^4 \nonumber \\
&\sim \log t^2 - (1/3)\pi (ts + t/s)\nne. 
\end{align}

 \hfill $\Box$
 
Recall our choice $T(0)=T(v,0) =(\frac{-1}{24} + \frac{1}{8 \pi v}) D$. 
The logarithmic singularities  of $\Xi(v,z,0)$ along $D$   have a different mutltiplicity than those of $T(0)$,  therefore the zero term $\Xi(v,z,0)$ does not extend to a Green function 
for $T(0)$ on $X$. In addition there are unwanted singularities along $D$, which we will substract as in section 2. \\

Again, let $\rho : X \to \R$ be a partition of unity w.r.t. $D$, i.e. $\rho$ is a smooth function and there are open set
$D \subset U_1 \subsetneq U_0$ with $\rho \equiv1$ on $U_1$ and $\rho  \equiv 0$ on $H \setminus U_0$.\\

\nn {\bf Theorem.}\label{ThemXi00} \emph{For any choice of $\rho$ the modified Green functions 
\begin{align}\label{mXi00}
 \mXi (v,z,0):=&\, \Xi (v,z,0)  - \rho(z) (\check{\Xi}(v,z,0) - \frac{t}{2v} \left(s + \frac{1}{s} \right) )\nne
\end{align}
are Green functions for $T(0) \subset X$ in the sense of Arakelov theory \`a la Gilet-Soul\'e. In particular one has} 
\begin{align}\label{mXi000}
 2\mXi (v,z,0) \sim & \,(-\pi /3 + 1/v)t(s+1/s) = (1/6 - (1/(2\pi v)))\log |q_1q_2| \nne
\end{align}

{\bf Proof.} Obvious with the previous expansion plus the identity
$$-\log|q_1 q_2| =  2 \pi t  \left(s + \frac{1}{s} \right)$$
\hfill $\Box$
 
 \nn {\bf Remark.} The additional summand $(1/2v)t(s+1/s)$ in (\ref{mXi00}) will reappear later and find another 
 explication in providing the modularity of the generating function belonging to the $\check{\Xi}(v,z,m)$\\

\newpage
\section{Generating series of differential forms and boundary functions}

Our choice for the modification/regularisation of $\mXi(v,z,m)$ turns out to be ``modular'',  in order to (partially) describe what we mean by that we set 
\begin{align*}
\varphi_\rho (v,z,m):= - dd^c \mXi (v,z,m).
\end{align*}


\nn {\bf Theorem.} \emph{For any choice of $\rho$ the generating function 
\begin{align*}
\sum_{m \in \mathbb Z} \varphi_\rho(v,z,m) q^m  
\end{align*} 
is a weight 2 modular form for the group $\SL_2(\mathbb Z)$ with values in} $ A^{1,1}(X).$\\

{\bf Proof.} We first observe the decomposition 
\begin{align*}
\varphi_\rho(v,z,m) =& - dd^c\mXi(v,z,m) = \varphi_{KM}(v,z,m) -  \varphi_{FM,\rho}(v,z,m)
\end{align*}
Here the first form is the Kudla-Millson form and a variant of the second form
\begin{align*}
\varphi_{FM,\rho}(v,z,m) :=& - dd^c(\rho (z)\check{\Xi }(v,z,m)) 
\end{align*}
has been considered by Funke and Millson (all this will be discussed in our appendix).
The claim follows now immediately from Proposition \ref{prop:KudlaMillson} 
and by applying the $dd^c$-Operator to the modular form of Proposition \ref{prop:Zagier-Funke} below.

\hfill $\Box$

\nn {\bf Proposition.} (Kudla-Millson $+\epsilon$)\label{prop:KudlaMillson}  
\emph{The generating series $\sum_m \varphi_{KM}(m) q^m  $ of the
Kudla-Millson forms is a modular form of weight $2$ for the group $\SL_2(\ZZ)$ with values in $A^{1,1}(X)$.}\\

{\bf Proof.} This is more or less a tautology, as it coincides with the theta series 
\begin{align}
\label{theta2}
\Theta _{\varphi_{KM}}(\tau ,z ) = \sum_{M\in L} \varphi_{KM}(v,z,M)q^{\det M}\nne 
\end{align} 
whose modularity for some subgroup of $\SL_2(\ZZ)$ is part of the general Kudla-Millson theory (\cite{KM1}).
In section \ref{Toy} of our appendix we shall give a direct proof leading directly to the group $\SL_2(\ZZ)$ 
independently of the work of Kudla Millson and based 
on Siegel's treatment of thetas for indefinite quadratic forms.
\hfill $\Box$\\

Already in Theorem \ref{ThemXi00} we had reason to modify the zero-term $\check{\Xi}(v,s,0).$ Hence we put  
\begin{align}
\label{XiW+}
\check \Xi ^+(v,z,0) := \check \Xi (v,z,0) - (1/2v)t(s+1/s)\nne
\end{align}
and
\begin{align}
\label{XiW++}
\check \Xi ^+(\tau ,z) := \check \Xi (\tau ,z) - (1/2v)t(s+1/s)\nne .
\end{align}
With this modification, we get a result which is based on the seminal Hirzebruch-Zagier paper \cite{HZ} and promoted by discussions with Funke on a 
forthcoming Funke-Millson paper \cite{FM3}.

\nn {\bf Proposition.} (Zagier, Funke$+\epsilon$) \label{prop:Zagier-Funke}\emph{For all $z\in \mathbb H^2,$ the generating series 
$$ 
\check \Xi ^+(\tau ,z) = \sum_{m\in \Z} \check \Xi (v,z,m)q^m- (1/2v)t(s+1/s)  \qquad (\tau=u+iv,\,q = exp(2 \pi i \tau))
$$
is a  non-holomorphic modular form of weight $2$ for $SL_2(\Z)$.}\\

\nn{\bf Remark.} In Hirzebruch-Zagier \cite{HZ} p.98 one can find the following result:\\

Take
 \begin{align*} 
{\rm erfc}(x) &:= 2/\sqrt \pi \int_x^\infty e^{-u^2}du\\[.3cm]
\beta (x) &:= 1/(16\pi )\int_1^\infty e^{-xu}u^{-3/2}du
= 1/(8\pi ) (e^{-x} - \sqrt{x\pi }\,{\rm erfc}\,(\sqrt x))
\end{align*}
and define for $\tau = u + iv \in \mathbb H,\, \lambda ,\lambda ' \in \R$

\begin{align*} 
U_\tau (\lambda ,\lambda ') &:= 2v^{-1/2} \beta (\pi v(\lambda -\lambda ')^2)e(\lambda \lambda '\tau )\\[.3cm]
V_\tau (\lambda ,\lambda ') &:= \begin{cases} (1/2) {\rm min}\,(\mid \lambda \mid ,\mid \lambda ' \mid )e(\lambda \lambda '\tau ) \,\,\,
&{\rm if}\,\lambda \lambda ' > 0,\\
 \,\,0 \,\,\quad &{\rm if}\,\lambda \lambda ' \leq  0,
 \end{cases}\\
W_\tau (\lambda ,\lambda ') 
&:= U_\tau (\lambda ,\lambda ') - V_\tau (\lambda ,\lambda ').
\end{align*} 
Then the Fourier transform
\begin{align*}
\hat W_\tau (\mu ,\mu ') &:= \int_{-\infty }^\infty  \int_{-\infty }^\infty W_\tau (\lambda ,\lambda ') 
e(-\lambda \mu - \lambda '\mu ') d\lambda d\lambda ',
 \end{align*}
satisfies the relation
\begin{align}
\label{W}
\hat W_\tau (\mu ,\mu ') = \tau ^{-2} W_{-1/\tau } (\mu ,\mu ').\nne
\end{align}
If $\ot$ is the maximal order in the real quadratic field $K$ with discriminant $D,$ 
then, as a consequence of \eqref{W}, the authors deduce  that 
\begin{align}
\label{Wseries}
\mathcal W (\tau ) := \sum _{\lambda \in \ot} W_\tau (\lambda ,\lambda ')\nne
\end{align}
is a non-analytic modular form of weight 2, level $D$ and Nebentypus $\chi _D$. 
Similar considerations are also found 
in \cite{FM2} Prop. 4.20) with $m=2$ and $\ell =2$, i.e., 
\begin{align}
\label{xiW}
\check \xi (\tau ,(t,s),(b,c))   
&= 8\pi t\,W_\tau (-b/s,cs). \nne
\end{align}
\\

{\bf Proof of Proposition \ref{prop:Zagier-Funke}.}    In our case \eqref{Wseries} is just given by    
\begin{align}
\label{XiW}
\check \Xi (\tau ,z) &= \sum_m \check \Xi (v,z,m)q^m\nonumber\\[.3cm]
&= \sum_{b,c\in \Z} (t/\sqrt v)(B(v\alpha (s)) - I(v\alpha (s))))e^{2\pi i(-bc)\,\tau }\nonumber\\[.3cm]
&= \sum_{b,c\in \Z} 8\pi t \,(U_\tau (-b/s,cs) - V_\tau (-b/s,cs))\nonumber\\
&= \sum_{b,c\in \Z}8\pi t\,W_\tau (-b/s,cs). \nne
\end{align}
Here we used $\alpha (s) = \pi (b/s+cs)^2$. Recall
convergence is known from Corollary \ref{Corgenbound}.
We want to use the equation (\ref{W}) to prove the modularity of  $\check \Xi $, but we will run into problems
as we will see that convergence of $\check \Xi (\tau ,z)$ is not strong enough to satisfy the usual Poisson summation formula
which is the key in proving modularity.
In fact we prove in Lemma \ref{Wlemma} that
we have a Poisson summation in the form

\begin{align}
\label{Poisson}
\sum_{b,c\in\mathbb Z}\hat W_\tau (-b/s,cs) = \sum_{b,c\in\mathbb Z} W_\tau (-b/s,cs) -( i/4 \pi  \tau )(s+1/s).\nne
\end{align}
Now  using the notation  \eqref{XiW++} 
$$
2\cdot \check{\Xi}^+(\tau,z) := 8\pi t \sum_{b,c\in\mathbb Z} W_\tau (-b/s,cs) - (1/v)t(s+1/s))
$$
and the above cited result (\ref{W}) on the Fourier transform shows that one has
$$
\begin{array}{rl}
2\cdot \check{\Xi}^+(-1/\tau,z) =& 8\pi t\sum_{b,c\in \mathbb Z} W_{-1/\tau }(-b/s,cs) - (|\tau|^2/v)t(s+1/s))\\[.3cm]
=& 8\pi t\tau ^2\sum_{b,c} \hat W_\tau (-b/s,cs)- (|\tau|^2/v)t(s+1/s)\\[.3cm]
=& 8\pi t\tau ^2\sum_{b,c\in\mathbb Z} W_\tau (-b,c) - 2it\tau (s+1/s)- (|\tau|^2/v)t(s+1/s))\\[.3cm]
=& \tau ^2(8\pi t\sum_{b,c\in\mathbb Z} W_\tau (-b,c) - (1/v)t(s+1/s))\\[.3cm]
=& \tau ^2 2\cdot \check{\Xi}^+(\tau,z) 
\end{array}
$$

and, finally, the claim follows. \hfill $\Box$

\nn {\bf Lemma.} \label{Wlemma} \emph{
We have a Poisson summation in the form
\begin{align*}
 \sum_{b,c\in\mathbb Z}\hat W_\tau (-b/s,cs) = \sum_{b,c\in\mathbb Z} W_\tau (-b/s,cs) 
 -( i/4 \pi  \tau )(s+1/s)
\end{align*}}

{\bf Proof.} With $B (\alpha ) := \int_1^\infty e^{-\alpha r} r^{-3/2}dr$   as in  \eqref{B}) and for fixed $s>0$, we slightly change notation and put
\begin{align*}
f_\tau (x,y) &:= W_\tau (x/s,ys) \\
&=\begin{cases} ((1/(8\pi \sqrt v))B(\pi v(x/s-ys)^2) - (1/2)\min\,(|x/s|,|ys|)) e^{2\pi i\tau xy}\,\, &{\rm if}\,\, xy > 0\\
 ((1/(8\pi \sqrt v))B(\pi v(x/s-ys)^2)e^{2\pi i\tau xy}\,\,\,\,&{\rm if}\,\, xy \leq  0.
 \end{cases}
\end{align*}

The function
\begin{align}
\label{Ftau}
F_\tau (x,y) := \sum_{m,n\in\Z} f_\tau (x+m,y+n)\nne
\end{align}
is periodic with period 1 in $x,y$ and therefore has a Fourier series
$$
\mathcal F _\tau (x,y) = \sum_{\mu ,\nu \in \Z } c(\mu ,\nu )e(\mu x+\nu y)
$$
with
$$
\begin{array}{rl}
c(\mu ,\nu ) =& \int_0^1\int_0^1 F_\tau (x,y))e(-\mu x-\nu y)\\[.3cm]
=& \hat f_\tau (\mu ,\nu ).
\end{array}
$$
 To get as usual a Poisson summation formula, one has to compare $F_\tau $ with its Fourier series $\mathcal F_\tau $ 
near $(0,0)$. As kindly pointed out by Funke, though $f_\tau $ is continuous, the associated periodic function  $F_\tau $ need not to be so also. We claim
\begin{align}
\label{claimF}
 F_\tau(x,y) \, \textit{is  discontinuous at } xy=0\textit{ in a small neighbourhood of } (0,0) .\nne
\end{align}

This leaves us with more a technical problem, since one may use a theorem  from Moore \cite{Mo} that generalizes the 
fact that the Fourier series at a point of non-continuity equals the arithmetic means of the limits from right and left from one to two dimensions. 
This theorem tells us that the value of the Fourier expansion of $F$ 
at $(0,0)$ 
is summable to one fourth 
of the sum of the four values of $F$ which come up as limits when approaching $(0,0)$ in the four quadrants.\\

We now want to prove \eqref{claimF}and determine those limits.\\

 At first we put
$$
f_1(x,y) := B(\pi v(x/s-ys)^2)e^{2\pi i\tau xy}
$$
and
$$
F_1(x,y) :=  \sum_{m,n\in\Z} f_1 (x+m,y+n).
$$
For $\alpha \geq 0$ one has
$$
B(\alpha)= \int_1^\infty e^{-\alpha u}u^{-3/2}du \leq e^{-\alpha }\int_1^\infty u^{-3/2}du = 2 e^{-\alpha }
$$
and
$$
e^{-\pi v(x/s-ys)^2} \mid e^{-2\pi i xy\tau } \mid \,= e^{-\pi v((x/s)^2+(ys)^2)}.
$$
Hence, for small $|x|,|y|$,  $F_1(x,y)$ is majorized by
$$
\sum_{m,n} e^{- \pi v(((x+m)/s)^2+((y+n)s)^2)}.
$$
This Majorant converges uniformly and hence $F_1(x,y)$ is
continuous in $(0,0)$ (by Weierstrass' Majorant Theorem). \\
 
Thus, we are left with the discussion of
\begin{align}\label{deff2}\nne
f_2 (x,y) &:= \begin{cases}\min\,(|x/s|,|ys|) e^{2\pi i\tau xy}\,\, &{\rm if}\,\, xy > 0\\
\,\,0\,\,&{\rm if}\,\, xy \leq  0
\end{cases}
\end{align}
and we need to study 
$$
F_2 (x,y) := \sum_{m,n \in\Z}\,f_2 (x+m,y+n)
$$
where throughout $s>0$ will be hold fixed if not specified otherwise.
If we show that  
\begin{align}\label{claimF2}
F_2 \,\textit{ is discontinuous at } xy =0 \textit{ in a small neighbourhood of } (0,0). \nne
\end{align}
we have proven \eqref{claimF}.   
We have
\begin{align}\label{claimF2sum}
F_2(x,y) = \sum_{m,n\in\Z,m\not=0,n\not=0}f_2(x+m,y+n)  + \sum_{m\in\Z,m\not=0}f_2(x+m,y) + \sum_{n\in\Z,n\not=0}f_2(x,y+n) + f_2(x,y).\nne
\end{align}
We put
\begin{align}\label{defgtau}
g_\tau (x,y) := \sum_{n \in\Z,n\not=0}\,f_2 (x,y+n)\nne
\end{align}
and
\begin{align}\label{defhtau}
h_\tau (x,y) := \sum_{m \in\Z,m\not=0}\,f_2 (x+m,y)\nne
\end{align}
and, hence, have
\begin{align}\label{claimF2sumgh}
F_2(x,y) = \sum_{m\in\Z,m\not=0}g_\tau (x+m,y)  + g_\tau (x,y) + h_\tau (x,y) + f_2(x,y).\nne
\end{align}
From \eqref {deff2} we know $f_2(0,0)= g_\tau (0,y) = h_\tau (x,0) = 0$ and therefore
\begin{align}
F_2(0,0) = \sum_{m,n\in\Z,m\not=0,n\not=0}f_2(m,n) = \sum_{m\in\Z,m\not=0}g_\tau (m,0).\nne
\end{align}
We want to compare this value with the limits of $F_2(x,y)$ if we approach $(0,0)$  in the 
four quadrants $Q_j$ with $x,y>0, -x,y>0, -x,-y>0$ and $x,-y>0$ respectively. Hence, 
in the sequel, all the time, we have  $-1 \ll x,y \ll 1$. \\

{\bf 1.} At first we look at $g_\tau :$ For $x>0$ and $n\in\mathbb N$ one has $f_2 (x,y-n) = 0$. 
Hence with $q = e^{2\pi i\tau }$, for these $x,y,$ we get

\begin{align*}
g_\tau (x,y) =& \sum_{n=1}^\infty \,f_2 (x,y+n)\\
=& \sum_{n=1}^\infty \,(x/s)e^{2\pi i\tau  x(y+n)}\\
=& (x/s)q^{xy+x}\,\sum_{n=0}^\infty \,q^{xn}\\ 
=& (x/s)q^{xy+x}\,\,/(1-q^x ).  
\end{align*}

Now via $q^x = 1 + 2\pi i\tau x+\dots $ we get as limiting values 
in the quadrants $Q_1$ and $Q_4$
\begin{align} 
\lim_{x\searrow 0, y \neq 0} g_\tau (x,y) &= -1/(2\pi is\tau )=: \beta _1/s\nne
\label{gtau}
\end{align}
Similarly, for $x<0$ and $n\in\mathbb N$ one has $f_2 (x,y+n) = 0$ and
\begin{align*}
g_\tau (x,y) =& \sum_{n=1}^\infty \,f_2 (x,y-n)\\
=& \sum_{n=1}^\infty \,(-x/s)e^{2\pi i\tau  x(y-n)}\\
=& (-x/s)q^{xy-x}\,\sum_{n=0}^\infty \,q^{-xn}\\ 
=& (-x/s)q^{xy-x}\,\,/(1-q^{-x} ),  
\end{align*}
i.e., in the quadrants $Q_2$ and $Q_3$ again we get as limiting values 
\begin{align} 
\lim_{-x\searrow 0, y \neq 0} g_\tau (x,y) &= -1/(2\pi is\tau )=: \beta _1/s\nne
\label{gtau-}
\end{align}

{\bf 2.} In the same way, for $y>0,$ one has
\begin{align*}
h_\tau (x,y) =& \sum_{n=1}^\infty \,f_2 (x+m,y)\\ 
=& \sum_{n=1}^\infty \,yse^{2\pi i\tau  (x+m)y}\\
=& ysq^{xy+y}\,\,/(1-q^x )
\end{align*}
and, hence in the quadrants $Q_1$ and $Q_2$
\begin{align}
\lim_{ y\searrow 0, x \neq 0}  h_\tau (x,y)= -s/(2\pi i\tau )= \beta _1s .\nne\label{htau}
\end{align}
and, with a similar reasoning as above, also in the quadrants $Q_3$ and $Q_4.$\\

{\bf 3.} The results from {\bf 1.} and {\bf 2.} together tell us that in each quadrant $h_\tau + g_\tau $ in the limit to $(0,0)$ go to
\begin{align*}
\beta _1(s+1/s).
\end{align*}

We are left to study 

\begin{align*}
\sum_{m\not=0}\,g_\tau (x+m,y) = \sum_{m,n\in\Z,m\not=0,n\not=0}f_2(x+m,y+n)
\end{align*}
By looking long enough the skillful reader will probably see that this sum is well behaved in the neighbourhood of $(0,0),$ i.e., 
continuous with $ \sum_{m\not=0}\,g_\tau (m,0) = F_2(0,0).$ If one accepts this, one is done: we have in each quadrant the limit value
$ \lim F_2(x,y) = F_2(0,0) - 1/(2\pi i\tau )(s+1/s)$ and 
by Moore's Theorem we get for the Fourier expansion $\mathcal F_2$ of $F_2$
\begin{align*}
\mathcal F_2(0,0) = F_2(0,0) - 1/(2\pi i\tau )(s+1/s).
\end{align*}
As a consequence of all this, going back to the notation (\ref{Ftau}), in our case ( $-(1/2)f_2$ entering in $f_\tau $), 
we have a Poisson summation in the form claimed above. \\

{\bf 4.} If one does not see the continuity of $\sum_{m\not=0}\,g_\tau (x+m,y),$ for $m>0,$ let us look at
$g_\tau (x+m,y)$. Again we remark that for small $x,y$ and $n>0$ we have $(x+m)(y-n)<0$ and 

\begin{align*}
\min((x+m)/s,(y+n)s) &=\begin{cases} (x+m)/s \,\,&{\rm if}\,\,m/s\leq ns \\
 (y+n)s \,\,&{\rm if}\,\,ns<m/s.
 \end{cases}
\end{align*}

In the hope to convince the reader without too much trouble, we only consider two special cases, namely $s=1$ and $s=1/2.$
For $s=1,$ with $\epsilon _m:= e^{2\pi i\tau (x+m)}$ we have
\begin{align*}
g_\tau (x+m,y) =& \sum_{n=1}^{m-1}(y+n)e^{2\pi i\tau (x+m)(y+n)} + (x+m)\sum_{n=m}^\infty e^{2\pi i\tau (x+m)(y+n)}\\[.3cm]
=& \epsilon _m^y ( \sum_{n=1}^{m-1}(y+n)\epsilon _m^n + (x+m)\epsilon _m^m\sum_{n=0}^\infty \epsilon _m^n)\\[.3cm]
=& \epsilon _m^y ( y\sum_{n=1}^{m-1}\epsilon _m^n + \sum_{n=1}^{m-1}n\epsilon _m^n + x\epsilon _m^m/(1-\epsilon _m)+ m\epsilon _m^m/(1-\epsilon _m)) \\ 
=& \epsilon _m^y ( y\epsilon _m\frac{1-\epsilon _m^{m-2}}{1-\epsilon _m} + \frac{\epsilon _m-m\epsilon _m^m+(m-1)\epsilon _m^{m+1}}{(1-\epsilon _m)^2} + \frac{x\epsilon _m^m}{(1-\epsilon _m)}+\frac{ m\epsilon _m^m}{(1-\epsilon _m)})
\end{align*} 

  Here one observes continuity in $x,y$ and 
 with $\epsilon = e^{2\pi i \tau m}=q^m$ for $x=y=0$
we are left with
\begin{align*}
g_\tau (m,0) =& \sum_{n=1}^{m-1}n\epsilon ^n + m\epsilon ^m/(1-\epsilon )\\[.1cm]
=& \epsilon (1-\epsilon ^{m}) /(\epsilon  - 1)^2\\
=& q^m\frac{1+q^m}{1-q^m}.
\end{align*} 
By summation $\sum_{m>0}\,g_\tau (x+m,y),$ again one has continuity and
\begin{align*}
\sum_{m\in\mathbb N}g_\tau (m,0) = \sum_{m\in\mathbb N}q^m\frac{1+q^m}{1-q^m}\\
\end{align*} 
and, as, by symmetry, for $m<0$ we get the same value, we end up with the value for $s=1$
\begin{align}\label{F200}
F_2(0,0) = \sum_{m\in\mathbb Z, m\not=0}g_\tau (m,0) = 2\sum_{m\in\mathbb N}q^m\frac{1+q^m}{1-q^m}.\nne
\end{align}  
For $s=1/2$ and $m>0$ one has
\begin{align*}
\min ((x+m)2,(y+n)/2) &=
\begin{cases} (x+m)2 \,\,&{\rm if} \,\,4(x+m) \leq (y+n)\\
 (y+n)/2 \,\,&{\rm if} \,\,4(x+m) > (y+n)
\end{cases}
\end{align*}}
and
$$
\begin{array}{rl}
g_\tau (x+m,y) =& \sum_{n=1}^{4m-1}(y+n)/2e^{2\pi i\tau (x+m)(y+n)} + 2(x+m)\sum_{n=4m}^\infty e^{2\pi i\tau (x+m)(y+n)}\\[.3cm]
=& \epsilon _m^y ( y/2\sum_{n=1}^{4m-1}\epsilon _m^n + \sum_{n=1}^{4m-1}n/2\epsilon _m^{n} + 2x\epsilon _m^{4m}/(1-\epsilon _m)+ 2m\epsilon _m^{4m}/(1-\epsilon _m)) \\[.3cm]
\end{array} 
$$
where
$$
\begin{array}{rl}
g_\tau (m,0) =& (1/2)\epsilon (1-\epsilon ^{4m}) /(\epsilon  - 1)^2.\\
\end{array} 
$$
Hence, as above, we have continuous participants as, again with similar considerations, in the remaining cases. 

\hfill $\Box$


Another application of the modularity result in Proposition  \ref{prop:Zagier-Funke} beside the 
modularity of $dd^c (\rho \,\check{\Xi}^+(v,z)$  
is situated in the context of the Hirzebruch Zagier Theorem.

\nn {\bf Proposition.} \label{PropthetadT}\emph{The following identity of generating series
\begin{align}
\label{thetadT}
\int_{T(1)} \sum_{m \in \ZZ} dd^c (\rho \,\check{\Xi}^+(v,z,m))q^m &=
\lim_{\varepsilon \rightarrow 0}\int_{\partial T_\varepsilon (1) } d^c \check{\Xi }(\tau ,z)  \nonumber \\
 &=  -1/(4\pi v) - (1/2)\sum_{b,c \in \mathbb Z , \, -bc>0} \min (|b|,|c|)q^{-bc} \nonumber\\
&\quad +\sum_{b,c \in \mathbb Z}(1/8\pi )v^{-1/2} \int_1^\infty  e^{-\pi v(b+c)^2r}r^{-3/2}dr)q^{-bc}\nne
\end{align}
is an identity of  modular forms of weight $2$ for $\SL_2(\ZZ)$ }
(c.f.\cite{FM2} Th.6.6.).\\

{\bf Proof. 1.} If $U_\varepsilon (\infty )$ is a small neighbourhood of $\infty \in T(1) ,$ and $T_\varepsilon = T\backslash U_\varepsilon $ 
one has the existence of a in $\varepsilon $ continuous 
function
$$
F(\varepsilon ) := \int _{T_\varepsilon (1)} dd^c \mXi (\tau ,z)
$$ 
because $\tilde \Xi (\tau ,z)$ is smooth on $T(1) \setminus U_\varepsilon $ (up to log-singularities outside of $U_\varepsilon $) 
and with Stokes' theorem we have
$$
F(\varepsilon ) = \int_{\partial T_\varepsilon (1)} d^c \mXi(\tau ,z).
$$ 
From the previous proposition we know that $\int_{T(1)} dd^c \Xi$ exists and hence one also has continuos in small $\varepsilon $ 
$$
f(\varepsilon ) =  \int _{T_\varepsilon (1) } dd^c (\rho \check{\Xi}) 
= \int_{\partial T_\varepsilon (1)} d^c (\rho \check{\Xi}).
$$
If we let $\varepsilon \rightarrow 0,$ one gets
$$
f(0) =  \int _{T(1)} dd^c (\rho \check{\Xi}) = \lim_{\varepsilon \rightarrow 0}\int_{\partial T_\varepsilon (1)} d^c \check{\Xi}.
$$ 
{\bf 2.} To evaluate the line integral we remind of (\ref{XiW++})
$$
\check{\Xi}^+(\tau ,z) = (1/2)\sum_{b,c} (t/\sqrt v)(B(v,z;b,c)-I(v,z;b,c))q^{-bc} - (t/v)(s+1/s) 
$$ 
and from the derivation formalism ((\ref{d'tFs}) in the appendix), we know

\begin{align*}
 d^c (t(B - I)) =& (1/4\pi)(t(B-I)\omega _{13} - st(B'-I') \omega _{24}\\
 =&  (1/8\pi)(-(B-I)(dx_1/s+sdx_2) - (B'-I') (dx_1-s^2dx_2).\\
\end{align*}

and 
\begin{align*}
d^c(t(s+1/s) =& (1/4\pi)(t(s+1/s)\omega _{13} - st(1-s^{-2}) \omega _{24}\\
=&- (1/8\pi)(s+1/s)(dx_1/s+sdx_2) - (1-1/s^2) (dx_1-s^2dx_2).
\end{align*}

The curve $T(1) $ is given on $\mathfrak H \times \mathfrak H$ by $z_1 = z_2 =: z.$ Hence we have $s=1,$ 
in both terms the second summand vanishes and the integral 
of the one-form $ d^c (t(B - I)) = d^c(tF(s))$
over the boundary $\partial T$ realized on the quotient $\Gamma \setminus \ D$ 
is given by
\begin{align}
\label{inttF}
\lim_{\varepsilon \rightarrow 0}\int_{\partial T_\varepsilon (1) }d^c (tF(s))) =  (1/(4\pi ))\int_0^1 (B - I) dx. \nne
\end{align}
and, moreover, one has
$$
\lim_{\varepsilon \rightarrow 0}\int_{\partial T_\varepsilon (1) }d^c t =  (1/4\pi ) \int_0^1 (1/4\pi )F(1)dx
$$
for $t \rightarrow \infty .$ (here we had to pay attention that the orientation coming from 
Stokes' theorem makes $\partial T_\varepsilon (1)$ the negative of $\int_0^1$). Since $x$ does not appear in the integrand and one has $s=1$, we 
have a trivial evaluation and get (again parallel to the results from \cite{FM2}) the result in the proposition.\hfill $\Box$\\

As corollary we recover in our special case a result of Funke's thesis.
  
\nn {\bf Corollary.}\label{CorthetaTt} \emph{The following identity of generating series
\begin{align}
\int_{T(1)} \sum_{m \in \ZZ} dd^c {\Xi}(v,z,m))q^m &=
- \int_{T(1)} 
\sum_{m \in \ZZ} 
 \varphi_{KM}(v,z,m)q^m  \nonumber \\
&=  \sum _{N \geq 0} H_1(N) q^N + 
\sum _{b,c\in\mathbb Z}(1/(8\pi \sqrt v)) \int_1^\infty  e^{-\pi v(b+c)^2r}r^{-3/2}dr)q^{-bc}.\nne
\label{thetaTt}
\end{align}
is an identity of modular forms of weight $2$ for $\SL_2(\ZZ)$.}\\ 
 
{\bf Proof.}  
Since $\mXi(v,z,m)$ is a Green function for $T(m)$ the the cohomology class of the form
$\varphi_\rho(v,z,m)$ satisfies 
\begin{align*}
[dd^c\mXi ] = [-\varphi_\rho(v,z,m)] = [ T(m) ] \in H^2(X,\mathbb R)
\end{align*}
actually for $m\neq 0$ this is an equality in $H^2(X,\mathbb Z)$. 

As a direct consequence  to theorem \ref{thm:hurwitz-hirzebruch} we derive

\begin{align}\label{Eisenstein}
\int_{T(1)} \sum_{m \in \ZZ} 
 - \varphi_\rho(v,z,m)q^m  &= 
[T(1)] \sum_{m \in \ZZ} [T(m)] q^m  = [T(1)][T(1)]  \cdot E_2(\tau,1)\nonumber\\
&= 2 \cdot E_2(\tau ,1)\nne
\end{align}

Finally, we recall from \cite{Hi}, p.82 (10), the nice identities for the Hurwitz class numbers for positive $N$
\begin{align}
\label{Hurwitz}
 2\sigma _1(N) =
H_1(N) + \min_{bc=N}( b, c)   = H_1(N) + 2 \sum_{d|N,d\leq \sqrt N} d   \nne
\end{align}
Splitting the integral over 
$\varphi_\rho(v,z,m)$
with respect to the decomposition  
\begin{align}
\label{hz-decomposition}
-\varphi_\rho(v,z,m) =&  dd^c\mXi(v,z,m) = - \varphi_{KM}(v,z,m) + dd^c(\rho (z)\check{\Xi }^+(v,z,m))\nne
\end{align}
the claimed identity is now straigtforward to see. \hfill $\Box$

\nn {\bf Remark.} \label{RemthetaTt} On the other hand, if one 
reads \eqref{hz-decomposition} in the other direction, e.g. 
if one proves corollary \ref{CorthetaTt} directly, 
 one would derive 
the alternative proof  for the famous Theorem from Hirzebruch-Zagier as presented
recently by Funke-Millson \cite{FM2}. In our appendix we also indicate such a direct proof (see Prop. \ref{PropthetaTt}).

\newpage
\section{Generating series for arithmetic special cycles}
 
 We assume the reader to be familiar with the basic notations and concepts
 of arithmetic intersection theory as it can be found in \cite{SABK} or \cite{Sou}. 
 
 Without any problems we can extend $X$ to the arithmetic threefold $\mathbb P^1_\Z \times  \mathbb P^1_\Z$, which by abuse of notation will denoted again by $X$. We let $T(N)$ now also denote the Zariski closure of $T(N)$ in $X$. One may define 
 $T(N)$ also by a special modular description, but this will not be needed in our context now.

We had observed that Kudla's Green function $\Xi (m)$ is not a Green function   for the divisor $T(m)$ on  $X$. Therefore,   
using the modification proposed before, we introduce the modified arithmetic special cycle
\begin{align*}
\hat{Z}_\rho (m) := (T(m),\tilde{\Xi }_\rho (v,z,m)) \in \widehat{CH}^1(X). 
\end{align*}

A reasonable hope would be that this modification is 
according to Kudla's conjectures modular also, 
by this we mean that the modified Kudla generating series
\begin{align}\label{kudla-gen-series}
\hat{\phi }_{K}(q) := \sum \hat{Z}_\rho (m)q^m \nne
\end{align}
is a modular form with coefficients in $\widehat{CH}^1(X)$ and we will show that this is indeed the case.
 
If \eqref{kudla-gen-series} were modular,  then
for all linear maps $L: \widehat{CH}^1(X)  \rightarrow \mathbb R,$ one would have 
\begin{align*}
\sum_{m\in \mathbb{Z}} L(\hat{Z}_\rho (m))q^m
\end{align*}
as a (nonholomorphic) $\mathbb R-$valued modular form in the usual sense. 
Our strategy is to show that is suffices to check the modularity just for one
suitable linear map.
 
\nn {\bf Remark.}  Let us denote by $\widehat{c}_1( \bar{\mathcal L})$ the first arithmetic Chern class of $\bar{\mathcal L}(12,12)$
and choose $L(-) = \widehat{c}_1( \bar{\mathcal L})^2\cdot (-)$. The  modified second  sub-conjecture of  
Kudla as mentioned in the introduction then takes the form
\begin{align}\label{eq:kudla-modi}
\sum_{m \in \mathbb{Z}} \widehat{c}_1( \bar{\mathcal L})^2\cdot \hat{Z}_\rho (m)\,q^m = \mathbb E'_2(\tau ,1) + f_\rho (\tau )\nne
\end{align}
with a certain modular form $f_\rho (\tau )$ and a suitable normalised non-holomorphic Eisenstein series $\mathbb E_2(\tau ,s)$. A proof of that
will be in a forthcoming paper.

 \nn{\bf Modularity Detour.}
Analogously to [BKK] p.168 we know the modularity of 
\begin{align}\label{phicaneis}
\hat{\phi }_{\rm FS}(q) := \sum _{N\in\mathbb Z}\widehat{T(N)}q^N    \nne
\end{align}
where
\begin{align}
\widehat{T(N)} = \begin{cases} 0 &{\rm if}\,N<0\\
-(1/24)\widehat{c}_1(\overline{\mathcal L}(12,12),\|\,\|_{\rm FS})+ (1/8\pi v)[D,g_D(\tau )] &{\rm if}\,N=0\\ 
[\mathcal T_N,g_N]&{\rm if}\,N>0\\ \end{cases}\nne
\end{align}
with
\begin{align*}
 g_D(\tau ) &= - \log \|\Delta (z_1)\Delta (z_2)\|^2_{\rm FS} \\\nonumber
  g_N(\tau ) &=  - \log \|\Psi_N (z_1,z_2)\|^2_{\rm FS}\nne
\end{align*}
where $\Psi_N(z_1,z_2) \in M_{\sigma (N)12,\sigma (N)12}$ has divisor $T(N)$ 
and the Fubini-Study metric of a  bi-modular form $f$ of weight $(12 k, 12 k)$ equals
\begin{align*}
\|f\|^2_{\rm FS}(z_1,z_2) = \frac{|f(z_1,z_2)|^2}{ (|E_4^3(z_1)|^2 +|\Delta(z_1)|^2)^k \cdot  (|E_4^3(z_2)|^2 +|\Delta(z_2)|^2)^k}.
\end{align*}

By construction one has
 \begin{align}
\widehat{T(N)} = \hat{c_1}(\bar{\mathcal L}(12,12))\begin{cases} 0 &{\rm if}\,N<0\\
(-(1/24)+ (1/8\pi v)) &{\rm if}\,N=0\\ 
\sigma (N)&{\rm if}\,N>0.\\ \end{cases}\nne
\end{align}
 
We therefore have the decomposition
\begin{align*} 
\hat{\phi }_{\rm FS}(q) := \sum \widehat{T(N)}q^N = \widehat{c_1}(\bar{\mathcal L}(12,12))\otimes E_2(\tau ,1)  \nne
\end{align*}
We want to use this for our proof of the modularity of $\hat{\phi }.$\\

\nn {\bf Lemma.}\label{lemker} \emph{There exist a family of functions $\{b(\tau,m)\}_{m \in \mathbb{Z}}$ such that 
\begin{align*}
\hat{\phi }_K(q) + \sum_{m \in \mathbb{Z}} b(\tau,m) q^m 
\end{align*}
is a $\widehat{CH}^1(X)$-valued modular form.}\\

{\bf Proof.} Given the modularity of  $\hat{\phi }_{\rm FS}(q)$, we also can look for the modularity of 
\begin{align}
\tilde{\phi }(q) := \hat{\phi }_{\rm FS}(q)- \hat{\phi }_K(q).  \nne
\end{align}
But now $\tilde{\phi}(q)$ is purely analytic object, 
indeed it is in the image of the ${\rm a}$-map in the exact sequence for the arithmetic Chow groups given in 
\cite{SABK}. The coefficients of $\tilde{\phi}(q)$ are   smooth functions.   We have shown already, that the $q$-series 
$dd^c \tilde{\phi}(q)$ is a differential form valued modular form. 
Moreover its coefficients are  smooth
differential forms homologous to zero. Thus we can apply the general 
mechanism of Green functions to invert the linear operator $dd^c$. 
Since these operators do not effect the $\tau$ variable, we get that
up to  
a family of functions $\{b(\tau,m)\}_{m \in \mathbb{Z}} \in \ker dd^c$
the $q$-series $\tilde{\phi }(q)$ is modular. 
 Now by the previous discussion $\hat{\phi }_{\rm FS}(q)$ is already modular
 and therefore the claim follows. \hfill $\Box$\\

We try to determine these functions $\{b(\tau ,m)\}_{m\in \mathbb Z}$ by  means of a linear form
that does not map the image of the ${\rm a}$-map to zero. We choose
as our linear from the height function  
\begin{align}
L: \widehat{CH}^1(X) &\longrightarrow \mathbb R,\,\, \hat{\alpha } \longmapsto {\rm ht}_{\hat{\alpha }}(P)    \nne
\end{align}
for a point $P \in X =  (\mathbb P^1_\mathbb Q\times \mathbb P^1_\mathbb Q)(\mathbb Q)$. It turns out to be suitable to choose 
\begin{align}
P = (i,\infty )    \nne
\end{align}
hence, we have to check the failure of
\begin{align}
L(\hat{\phi }_K) =\sum_{m\in\Z} {\rm ht}_{\hat{Z}_\rho (m)}(i,\infty )q^m.    \nne
\end{align}
being a modular form. 

\nn {\bf Proposition.}\label{Prophtm} For $m\not=0$ one has
\begin{align}\label{eqhtm}
{\rm ht}_{\hat{Z}_\rho (m)}(i,\infty ) = \begin{cases}
0 & {\rm if }\,\, m<0 \\ 
4\pi \sigma (m) & {\rm if }\,\, m>0. \nne
\end{cases}
\end{align}
{\bf Proof.} Using [BBK] (1.16) p.21, we get
\begin{align}
{\rm ht}_{\hat{Z}_\rho (m)}(i,\infty ) = (T(m)\cdot (i,\infty ))_{{\rm fin}} + \tilde{\Xi }_\rho (m)(i,\infty ) \nne
\end{align}
By closer inspection, one remarks that there is no geometric intersection at the finite places, thus the first term on the right side vanishes. For $m<0$ the second  term vanishes also. 
For $m>0$ we deduce from Theorem \ref{ThemXim}  
\begin{align*}
\tilde{\Xi }_\rho (v,z,m) = \Xi (v,z,m)  - \rho (z)t/\sqrt v (B(v,s,m)-I(v,s,m))
\end{align*}
and Theorem \ref{ThnearD} tells us that near the boundary divisor $D$ one has up to a smooth function vanishing on $D$
\begin{align*}
\Xi (v,z,m) \sim (1/2)\sum_{-bc=m} -\log|q_1^{|c|}-q_2^{|b|}|^2 - t/\sqrt v (B(v,s,m)-I(v,s,m)).
\end{align*}
Hence we get
\begin{align*}
\tilde{\Xi }_\rho (v,(i,\infty ),m) &= (1/2)\sum_{-bc=m} -\log|e^{-2\pi |c|}|^2\\
&= (1/2)\sum_{-bc=m} 4\pi |c| = 4\pi \sigma (m). 
\end{align*}
\hfill$\Box$\\

It remains to treat the case $m=0.$ 

\nn {\bf Proposition.}\label{Propht0} For $m=0$ one has
\begin{align}\label{eqht0}
{\rm ht}_{\hat{Z}_\rho (0)}(i,\infty ) = 4\pi  c_0 .  \nne
\end{align}
{\bf Proof.} 
Here we have to choose carefully a representative for the class 
$\hat{Z}_\rho (0) \in \widehat {CH}^1(X)$ such that we again may use formula (1.16) in [BBK] for ${\rm ht}_{\hat{Z}_\rho (0)}.$ 
As usual, we have 
\begin{align*}
T(0) = c_0D = c_0 (\mathbb P^1\times \infty + \infty \times \mathbb P^1),\,\, c_0 = -(1/24) + (1/8\pi v) 
\end{align*}
and, with $\varrho = e^{\pi i/3},$ now we take (using again the notation for $\widehat {{\rm div}}$ as in [Sou])
\begin{align}\label{hatZ0}
\hat{Z}_\rho (0) &= [T(0),\hat{\Xi }_\rho (0)] + c_0 \widehat {{\rm div}}(j(z_2))
\notag \\
&= [T(0),\hat{\Xi }_\rho (0)] + c_0[\mathbb P^1\times \varrho -   \mathbb P^1 \times \infty, -\log |j(z_2|^2] \notag \\
&= [c_0(\infty \times \mathbb P^1 - \mathbb P^1 \times \varrho ),\tilde{\Xi }_\rho (v,z,0) - c_0\log |j(z_2|^2] \nne
\end{align}
and, via (1.16) in [BBK], get 
\begin{align*}
{\rm ht}_{\hat{Z}_\rho (0)}(i,\infty ) = c_0((\infty \times \mathbb P^1 - \mathbb P^1 \times \varrho )\cdot (i,\infty ))_{{\rm fin}}
+ (\tilde{\Xi }_\rho (v,z,0)- c_0\log |j(z_2|^2)(i,\infty )  
\end{align*}
Again we inspect that the  first term on the right hand side vanishes 
as there is no geometroic intersection. For the evaluation of the second term, we 
take $z_2 = iT$ and, hence, for large $T$ have 
\begin{align}\label{jot}
\log |j(z_2)| = \log |q^{-1}(1+744q+\dots)| \simeq \log e^{2\pi T} = 2\pi T
\nne
\end{align}

From Theorem \ref{ThemXi00} we know
\begin{align*}
\tilde{\Xi }_\rho (v,z,0) &= \Xi (v,z,0) - \rho (z)(\check{\Xi }(v,z,0) -(t/2v)(s+1/s))\\
&\sim (1/12 - 1/(4\pi v))\log |q_1q_2| 
\end{align*}
i.e., for $z=z_T= (i,iT)$
\begin{align*}
\lim _{T\rightarrow \infty }\tilde{\Xi }_\rho (v,z_T,0) 
&= \lim _{T\rightarrow \infty }(1/12 - 1/(4\pi v))(-2\pi - 2\pi T) = \lim _{T\rightarrow \infty }4\pi c_0(1 + T). 
\end{align*}
If we join this with (\ref{hatZ0}) and (\ref{jot}),  the $T-$terms cancel and we get the desired result
\begin{align}
{\rm ht}_{\hat{Z}_\rho (0)}(i,\infty ) = 4\pi c_0. \nne
\end{align} \hfill $\Box$

 \nn {\bf Main Theorem.} \emph{For any choice of $\rho$   the generating function of modified special arithmetic cycles
\begin{align*} 
\sum_{m \in \ZZ}  \widehat{Z}_\rho(v,z,m) q^m  
\end{align*} 
is an $\widehat{CH}^1(X) $-form valued  weight $2$ modular form for the group $\SL_2(\ZZ)$, where  $q = e^{2\pi i\tau }$ and $v = \Ima \tau >0$.}\\

{\bf Proof.} 
From the Propositions \ref{Prophtm} and \ref{Propht0} we learn that 
\begin{align*}
L(\hat{\phi }_K) =\sum_{m\in\Z} {\rm ht}_{\hat{Z}_\rho (m)}(i,\infty )q^m = 4\pi E_2(\tau ,1)    \nne
\end{align*}
is already modular. Hence, in Lemma \ref{lemker} we have all $b(\tau ,m)=0$ and the claim follows. 
  \hfill $\Box$

\newpage
\addcontentsline{toc}{section}{\bf Appendix}

{\bf \huge Appendix}

\hspace{1cm}
 
\appendix
 \setcounter{section}{0}
\section{A dictionary to the orthogonal point of view}

In the appendix we give a dictionary to the orthogonal point of view and prove some 
more or less well known but not well documented facts concerning this orthogonal world.\\

There are several ways to introduce coordinates into our calculations. 
We describe two of these which come up in our text.\\

{\bf Our orthogonal world}\\

Let $(V,q)$ be the quadratic space of signature (2,2) given by
\begin{align}
\label{V}\nne
V&= \{ M =\left(
 \begin{smallmatrix} a & b\\ c & d\end{smallmatrix}\right) \in 
\operatorname{M}_2(\R) \}
\end{align}
and the forms 
\begin{align}
\label{q} \nne
q(M)=\det(M)= (1/2)(M,M), \qquad (M,M')= ( ad' + a'd - bc' - b'c).
\end{align}
\nn {\bf Remark.} Sometimes it will be convenient to use the following  identifications.

\begin{align}
\label{a} 
M \equiv \, {\bf a} =& {}^t(a,b,c,d) = {}^t(y_1, - y_2,y_3,y_4) \nonumber\\
=& (1/\sqrt 2){}^t(x_1+x_4,-x_2-x_3,x_2-x_3,x_1-x_4)\nonumber\\ 
x =& {}^t(x_1,x_2,x_3,x_4) = (1/\sqrt 2){}^t(y_1+y_4,y_2+y_3,y_2-y_3,y_1-y_4)\nne\\
=& (1/\sqrt 2){}^t(a+d,c-b,-c-b,a-d).\nonumber
\end{align}

With
$$
Q =  \begin{pmatrix}1&&&\\&1&&\\&&-1&\\&&&-1 \end{pmatrix},\quad
\tilde Q =  \begin{pmatrix}&&&1\\&&1&\\&1&&\\1&&& \end{pmatrix}.
$$
one has
$$
(M,M) = 2\det M = 2(ad-bc) = 2(y_1y_4+y_2y_3) = {}^ty \tilde Q y = {}^tx Q x = (x_1^2+x_2^2-x_3^2-x_4^2).
$$
\nn {\bf Definition.} For the identity components of the corresponding orthogonal groups we write
\begin{align}
\label{G}\nne 
G := \Or _0(Q) = \SO_0(2,2) \simeq \tilde G := \Or _0(\tilde Q).
\end{align}
As usual (e.g. \cite{BF}) the symmetric space associated with $G$ may be identified with the
set 
\begin{align}
\label{D}\nne 
\mathbb D=\{\mbox{ oriented negative 2-planes } X \subset V\,\},
\end{align}
i.e. $ X = \,<M_1,M_2>,\, M_j \in V, \,\det M_j < 0, \,(M_1,M_2) = 0.$  
 It is well-known that $\mathbb D$ has
two connected components $\mathbb D^+$ and $\mathbb D^-$, each of it is isomorphic to
the product of two upper half planes $\mathbb H^2$.  There are several ways to fix coordinates. In the sequel we use
the isomorphism $\mathbb H^2 \to \mathbb D^+$ (and hence neglect the +) given by 
$$
(z_1,z_2) \longmapsto <\Rea Z,\Ima Z>,
$$
 where 
\begin{align}
\label{Z} \nne
Z=\begin{pmatrix} z_1 z_2 & z_1\\ z_2 &1\end{pmatrix}.
\end{align}
We will dwell on this a bit: One has the fact that for $\gamma =( \gamma _1,\gamma _2) \in \bar G = \operatorname{SL}_2(\R)^2$ 
\begin{align}
\label{M'}\nne 
M \longmapsto  M' = M^\gamma := \gamma_1 M\, {}^t\gamma_2
\end{align}
resp.
\begin{align}
\label{Agamma}\nne 
\iota :\,{\bf a} \longmapsto {\bf a'} = A(\gamma ){\bf a}
\end{align}
is an isometry of $V.$ \\
\nn {\bf Remark.} We have a homomorphism
$$
\iota ':\, \bar G = \operatorname{SL}_2(\R)^2 \longrightarrow \tilde G, \quad \gamma \longmapsto A'(\gamma )
 $$
where (as given by an easy calculation) for 
$$
\gamma _1 = \begin{pmatrix} \alpha _1 & \beta _1\\ \gamma _1 & \delta _1\end{pmatrix}, 
\gamma _2 = \begin{pmatrix} \alpha _2 & \beta _2\\ \gamma _2 & \delta _2\end{pmatrix}
$$
one has
\begin{align}
\label{A'} \nne
 A'(\gamma) = 
\begin{pmatrix} \alpha _1 \alpha _2& -\alpha _1\beta _2&\beta _1\alpha _2&\beta _1\beta _2\\
-\alpha _1 \gamma _2 & \alpha _1\delta _2&-\beta _1\gamma _2&-\beta _1\delta _2\\
\gamma _1\alpha _2&-\gamma _1\beta _2&\delta _1\alpha _2&\delta _1\beta _2\\
\gamma _1\gamma _2&-\beta _1\delta _2&\delta _1\gamma _2&\delta _1\delta_2
\end{pmatrix} .
\end{align}
As special cases we have with $ t = \alpha _1\alpha _2, s = \alpha _1/\alpha _2$
$$
\begin{array}{rcl}
 (\begin{pmatrix} \alpha _1 & \\ & \alpha _1^{-1}\end{pmatrix},\begin{pmatrix} \alpha _2 & \\ & \alpha _2^{-1}\end{pmatrix})
&\longmapsto & \begin{pmatrix} \alpha _1\alpha _2 & &&\\ & \alpha _1\alpha _2^{-1}&&\\
&&\alpha _1^{-1}\alpha _2&\\&&&(\alpha _1\alpha _2)^{-1}\end{pmatrix} 
=: a(t,s),\\
(\begin{pmatrix} 1& \beta _1\\ & 1\end{pmatrix},(\begin{pmatrix} 1 & \\ & 1\end{pmatrix})) &\longmapsto & 
\begin{pmatrix} 1 && \beta _1 &\\ & 1 & &-\beta_1\\
&&1&\\&&&1 \end{pmatrix} =:n'(0,\beta _1)\\
((\begin{pmatrix} 1 & \\ & 1\end{pmatrix},\begin{pmatrix} 1& \beta_ 2\\ & 1\end{pmatrix})
&\longmapsto & \begin{pmatrix} 1 &-\beta _2& &\\ &1 & &\\
&& 1 &\beta _2\\&&&1 \end{pmatrix} =:n'(\beta _2,0)\\
\end{array}
$$
\nn {\bf Remark.} \label{Remiota }The matrix $A(\gamma )$ describing the map $\iota $ (\ref{Agamma}) in the ${\bf a-}$coordinates 
is the same as $A'(\gamma )$ above but without the minus-signs. Moreover we write $n(0,\beta _2),n(\beta _1,0)$ 
by leaving out the minus signs in the corresponding $n'-$elements.\\

The homomorphism $\bar G \rightarrow \tilde G$ induces an isomorphism of the symmetric spaces
$$
\bar G/\bar K = \mathbb H^2 \longrightarrow \tilde G/\tilde K,
$$
namely, with $z_j = x_j +y_j i \in \mathbb H $ and
\begin{align}
\label{gz} \nne
g_{z_j} = \begin{pmatrix} y_j^{1/2} & x_jy_j^{-1/2}\\ &y_j^{-1/2} \end{pmatrix} \in \SL(2,\R)
\end{align}
(j = 1,2) one has 

\begin{align}
\label{Az} \nne
A(z) := A(g_{z_1},g_{z_2}) = \left(\begin{smallmatrix}
 \sqrt{y_1 y_2} &  \sqrt{y_1/y_2} x_2  &  \sqrt{y_2/y_1} x_1 &  \sqrt{y_1 y_2}^{-1} x_1 x_2\\ 
0 & \sqrt{y_1/y_2} &0 & \sqrt{y_1 y_2}^{-1} x_1\\
0& 0 & \sqrt{y_2/y_1} &\sqrt{y_1 y_2}^{-1} x_2 \\
0& 0 & 0 &\sqrt{y_1 y_2}^{-1}
\end{smallmatrix}\right).
\end{align}

 Now, fixing as base point of
$\mathbb D$ the negative 2-plane spanned by $M_1 := \left(\begin{smallmatrix} 1 & \\
    & -1\end{smallmatrix}\right)$ and $M_2 = \left( \begin{smallmatrix} & -1\\ -1&
    \end{smallmatrix}\right)$
we get 
$$
g_{z_1}M_1{}^t g_{z_2} =(y_1y_2)^{-1/2}\begin{pmatrix} y_1y_2-x_1x_2& -x_1\\ -x_2& -1 \end{pmatrix} = - (y_1y_2)^{-1/2}\Rea Z
$$
$$
g_{z_1}M_2{}^t g_{z_2} =(y_1y_2)^{-1/2}\begin{pmatrix} -x_1y_2-x_2y_1& -y_1\\ -y_2& 0 \end{pmatrix} = - (y_1y_2)^{-1/2}\Ima Z
$$
where $Z$ is given as above by $Z = (\begin{smallmatrix} z_1z_2&z_1\\z_2&1\end{smallmatrix}).$ This explains the formula for the isomorphism 
\begin{align}\label{isom-H2D}\nne
\mathbb H^2 \longrightarrow \mathbb D \notag \\
z=(z_1,z_2)&\mapsto X(z):=\left< g_{z_1} \left(\begin{smallmatrix} 1 & \\  & -1\end{smallmatrix}\right)
 \,{}^t g_{z_2},   g_{z_1} \left(\begin{smallmatrix}  & -1\\ -1 & \end{smallmatrix}\right)
 \,{}^t g_{z_2} \right>.\nonumber
\end{align}
 also given above. We observe the relations

\begin{align*}
-2y_1 y_2 &=  (\Rea Z, \Rea Z) = ( \Ima Z, \Ima Z) \\ 
0&= (\Rea Z, \Ima Z) 
\end{align*}

We also note that $(i,i) \in \mathbb H^2$ corresponds to our basepoint $(M_1,M_2) \in \mathbb D$.\\

{\bf The majorant and its kernel $R$ }\\

We shall decompose $M \in V$ with respect to  the negative 2-plane
$X(z)$ into its positive and negative parts:
\begin{align*}
M = M' + \alpha \Rea Z + \beta  \Ima Z,
\end{align*}
i.e., such that
\begin{align*}
(M', \Rea Z) &= (M' , \Ima Z) = 0\\
\alpha &= - \frac{(M,\Rea Z)}{2y_1 y_2}\\
\beta &= - \frac{(M,\Ima Z)}{2y_1 y_2}.
\end{align*}
Now, we see
\begin{align*}
(M,M) &= (M',M') + \alpha^2 (\Rea Z,\Rea Z) + \beta^2 (\Ima Z,\Ima Z) \\
&= (M',M') - \left( (M, \Rea Z)^2 + (M, \Ima Z)^2 \right)/ (2y_1y_2).
\end{align*}
\nn {\bf Remark.} The (positive definite!){\it majorant} with respect to $X(z)$ is given by
\begin{align}
\label{maj} 
(M,M)_z &= (M',M') - \alpha^2 (\Rea Z,\Rea Z) - \beta^2 (\Ima Z,\Ima Z) \nonumber\\
&= (M',M') + \left( (M, \Rea Z)^2 + (M, \Ima Z)^2 \right)/ (2y_1y_2)\nonumber\\
&=(M,M) + \left( (M, \Rea Z)^2 + (M, \Ima Z)^2 \right)/ (y_1y_2) \nne\\
&=: (M,M) + 2 R(z,M),\nonumber
\end{align}
where
\begin{align}
\label{R} \nne
R(z,M) = R(z_1,z_2,M) = 
\frac{|a  - b z_2 - c z_1 +d  z_1 z_2 |^2}{2 y_1y_2}.
\end{align}
The last equality is easily checked by a straightforward
calculation. Lacking a better expression, we call $R$ the {\it kernel} of the majorant.\\

Note that at the base point we have
\begin{align*}
(M,M)_0 &:= (M,M)_{(i,i)} 
=  2(ad -bc) + ((a-d)^2 +(b+c)^2)\\
&=  ( a^2 +b^2 +c^2 +d^2).
\end{align*}
        
As we shall use it later, we note\\

\nn {\bf Key Remark}: One has
\begin{align}
\label{Key} \nne
{}^t(A(z)^{-1}{\bf a})A(z)^{-1}{\bf a} = 2R(z,M) + (M,M).
\end{align}
Here we use the notation from above: we have in the ${\bf a}$-coordinates
\begin{align}
\label{t,s} 
A(z) =& n(x_2,x_1)a(t,s), \nonumber\\t =& \sqrt{y_1y_2}, s = \sqrt{y_1/y_2}\nne
\end{align}
with
$$
n(\mu _1,\mu _2) = \begin{pmatrix}
 1 &   \mu_1  &   \mu_2 &   \mu_1 \mu_2\\ 
 & 1 & &  \mu_2\\
&  & 1 & \mu_1 \\
&  &  &1
\end{pmatrix},\quad
a(t,s) =
\begin{pmatrix}
 t &  && \\
&s && \\
& & s^{-1} & \\
&&& t^{-1} 
\end{pmatrix}. 
$$
Hence we get 
\begin{align}
\label{a'}\nne
{\bf a'} = A(z)^{-1}{\bf a} = \begin{pmatrix}
 t^{-1}(a - x_2b - x_1c + x_1x_2d)\\ 
 s^{-1}(b - x_1d)\\
s(c - x_2d) \\
td
\end{pmatrix}.
\end{align}
From here the formula in the Remark can be verified directly 
(but this is also clear from the general majorant framework!).\hfill$\Box$\\

\nn {\bf Corollary}: For $g \in \bar G$ one has the invariance property
\begin{align}
\label{inv} \nne
R(g(z),M^g) = R(z,M).
\end{align}
{\bf Proof}:
We have
$$
\begin{array}{rl}
A(g(z))^{-1}{\bf a} =& (A(g)A(g_z))^{-1}{\bf a}\\[.3cm]
=& A(g_z)^{-1}A(g)^{-1}{\bf a}\\[.3cm]
=& A(z)^{-1}{\bf a}^{-g}\\[.3cm]
\end{array}
$$
with ${\bf a}^{-g}$ describing the entries of $M^{-g}$ and
$$
\begin{array}{rl}
{}^t(A(g(z))^{-1}{\bf a})(A(g(z))^{-1}{\bf a}) =& 2 R(g(z),M) + (M,M)\\[.3cm]
=& {}^t(A(z)^{-1}{\bf a}^{-g})(A(z)^{-1}{\bf a}^{-g})\\[.3cm]
=& 2 R(g(z),M^{-g}) + (M^{-g},M^{-g}),\\[.3cm]
\end{array}
$$
i.e.
$$
R(g(z),M) = R(z,M^{-g}).
$$
\hfill$\Box$\\
{\bf Another description of $\mathbb D$}\\

If one changes the action of the group, one has a description of 
the space $\mathbb D$ which slightly differs from the one given above 
which we introduced following Funke-Millson and others.\\

We again have
$$
\mathbb D=\{ \mbox{ oriented negative 2-planes } X \subset V\,\}.
$$
$$
\bar G = \SL(2,\R)^2
$$  
and this time take the action of $\bar G$ on $V$ given by
\begin{align*}
 M \longmapsto M' = M^\gamma = \gamma_1 M\, \gamma_2^{-1}
\end{align*}
and hence get a homomorphism of $\bar G$ to $\Or(\tilde Q)$ given by
$\gamma = ( \gamma_1, \gamma_2) \mapsto A(\gamma)$  with
$$
\gamma _1 = \begin{pmatrix} \alpha _1 & \beta _1\\ \gamma _1 & \delta _1\end{pmatrix}, 
\gamma _2 = \begin{pmatrix} \alpha _2 & \beta _1\\ \gamma _2 & \delta _2\end{pmatrix}
$$
and
\begin{align}
\label{A-1}\nne
\gamma := (\gamma _1,\gamma_2) \longmapsto A(\gamma) = 
\begin{pmatrix} \alpha _1 \delta _2& -\alpha _1\gamma _2&\beta _1\delta _2&-\beta _1\gamma _2\\
-\alpha _1 \beta _2 & \alpha _1\alpha _2&-\beta _1\beta _2&-\beta _1\alpha _2\\
\gamma _1\delta _2&-\gamma _1\gamma _2&\delta _1\delta _2&-\delta _1\gamma _2\\
-\gamma _1\beta _2&\gamma _1\alpha _2&-\delta _1\beta _2&\delta _1\alpha_2
\end{pmatrix}. 
\end{align}

 Now, as above, fixing as base point of
$\mathbb D$ the negative 2-plane spanned by $M_1 := \left(\begin{smallmatrix} 1 & \\
    & -1\end{smallmatrix}\right)$ and $M_2 := \left( \begin{smallmatrix} & -1\\ -1&
    \end{smallmatrix}\right)$
we get 
$$
g_{z_1}M_1g_{z_1}^{-1} =(y_1y_2)^{-1/2}\begin{pmatrix} y_1& -x_1y_2-x_2y_1\\0& -y_2\end{pmatrix} = - (y_1y_2)^{-1/2}\Re \tilde Z
$$
$$
g_{z_1}M_2 g_{z_2}^{-1} =(y_1y_2)^{-1/2}\begin{pmatrix} -x_1&x_1x_2-y_1y_2-x_2\\-1& x_2 \end{pmatrix} = - (y_1y_2)^{-1/2}\Im \tilde Z
$$
with
$$
\tilde Z = \begin{pmatrix} - \bar z_1 & \bar z_1\bar z_2 \\ -1 & \bar z_2\end{pmatrix}.
$$
This explains an isomorphism 
\begin{align}\label{isom-H2D'}
\mathbb H^2 &\to \mathbb D \nonumber \\
z=(z_1,z_2)&\mapsto X(z):=\left< g_{z_1} \left(\begin{smallmatrix} 1 & \\  & -1\end{smallmatrix}\right)
 \,g_{z_2}^{-1},   g_{z_1} \left(\begin{smallmatrix}  & -1\\ -1 & \end{smallmatrix}\right)
 \,g_{z_2}^{-1} \right>.\nne
\end{align}
 We observe the relations
\begin{align*}
-y_1 y_2 &=  (\Re \tilde Z, \Re \tilde Z) = ( \Im \tilde Z, \Im \tilde Z) \\ 
0&= (\Re \tilde Z, \Im \tilde Z) 
\end{align*}
and, with the same procedure as every year, this time come to the majorant
$$
(M,M)_z =(M,M) + 2 \tilde R(z,M),
$$
where
\begin{align}
\label{RW}
\tilde R(z,M) :&= \tilde R(z_1,z_2,M) = ((M,\Re \tilde Z)^2+ (M,\Im \tilde Z)^2)/( y_1y_2)\nonumber\\[.3cm]
&=\frac{|az_2  + b  - c z_1 z_2 - dz_1|^2}{2 y_1y_2}.\nne
\end{align}
We have $R(z,M) = 0$ exactly for $z_1 = M(z_2). $\\

\nn {\bf Remark.} The coordinates in this description are related to the {\it old} coordinates by
$$
z_1 \longmapsto \bar z_1, \,\,z_2 \longmapsto -1/\bar z_2.
$$
This can be seen as follows: We denote the new coordinates by $w = (w_1,w_2).$ 
A point of $\mathbb D$ is given in the old coordinates by the plane
$$
X(z) = <\Rea Z, \Ima Z >, \,Z = \begin{pmatrix}  z_1z_2 &  z_1 \\ z_2& 1\end{pmatrix}
$$
and in the new coordinates by
$$
\tilde X(w) = < \Rea W, \Ima W >, \, W = \begin{pmatrix} - \bar w_1 & \bar w_1\bar w_2 \\ -1 & \bar w_2\end{pmatrix}.
$$
Hence one has with real constants
$$
\Rea W = \alpha \Rea Z + \beta \Ima Z,\,\,\Ima W = \gamma \Rea Z + \delta \Ima Z,
$$
i.e.,
$$
\begin{array}{rl}
W = \Rea W + i \Ima W =& (\alpha +i\gamma )\Rea Z + (\beta +i\delta )\Ima Z\\
=& (\alpha + \delta +i(\gamma -\beta )Z/2 + (\alpha -\delta + i(\beta + \gamma )\bar Z\\
=& (\alpha +i\gamma )Z =: \eta Z
\end{array}
$$
if we take $\alpha =\delta , \gamma = - \beta .$ This leads to
$$
-\bar w_1 = \eta z_1z_2, \bar w_1\bar w_2 = \eta z_1, -1 = \eta z_2, \bar w_2 = \eta, 
$$
i.e.,
\begin{align}
\label{newcoord}\nne
w_2 = - 1/\bar z_2, \,\, w_1 = \bar z_1.
\end{align}
\hfill$\Box$\\
With $w_j =:u_j + iv_j$ one has
$$
t^2 = y_1y_2 = v_1v_2/\mid w_2 \mid ^2
$$
and
$$
dz_1 = d\bar w_1, \,dz_2 = d\bar w_2/\bar w_2^2.
$$
We get the nice transformation formulae
$$
\frac{dz_j\wedge d\bar z_j}{y_j^2} = - \frac{dw_j\wedge d\bar w_j}{v_j^2}, \,j = 1,2
$$
and the less nice formula
$$
\frac{dz_1\wedge d\bar z_2 \pm dz_2 \wedge d\bar z_1}{y_1y_2} = \frac{\bar w_2^2d\bar w_1\wedge dw_2 \pm w_2^2d\bar w_2^2\wedge dw_1}{v_1v_2\mid  w_2\mid ^2}.
$$
In the sequel we will stick to the {\it old} coordinates from the previous subsection 
as one can change results to the {\it alternative} coordinates 
by the formulae just obtained.\\

{\bf Special cycles}\\

The zero-set of our function $R(z,M)$ has an important geometric meaning.\\

\nn {\bf Remark.} $R(z,M) = 0$ is exactly the case for $z_1 = \frac{-bz_2+a}{-dz_2+c},$ i.e.,
$$
z \in T_M := \{ (z_1,z_2); z_1 = M\,S(z_2)\}, \,\,S = (\begin{smallmatrix} &1\\ -1 &\end{smallmatrix}) \}.
$$

\nn {\bf Remark.} For $\det M =m$ one also finds in the literature the description
$$
T_m = \{z \in \mathbb H^2, z \perp M \}.
$$
This is another way of saying that the Matrix $Z = \begin{pmatrix} z_1 z_2 & z_1\\ z_2 &1\end{pmatrix}$ associated to z above has
\begin{align}\label{eqTm}
(Z,M) = (\begin{pmatrix} z_1 z_2 & z_1\\ z_2 &1\end{pmatrix},\begin{pmatrix} a&b\\ c &d\end{pmatrix}) = a +z_1z_2d -bz_2-cz_1 = 0.\nne
\end{align}
If we change $z_2 = - 1/\tilde z_2,$ for $\det M =m \not= 0$ we get
$$
R(z,M) = 0 \Longleftrightarrow z_1 = M(\tilde z_2),
$$
i.e. $R(z,M) = 0$ is exactly the case for the Hirzebruch-Zagier divisor $T_M$ here in the 
form $T_M = \{(z_1,z_2 = -1/\tilde z_2), z_1 = M(\tilde z_2)\}.$\\

\nn {\bf Remark.} The Hirzebruch-Zagier divisor in the usual form $\tilde T_M = \{(\tilde z_1,\tilde z_2), \tilde z_1 = M(\tilde z_2)\}$
would have shown up 
as the zero-locus of $\tilde R(\tilde z,M)$ if we had chosen our second coordinization of $\mathbb D,$ namely the one coming from the 
action of $\bar G$ on $V$ given by 
$$
M \longmapsto \tilde M' = \gamma _1M\gamma _2^{-1} 
$$

Here, as above in (\ref{RW}), we come to 
$$
\tilde R(\tilde z,M) =  \frac{|a\tilde z_2  + b  - c \tilde z_1 \tilde z_2 - d\tilde z_1|^2}{2 \tilde y_1\tilde y_2}. 
$$
But anyway, if one looks at the divisor $T(m)$ on $X$ as defined in (\ref{T(N)}) in section 1, clearly
it is independent of the coordinates used to its description.\\

{\bf Boundary}\\

It is quite natural that the space $X_0 = \tilde{G}/\tilde{K} \simeq (SL(2,\mathbb Z)\times \SL(2,\mathbb Z))\backslash (\mathbb H \times \mathbb H)$ is 
compactified by taking as boundary $D = \mathbb P_1 \times \{\infty \}\cup \{\infty \}\times \mathbb P_1.$ This is consistent 
with the general compactification theory.  $X_0$ is covered by
$\mathbb D = \tilde{G}/\tilde{K} = \mathbb H \times \mathbb H$ and following the uniform 
construction of a reductive Borel-Serre compactification as in \cite{BJ} p. 338ff we 
look at the standard rational parabolics in $\tilde{G}.$ These are ${\bf P_0, P_1, P_2},$ the images in $\tilde{G}$ of respectively
${\bf B}\times {\bf B}, {\bf B}\times \SL(2,\mathbb Q), \SL(2,\mathbb Q)\times  {\bf B} , \,{\bf B}$ the usual upper triangular 
subgroup of $\SL(2,\mathbb Q).$ For each group we take the Langlands decomposition 
$P = N_PA_PM_P$ and with $K_P = M_P \cap \tilde{K}$ and get the boundary component $e({\bf P}) = P/N_PA_PK_P.$ 
\nn {\bf Remark.}  In our case we have 
$$
\tilde P_0 = N_0A_0,\,\tilde P_1 = N_1A_1M_1,\,\tilde P_2 = N_2A_2M_2
$$
with 
\begin{align*}
N_0 =& \{n'(x_2,x_1), x_1, x_2 \in \mathbb R \}, A_0 =\{a((y_1y_2)^{1/2},(y_1/y_2)^{1/2}),\\
N_1 =& \{n'(0,x_1), x_1 \in \mathbb R \}, A_1 =\{a(y_1^{1/2},y_1^{1/2}),\, y_1>0 \},\,M_1 \simeq \SL(2,\mathbb R),\\
N_2 =& \{n'(x_2,0), x_2 \in \mathbb R \}, A_2 =\{a(y_2^{1/2},y_2^{-1/2}),\, y_2>0 \},\,M_2 \simeq \SL(2,\mathbb R)
\end{align*}
and
$$
e({\bf P_0})\simeq  \{{\rm pt }\}, \,\,e({\bf P_1}) \simeq  \mathbb H,\, \, e({\bf P_1}) \simeq  \mathbb H
$$
and adding these to $\mathbb D$ get the partial compactification which quotients to our $X.$ \\

In our coordinates, going to $e({\bf P_1})$ is given by $y_1 \rightarrow \infty $ and 
 going to $e({\bf P_2})$ given by $y_2 \rightarrow \infty.$ Hence, $t = \sqrt {y_1y_2},$ as used above (and later on) on several occasions, covers 
 both cases.\\

{\bf Lie algebras and differentials}\\

Later on we shall realize the Kudla-Millson prescription for the construction of the Schwartz form
as a two form on $\mathbb D.$ Hence one needs differentials on $\mathbb D.$ Our coordinates on $\mathbb D$ are 
$z_j = x_j+iy_j, j=1,2,$ resp. $\mu _1=-x_2, \mu _2=x_1, t=\sqrt {y_1y_2}, s=\sqrt {y_1/y_2}$ appearing as parameters in 
 our parabolic subgroup of $\tilde G = \Or(\tilde Q)$
$$
\tilde{P_0}  = \{ n(\mu _1,\mu _2)a(t,s) ; \,\mu _1,\mu_2 \in \R, t,s \in \R_{>0} \}
$$
if
$$
n'(\mu _1,\mu _2) = \begin{pmatrix}
 1 &   \mu_1  &   \mu_2 &   -\mu_1 \mu_2\\ 
 & 1 & &  -\mu_2\\
&  & 1 & -\mu_1 \\
&  &  &1
\end{pmatrix},\quad
a(t,s) =
\begin{pmatrix}
 t &  && \\
&s && \\
& & s^{-1} & \\
&&& t^{-1} 
\end{pmatrix}. 
$$

\nn {\bf Lemma.}\label{leminvdiff} \emph{The left-invariant differentials on $\tilde P_0$ are given by}
\begin{align}
\label{invdiff}
\Omega ^1_\mathbb D =  \langle \nu _1 = dt/t, \nu _2 = ds/s, \nu _3 = d\mu _1/(t/s), \nu _4 = d\mu _2/(ts) \rangle.\nne
\end{align}
{\bf Proof.} This is an easy consequence of the relation
$$
\begin{array}{rcl}
g'g &=& n'(\mu' _1,\mu' _2)a(t',s')n'(\mu _1,\mu _2)a(t,s)\\
&=& n'(\mu' _1,\mu' _2)n'(\mu _1t'/s',\mu _2t's')a(t',s')a(t,s)\\
&=& n'(\mu' _1 + (t'/s')\mu_1,\mu'_2 + t's'\mu _2)a(t't,s's)\\
&=&: n'(\mu'' _1,\mu'' _2)a(t'',s'')
\end{array}
$$
as one has
$$
dt''/t'' = dt/t = \nu _1
$$
etc. \hfill$\Box$\\

\nn {\bf Remark.} With $\tilde\p_0 = \Lie \tilde P_0$ one has
$$
\Omega ^1_\mathbb D \simeq \tilde\p_0^\ast .
$$
Now, in our case, the Kudla-Millson construction of their Schwartz form is based on elements $\omega _{ij},$ 
duals to elements $X_{ij}$ from Lie algebra $\g = \ot(2,2)$ of $G = \Or(2,2)$ 
which are defined as follows. For
$$
\begin{array}{rcl}
\mathfrak{g} &=&\operatorname{Lie}(G)= 
\mathfrak{k} +\mathfrak{p}\\[0.4 cm]
&=& \langle \begin{pmatrix} X&\\&Y \end{pmatrix}: X,Y \mbox{ skew-symmetric} \rangle +
\langle \begin{pmatrix}&Z\\{}^tZ& \end{pmatrix}: Z \in M_{2,2}(\mathbb R)\, \rangle.
\end{array}
$$
as usual, we use the identification
$$
\rho : \wedge ^2V \longrightarrow \ot(V)
$$
for $v,v',v'' \in V$ given by
$$
\rho (v\wedge v') v'' = (v,v'')v' -(v',v'')v.
$$
\nn {\bf Definition.} \label{DefXij} We denote $V = \langle e_1,e_2,e_3,e_4 \rangle$ and put
$$
\begin{array}{lrclrl}
X_{12} &:=& e_1\wedge e_2 = (\begin{smallmatrix}&1&&\\-1&&&\\&&&&\\&&&& \end{smallmatrix}),\quad &
X_{34} &:=& e_3\wedge e_4 = (\begin{smallmatrix}&&&&\\&&&&\\&&&&1\\&&&-1& \end{smallmatrix}),\\
X_{14} &:=& e_1\wedge e_4 = (\begin{smallmatrix}&&&&1\\&&&&\\&&&&\\1&&&& \end{smallmatrix}),\quad &
X_{23} &:=& e_2\wedge e_3 = (\begin{smallmatrix}&&&&\\&&&1&\\&1&&&\\&&&& \end{smallmatrix}),\\
X_{13} &:=& e_1\wedge e_3 = (\begin{smallmatrix}&&&1&\\&&&&\\1&&&&\\&&&& \end{smallmatrix}),\quad &
X_{24} &:=& e_2\wedge e_4 = (\begin{smallmatrix}&&&&\\&&&&1\\&&&&\\&1&&& \end{smallmatrix}).
\end{array}
$$
One has $\mathfrak{p} = \langle X_{14},X_{23},X_{13},X_{24}\rangle$ and hence 
$\mathfrak{p}* = \langle \omega _{14},\omega _{23},\omega _{13},\omega _{24}\rangle.$
Here our aim is a map $\sigma ^\ast : \p^\ast \rightarrow \tilde\p_0^\ast $ and for this we get\\

\nn {\bf Proposition.}\label{propsigomeg} \emph{One has}
$$
\begin{array}{rcl}
\sigma ^\ast (\omega _{13}) &=& (1/2)(-\nu _4 + \nu _3)  = (1/2)(d\mu _1/(t/s) - d\mu _2/(ts))\\
&=& - (1/2)(dx_2/(t/s) + dx_1/(ts))\\
&=& -(1/2)(dx_1/y_1 + dx_2/y_2), \\[0.4 cm]
\sigma ^\ast (\omega _{24}) &=& (1/2)(\nu _4 + \nu _3)  = (1/2)(d\mu _1/(t/s) + d\mu _2/(ts)) \\
&=& (1/2)(-dx_2/(t/s) + dx_1/(ts)) \\
&=& (1/2)(dx_1/y_1 - dx_2/y_2),\\[0.4 cm]
\sigma ^\ast (\omega _{14}) &=& \nu _1  = dt/t = (1/2)(dy_1/y_1 + dy_2/y_2),\\[0.4 cm]
\sigma ^\ast (\omega _{23}) &=& \nu _2  = ds/s = (1/2)(dy_1/y_1 - dy_2/y_2),\\[0.4 cm]
\end{array}
$$
{\bf Proof.} Parallel to $\wedge ^2V \simeq \ot(2,2) = \ot(V) = \mathfrak g = \mathfrak k +\p$ 
one has $\wedge ^2\tilde V \simeq  \ot(\tilde V) = \tilde {\mathfrak g} = \tilde {\mathfrak k} +\tilde {\p}$ 
where $\tilde V = \langle u_1,u_2,u'_2,u'_1\rangle$ with
$$
(u_1,u_2,u'_2,u'_1)= (e_1,e_2,e_3,e_4)C \,\,\mbox{and} \,\, \tilde Q = {}^tCQC .
$$
$$
C := (1/\sqrt 2) \,\,\begin{pmatrix}1&&&1\\&1&1&\\&1&-1&\\1&&&-1 \end{pmatrix}
$$

Matrices $X$ acting on $V$ transform to matrices $\tilde X = {}^tCXC$ acting on $\tilde V.$
We get\\
$$
\begin{array}{lrclrcl}
X_{12} &:=& e_1\wedge e_2 = \left(\begin{smallmatrix}&-1&&\\1&&&\\&&&&\\&&&& \end{smallmatrix}\right)&\longmapsto &
\tilde X_{12} &=& -(1/2)\left(\begin{smallmatrix}&1&1&\\-1&&&-1\\-1&&&-1\\&1&1& \end{smallmatrix}\right),\\[0.4cm]
X_{34} &:=& e_3\wedge e_4 = \left(\begin{smallmatrix}&&&&\\&&&&\\&&&&1\\&&&-1& \end{smallmatrix}\right)&\longmapsto &
\tilde X_{34} &=& (1/2)\left(\begin{smallmatrix}&-1&1&\\1&&&-1\\-1&&&1\\&1&-1& \end{smallmatrix}\right),\\[0.4cm]

X_{14} &:=& e_1\wedge e_4 = \left(\begin{smallmatrix}&&&&1\\&&&&\\&&&&\\1&&&& \end{smallmatrix}\right)&\longmapsto &
\tilde X_{14}  &=& \left(\begin{smallmatrix}1&&&&\\&&&&\\&&&&\\&&&&-1 \end{smallmatrix}\right),\\[0.4cm]
X_{23} &:=& e_2\wedge e_3 = \left(\begin{smallmatrix}&&&&\\&&&1&\\&1&&&\\&&&& \end{smallmatrix}\right)&\longmapsto &
\tilde X_{23} &=& \left(\begin{smallmatrix}&&&\\&1&&\\&&-1&\\&&& \end{smallmatrix}\right),\\[0.4cm]
X_{13} &:=& e_1\wedge e_3 = \left(\begin{smallmatrix}&&&1&\\&&&&\\1&&&&\\&&&& \end{smallmatrix}\right)&\longmapsto &
\tilde X_{13} &=& (1/2) \left(\begin{smallmatrix}&1&-1&&\\1&&&1\\-1&&&-1\\&1&-1& \end{smallmatrix}\right),\\[0.4cm]

X_{24} &:=& e_2\wedge e_4 = \left(\begin{smallmatrix}&&&\\&&&1\\&&&\\&1&& \end{smallmatrix}\right)&\longmapsto &
\tilde X_{24} &=& (1/2)\left(\begin{smallmatrix}&1&1&\\1&&&-1\\1&&&-1&\\&-1&-1& \end{smallmatrix}\right).
\end{array}
$$

Using again as at the beginning the identification
$
\rho : \wedge ^2\tilde V \longrightarrow \ot(\tilde V)
$
for $v,v',v'' \in \tilde V$ given by
$
\rho (v\wedge v') v'' = (v,v'')v' -(v',v'')v,
$
 this time we get as matrices acting on $\tilde V = \langle u_1,u_2,u'_2,u'_1 \rangle$ 
$$
\begin{array}{lcrll}
\tilde U_{14} &:=& u_1\wedge u'_1 &=& \left(\begin{smallmatrix}-1&&&\\&&&\\&&&\\&&&1 \end{smallmatrix}\right) = - \tilde X_{14},\\[0.4cm]
\tilde U_{23} &:=& u_2\wedge u'_2 &=& \left(\begin{smallmatrix}&&&\\&-1&&\\&&1&\\&&& \end{smallmatrix}\right) = -\tilde X_{23},\\[0.4cm]
\tilde U_{12} &:=& u_1\wedge u_2 &=& \left(\begin{smallmatrix}&&-1&\\&&&1\\&&&\\&&& \end{smallmatrix}\right)=
(1/2)((\tilde X_{13}- \tilde X_{24})+(\tilde X_{12}-\tilde X_{34}))\\[0.4cm]
\tilde U_{13} &:=& u_1\wedge u'_2 &=& \left(\begin{smallmatrix}&-1&&\\&&&\\&&&1\\&&&& \end{smallmatrix}\right)=
 (1/2)((\tilde X_{12}+ \tilde X_{34})-(\tilde X_{13}+\tilde X_{24})),\\[0.4cm]
\tilde U_{24} &:=& u_2\wedge u'_1 &=& \left(\begin{smallmatrix}&&&\\-1&&&\\&&&\\&&1& \end{smallmatrix}\right) =
- (1/2)((\tilde X_{13} + \tilde X_{24})+(\tilde X_{12}+\tilde X_{34})),\\[0.4cm]

\tilde U_{34} &:=& u'_2\wedge u'_1 &=& \left(\begin{smallmatrix}&&&\\&&&\\-1&&&\\&1&& \end{smallmatrix}\right)=
 (1/2)((\tilde X_{13}+ \tilde X_{34})-(\tilde X_{12}+\tilde X_{24})).\\[0.4cm]
\end{array}
$$

Using this notation, for $\tilde \p_0 = \Lie \tilde P_0$ we get
$$
\tilde \p_0 = \langle \left(\begin{smallmatrix}1&&&\\&&&\\&&&\\&&&-1 \end{smallmatrix}\right) =: -\tilde U_{14},
\left(\begin{smallmatrix}&&&\\&1&&\\&&-1&\\&&& \end{smallmatrix}\right) =: -\tilde U_{23}, 
\left(\begin{smallmatrix}&1&&\\&&&\\&&&-1\\&&& \end{smallmatrix}\right) =: -\tilde U_{13},
\left(\begin{smallmatrix}&&1&\\&&&-1\\&&&&\\&&&& \end{smallmatrix}\right) =: -\tilde U_{12} \rangle
$$
We have
$$
\tilde \g = \tilde \ka + \tilde \p = \langle \tilde X_{12}, \tilde X_{34} \rangle + 
\langle \tilde X_{14}, \tilde X_{23},\tilde X_{13},\tilde X_{24}\rangle
$$
and, by composing the injection
 $$
\tilde \p_0 \longrightarrow \tilde \g,
 $$
the surjection
 $$
\tilde \g \longrightarrow  \tilde \g/\tilde \ka \simeq \tilde \p,
 $$
 and the identification $\tilde\p \rightarrow \p$
we get the map 
 $$
\sigma:\,\, \tilde \p_0\longrightarrow \p.
 $$
Obviously $\sigma $ maps $\tilde X_{14} = -\tilde U_{14}$ to $X_{14}$ and $\tilde X_{23} =-\tilde U_{23}$ to $X_{23}$ and 
one has
$$
\sigma (\tilde U_{12} ) = (1/2)( X_{13} -  X_{24})
,\,\,\sigma (\tilde U_{13}) = -(1/2)( X_{13} +  X_{24}).
$$
This comes out as follows: One has
$$
\begin{array}{rcl}
\tilde \ka &=&  \langle \tilde X_{12}, \tilde X_{34} \rangle \\
&=&  \langle \tilde X_{12}+ \tilde X_{34}, \tilde X_{12}- \tilde X_{34} \rangle \\
\end{array}
$$
and hence
$$
\begin{array}{rcl}
\tilde U_{12}  &=& (1/2)((\tilde X_{13} - \tilde X_{24}) + (\tilde X_{12} - \tilde X_{34})\\
&\equiv &  (1/2)((\tilde X_{13} - \tilde X_{24}) \mod \tilde \ka\\[0.4 cm]
\tilde U_{13} &=& -(1/2)((\tilde X_{13} + \tilde X_{24})) -(\tilde X_{13} + \tilde X_{24}))\\
&\equiv &  -(1/2)((\tilde X_{13} + \tilde X_{24}) \mod \tilde \ka\\ .
\end{array}
$$
We take the dual $\sigma ^\ast $ of $\sigma $ and get the proposition.\hfill$\Box$\\

\nn {\bf Remark.} In the future we will abuse notation and skip the $\sigma ^\ast .$\\

Hence we get the formulae which will be needed in the construction of the Schwartz forms.\\

\nn {\bf Corollary.} \label{Cordiff}\emph{One has
\begin{align}\label{diff13}
\Omega _1 :&= \omega _{13}\wedge \omega _{14} \nonumber\\
&= -(1/4)(dx_1\wedge dy_1/y_1^2 + dx_2\wedge dy_2/y_2^2 +
(dx_1\wedge dy_2 + dx_2\wedge dy_1)/(y_1y_2))\nne\\
&= -(i/8)(dz_1\wedge d\bar z_1/y_1^2 + dz_2\wedge d\bar z_2/y_2^2 + (dz_1\wedge d\bar z_2 + dz_2\wedge d\bar z_1)/(y_1y_2))\nonumber
\\[0.4 cm]
\Omega _3 :&= \omega _{23}\wedge \omega _{24} \nonumber\\
&= -(1/4)(dx_1\wedge dy_1/y_1^2 + dx_2\wedge dy_2/y_2^2 -
(dx_1\wedge dy_2 + dx_2\wedge dy_1)/(y_1y_2))\nne\\
 &= -(i/8)(dz_1\wedge d\bar z_1/y_1^2 + dz_2\wedge d\bar z_2/y_2^2 
- (dz_1\wedge d\bar z_2 + dz_2\wedge d\bar z_1)/(y_1y_2))\nonumber
\end{align}
and}
\begin{align}\label{diff22}
\Omega _2 :&= \omega _{13}\wedge \omega _{24} + \omega _{23}\wedge \omega _{14}\nonumber\\
 &= (1/2)(dx_1\wedge dx_2 + dy_1\wedge dy_2)/(y_1y_2)\nne\\
&= (1/4)(dz_1\wedge d\bar z_2 - dz_2\wedge d\bar z_1)/(y_1y_2) \nonumber
\end{align}

{\bf Some Derivatives}\\

As we shall use this later, we still stay for a moment with differentials and assemble some 
material which may have some interest on its own. One has
$$
dd^c = - \frac{1}{2\pi i}\partial \bar \partial 
$$
and in our situation

\begin{align}
d' = 4\pi i d^c :=& \partial -\bar \partial \nonumber\\=& (1/2)(\partial _{x_1} - i \partial _{y_1})(dx_1 + idy_1) -
(1/2)(\partial _{x_2} + i \partial _{y_2})(dx_2 - idy_2) \nonumber\\&+
(1/2)(\partial _{x_2} - i \partial _{y_2})(dx_2 + idy_2) -
(1/2)(\partial _{x_1} + i \partial _{y_1})(dx_1 - idy_1)\nonumber\\[0.4 cm]
=& i(\partial _{x_1}dy_1 - \partial _{y_1}dx_1 +\partial _{x_2}dy_2 - \partial _{y_2}dx_2)\nonumber\\[0.4 cm]
=& i(\partial _{\mu _2}(sdt +tds) - (1/2)((1/s)\partial _t +(1/t)\partial _s)d\mu _2 \nne\\
& - \partial _{\mu _1}((1/s)dt - (t/s^2)ds) + (1/2)(s\partial _t - (s^2/t)\partial _s)d\mu _1)\nonumber\\[0.3 cm]
=& i((st\partial _{\mu _2}- (t/s)\partial _{\mu _1})\omega _{14} + (st\partial _{\mu _2}+ (t/s)\partial _{\mu _1})\omega _{23}
+t\partial _t \omega _{13} - s\partial _s \omega _{24})\nonumber
\end{align}
and
\begin{align}
d := \partial +\bar \partial =& 
 (\partial _{x_1}dx_1 + \partial _{y_1}dy_1 +\partial _{x_2}dx_2 + \partial _{y_2}dy_2)\nonumber\\[0.4 cm]
=& (\partial _{\mu _2}d\mu _2 + (1/2)((1/s)\partial _t +(1/t)\partial _s)(sdt +tds) \nonumber\\
& + \partial _{\mu _1}d\mu _1 + (1/2)(s\partial _t - (s^2/t)\partial _s)((1/s)dt - (1/s^2)ds)\nonumber\\[0.4 cm]
=& \partial _{\mu _1}d\mu _1 + \partial _{\mu_2}d\mu _2 + \partial _tdt + \partial _sds\nne\\
=& ((-st\partial _{\mu _2}+ (t/s)\partial _{\mu _1})\omega _{13} + (st\partial _{\mu _2}+ (t/s)\partial _{\mu _1})\omega _{24}
+t\partial _t \omega _{14} + s\partial _s \omega _{23}.\nonumber
\end{align}
By application to special types of functions $\phi = \phi (t,s,x_1,x_2)$ we get\\

\nn {\bf Lemma.} i) \emph{For $\phi = tF(s)$ one has
\begin{align}
\label{d'tFs}
d^c(tF(s))  &= (1/8\pi )(-((1/s)F + F')d\mu _2  + (sF -s^2F')d\mu _1\nonumber\\
 &= (1/4\pi )(tF(s) \omega _{13} - stF'(s) \omega _{24})\nonumber\\
&=(1/8\pi )(-F(dx_1/s+sdx_2)-sF'(dx_1/s-sdx_2)),\nne\\
d (tF(s)) &= F(s)dt + tF'(s) ds\nonumber\\
&= tF(s) \omega _{14} + stF'(s) \omega _{23},\nonumber
\end{align}
and}

\begin{align}
\label{ddtFs}
d d^c (tF(s)) 
 &= (1/8\pi ) (F - sF' -s^2F'')(ds/s \wedge sd\mu _1 + ds/s \wedge d\mu _2/s)\nonumber\\[.2cm]
&= (1/4\pi) t(F - sF' - s^2F'') \omega _{23}\wedge \omega _{24}\nne.
\end{align}

ii)  \emph{For $\phi =F(t)$ one has}

\begin{align}
\label{d'Ft}
d'F(t) &= -(i/2s) F' d\mu _2 + (i/2)sF'd\mu _1\nne\\[.3cm]
dd'F(t) &= - (i/2)(F''dt\wedge (d\mu _2/s-sd\mu _1) - F'ds/s\wedge (d\mu _2/s-sd\mu _1))\nonumber\\[.3cm]
&=i(tF'\omega _{23}\wedge \omega _{24} + t^2F''\omega _{14}\wedge \omega _{13})\nonumber\\[.1cm]
dd^cF(t) &= (1/4\pi )(tF'\omega _{23}\wedge \omega _{24} + t^2F''\omega _{14}\wedge \omega _{13}.\nonumber
\end{align}
\nn {\bf Corollary.}  \emph{For $F(t)=t$ one has
\begin{align}
\label{dd^ct}
dd^c t = (1/4\pi )t \omega _{23}\wedge \omega _{24}\nne.
\end{align}
and for }$F(t)= \log t$
\begin{align}\label{dd^cl}
dd^cF &= (1/4\pi )(\omega _{23}\wedge \omega _{24} - \omega _{14}\wedge \omega _{13}) = (1/4\pi )(\Omega _1+\Omega _2 )\nonumber\\
&= -(i/(16\pi ))(dz_1\wedge d\bar z_1/y_1^2 + dz_2\wedge d\bar z_2/y_2^2)\nne, 
\end{align}
i.e., up to a factor, the {\it K\"ahler form} for our $\mathbb D.$

\section{The Kudla-Millson Schwartz form}

Here we assemble some material round about realizations of the general Kudla-Millson theory of Schwartz forms 
and their associated theta forms.\\

{\bf The $\Or(2,2)-$Schwartz form}\\

 Generalizing the standard Schwartz function $\varphi _0 \in \mathcal S(\mathbb R^4),$ with
$$
\varphi _0(x) = \exp (-\pi (x_1^2+x_2^2+x_3^2+x_4^2))
$$
the general Kudla-Millson prescription (\cite{KM1},\cite{KM2}, in particular \cite{KM3} p.147) proposes as Schwartz form  
 $\varphi _{KM} \in \mathcal S(\mathbb R^4)\otimes \wedge ^2\p,$  
in our situation
\nn {\bf Definition.}\label{DefSchwartzform}
\begin{align}
\label{Schwartzform}
\varphi _{KM}(x) = \varphi _{(2,2)}(x)  = & 2(x_1^2-(1/(4\pi)))\varphi _0(x)\omega _{13}\wedge \omega _{14}\nonumber\\
                               &+  2x_1x_2 \varphi _0(x)(\omega _{13}\wedge \omega _{24}+\omega _{23}\wedge \omega _{14})\nne\\
                               &+  2 (x_2^2-(1/(4\pi)))\varphi _0(x)\omega _{23}\wedge \omega _{24}\nonumber.
\end{align}

The $\omega _{ij}$ are the duals of elements $X_{ij}$ from the Lie algebra $\g = \ot(2,2)$ of $G$ 
which we fixed above in the definition \ref{DefXij}.
With a slight abuse of notation we abbreviate this by
$$
\varphi _{KM} = c_1 \varphi _0(x) \Omega _1 + c_2 \varphi _0(x) \Omega _2 +c_3 \varphi _0(x) \Omega _3.
$$
We  translate to the {\bf a}-coordinates related to the $x-$coordinates by (\ref{a})
$$
x_1 = (1/\sqrt 2)(a+d), \,x_2 = (1/\sqrt 2)(c-b), \, x_3 = - (1/\sqrt 2)(c+b), \,x_4 = (1/\sqrt 2)(a-d).
$$ 
\nn {\bf Remark.} It is easy to see that one has
$$
\begin{array}{rcl}
\varphi _0(x)&=&\exp(- \pi(x_1^2+x_2^2+x_3^2+x_4^2))\\ 
&=&\exp(-\pi  (a^2+b^2+c^2+d^2))
\end{array}
$$
and
$$
\begin{array}{lcl}
c_1 &  = & (a + d)^2 -(1/(2\pi )),\\
c_2 & = & (a+d)(c-b),\\
c_3 & = &  (c - b)^2 - (1/(2\pi )).
\end{array}
$$

From objects living on $V$ we come to objects living on the symmetric space $\mathbb D$ by making the orthogonal group act 
(by left inverse) on the respective coordinates. As we have chosen $(z_1,z_2)$ as parameters for $\mathbb D$ 
(one has to distinguish the components $x_j,y_j$ of these $z$ from the coordinates used for $V$)
by the Key-Remark (\ref{Key})
$$
\varphi _0(A(z)^{-1}{\bf a}) = \exp(-\pi(2R(z,M) + (M,M)))=: {\varphi _0}'
$$
and similarly by the formulae (\ref{a'})
$$
\begin{array}{rcl}
{c_1}' &  = & (a' + d')^2 -(1/(2\pi )),\\
{c_2}' & = & (a'+d')(c'-b'),\\
{c_3}' & = & (c' - b')^2 - (1/(2\pi )).
\end{array}
$$
where (never forget $t^2 = y_1y_2, s^2 = y_1/y_2$) 
\begin{align}
\label{a'+d'}
a'+d' =& \,t^{-1}(a-x_2b-x_1c+x_1x_2d +d t^2)\nonumber,\\
 c'-b' =& \,t^{-1}(y_1(c-x_2d) - y_2(b-x_1d))\nne
\end{align}
We write $R(z,M) =: (1/2y_1y_2)A\bar A$ with
$$
A = (a -bz_2 - cz_1 +dz_1z_2)
$$
hence 
$$
\Rea A = a-x_2b-x_1c+d(x_1x_2-y_1y_2), \quad \Ima A = -by_2 -cy_1 + d(x_1y_2 +x_2y_1).
$$
We introduce
$$
D := y_1(c-x_2d) - y_2(b-x_1d) = \Ima A - 2y_1(dx_2 -c)
$$
and get the expressions
$$
\begin{array}{rcl}
{c_1}' &  = & (t^{-2}(\Rea A + 2dt^2)^2 -(1/(2\pi ))),\\
{c_2}' & = & t^{-2}(\Rea A + 2dt^2)D,\\
{c_3}' & = & (t^{-2} D^2 - (1/(2\pi ))).
\end{array}
$$

\nn {\bf Notation.} Hence, we shall write
\begin{align}
\label{phi'}\nne
\varphi _{KM}(z,M) = ({c_1}'\Omega _1 + {c_2}'\Omega _2 + {c_3}'\Omega _3)\varphi '_0.
\end{align}
Here we use the expressions (\ref{diff13}) and (\ref{diff22}) for the two-forms in our parameters, namely
$$
\begin{array}{rcl}
\Omega _1 &=& (1/8i)(dz_1\wedge d\bar z_1/y_1^2 + dz_2\wedge d\bar z_2/y_2^2 
+ (dz_1\wedge d\bar z_2 + dz_2\wedge d\bar z_1)/(y_1y_2))\\[0.4 cm]
\Omega _2 &=& (1/4)(dz_1\wedge d\bar z_2 - dz_2\wedge d\bar z_1)/(y_1y_2)\\[0.4 cm]
\Omega _3 &=& (1/8i)(dz_1\wedge d\bar z_1/y_1^2 + dz_2\wedge d\bar z_2/y_2^2 
- (dz_1\wedge d\bar z_2 + dz_2\wedge d\bar z_1)/(y_1y_2)).\\
\end{array}
$$
Moreover we introduce the following notation: We write

$$
\begin{array}{rl}
\Omega _1 =&: (1/8i)(\Omega + \Omega ^+)\\[.2cm]
\Omega _3 =&: (1/8i)(\Omega - \Omega ^+)\\[.2cm]
\Omega _2 =&: (1/4) \Omega ^-
\end{array}
$$
and
\begin{align}
\label{phiz}
\varphi _{KM}(x) =& ((c_1+c_3)(1/8i)\Omega +(c_1-c_3)(1/8i)\Omega ^+ + c_2(1/4)\Omega ^-)\varphi _0 \nne\\[.3cm]
=&: (c_{11}\frac{dz_1\wedge d\bar z_1}{y_1^2} + c_{22}\frac{dz_2\wedge d\bar z_2}{y_2^2} + 
c_{12}\frac{dz_1\wedge d\bar z_2}{y_1y_2} + c_{21}\frac{dz_2\wedge d\bar z_1}{y_1y_2})\varphi _0 \nonumber
\end{align}
with
$$
\begin{array}{rl}
c_{11} =& (1/8i)((a+d)^2 + (c-b)^2 -(1/\pi ))\\[.2cm]
c_{22} =& (1/8i)((a+d)^2 + (c-b)^2 -(1/\pi ))\\ [.2cm]
c_{12} =& (1/8i)((a+d)^2 - (c-b)^2 + 2i(a+d)(c-b))\\ [.2cm]
c_{21} =& (1/8i)((a+d)^2 - (c-b)^2 - 2i(a+d)(c-b)).\\[.2cm]
\end{array}
$$
For $\varphi _{KM}(z,M)$ we write as above but replace $a$ by $a'$ etc., i.e., 
these formulae above represent the special value of $\varphi _{KM}(z,M)$ for $z = (z_1,z_2)$ with $z_1 = z_2 = i.$\\

By a small calculation for the transformed coefficients we get
the nice expression
$$
c'_{11}  = c'_{22} = (1/(8i))(2R(z,M) + 2(M,M) -1/\pi )
$$
while the other coefficients do not look that nice. For eventual further use, we state
\begin{align}
\label{Omega}
\frac{dz_1\wedge d\bar{z}_1 }{y_1^2} + \frac{dz_2\wedge d\bar{z}_2 }{y_2^2} =& \,4i(\Omega _1 + \Omega _3) = \Omega \nonumber\\
\frac{dz_1\wedge d\bar{z}_2 }{y_1y_2} =&  \,2i(\Omega _1-\Omega _3 + 2\Omega _2) \nonumber\\
\frac{dz_2\wedge d\bar{z}_1 }{y_1y_2} =& \, 2i(\Omega _1-\Omega _3 - 2\Omega _2) \nne.
\end{align}


\nn {\bf Remark.} The Schwartz form is a real form and as such it looks like this
\begin{align}
\label{phireal}
\varphi _{KM}(z,M) = - (1/4)\varphi '_0&\big(((a'+d')^2 + (c'-b)^2 - (1/\pi ))(\frac{dx_1\wedge dy_1}{y_1^2}+\frac{dx_2\wedge dy_2}{y_2^2})\nonumber   \\
\!+&((a'+d')^2 - (c'-b')^2)\frac{dx_1\wedge dy_2-dy_1\wedge dx_2}{y_1y_2}\nonumber\\
 -&2(a'+d')(c'-b')\frac{dx_1\wedge dx_2+dy_1\wedge dy_2}{y_1y_2}\big)\nne
\end{align}
And we get immediately for the special case $M = 0$ in our Schwartz form $\varphi _{KM},$ i.e.
$$
\varphi _{KM}(z,0) = -(1/2\pi )(\Omega _1 + \Omega _3) = (i/8\pi )(\frac{dz_1\wedge d\bar z_1}{y_1^2} + \frac{dz_2\wedge d\bar z_2}{y_2^2}),
$$ 
ans by (\ref{dd^cl}) one has 
\begin{align}
\label{dd^clog}\nne
dd^c \log (t^2) = - \varphi _{KM}(z,0).
\end{align}
As clear from the discussion of generating functions and modularity in our main part, one has to comply with a symplectic variable $\tau =u+iv\in \HH.$ 
Via the Weil representation (we will return to this further on), one has the standard way to replace the variables $a,b,c,d$ by $\sqrt v a$ etc. and 
hence we introduce notations for Schwartz functions and forms as
\begin{align}
\label{Sch}
\varphi _0(v,z,M) =& \,e^{-2\pi vR(z,M)}\nonumber\\
\varphi _\infty (\tau ,z,M) =& \,e^{\pi i(M,M)_{\tau ,z}} \nonumber\\
=& \,e^{\pi i(u(M,M) + iv(M,M)_{z})} = \,e^{\pi i(u(M,M) + iv(2\det M+2R(z,M))}\nonumber\\
=& \,e^{-2\pi vR(z,M)}e^{2\pi i\tau \det M},\nne\\
\varphi _{KM}(v,z,M) =& \,\sum c'_j \Omega _j\varphi _0(v,z,M)\nonumber.
\end{align} 
\nn {\bf Remark.} Finally, we state the nice formulae
\begin{align}
\label{2phi} 
\varphi _{KM}(v,z,M)\wedge& \varphi _{KM}(v,z,M) \nonumber\\
=& (1/16\pi )(v(R(z,M)+(M,M))-1/(4\pi)) \varphi _0(v,z,M)^2 d\mu   \nne 
\end{align}
with
$$
d\mu := \frac{dz_1\wedge d\bar z_1\wedge dz_2\wedge d\bar z_2}{(y_1y_2)^2},
$$
and
$$
\varphi _{KM}(v,z,M)\wedge \Omega = (1/(2i))(v(R(z,M) + 2\det M) - 1/(2\pi ))\varphi _0(v,z,M)d\mu .
$$\\

{\bf The Schwartz form and its relation to the $\Or(1,2)-$case}\label{subs:Schwartzform(1,2)}\\

One has a natural map (which we shall fix soon)
$$
\Or(1,2)/\Or(2) \simeq \mathbb D_1 \longrightarrow \mathbb D = \Or(2,2)/\Or(2)\times \Or(2)
$$
and, hence we can restrict our $\Or(2,2)-$Schwartz form $\varphi _{KM}=\varphi _{(2,2)}$ on $\mathbb D$ to $\mathbb D_1.$ By the general Kudla-Millson 
theory, as for instance stated in Theorem 2.1 in [Fu], one expects to get the well known and 
elaborated $\Or(1,2)-$Schwartz form $\varphi _{(1,2)}.$ Though, in principle, this is superfluous, we make this explicit in a
very elementary way, as it presents a nice exercise. We want to give sense to and prove the following\\

\nn {\bf Proposition.}\label{prop:(2,2)(1,2)} \emph{For 
$$
M = \begin{pmatrix}  a&b\\c&d \end{pmatrix},\,\,X := \begin{pmatrix} \tilde a&b\\c&-\tilde a \end{pmatrix}
$$
with $\tilde{a} = (1/2)(a-d),\,\tilde{d} = (1/2)(a+d)$ one has}

\begin{align}
\varphi _{KM}(v,(\bar w,-1/\bar w),M) e^{2\pi i\tau \det M} = e^{2\pi i\tau \tilde d^2}\varphi _{(1,2)}(\tau ,w,X)\nne.
\label{phiw}
\end{align}

{\bf Proof : 1.} The $\Or(1,2)-$case as discussed for instance in [Fu] and \cite{BF} uses a vector space 
$$
V_0 = \{ X \in M_2(\R), {\rm tr} X = 0 \}
$$ 
with quadratic form $q(X) = \det X = -\tilde a^2 - bc$  for 
$$
X := \begin{pmatrix} \tilde a&b\\c&-\tilde a \end{pmatrix}
$$
and coordinates $z =x+iy \in \mathbb H$  stemming from the action of $\SL_2(\R)$ on $V$ 
given by $X \mapsto \gamma X \gamma ^{-1}.$ Hence the associated Schwartz form is (\cite{BF} (3.7))
$$
\varphi (\tau ,z,X) := \varphi _{(1,2)}(\tau ,z,X) = (v(X,X(z))^2 -(1/2\pi )) e^{\pi i(X,X)_{\tau ,z}} \omega 
$$
where
$$
\begin{array}{rl}
(X,X(z)) =& -(1/y)(cz\bar z-2\tilde ax-b)\\[.3cm]
(X,X)_{\tau ,z} =& u(X,X) + iv(X,X)_z \\[.3cm]
(X,X)_z =& 2\det X + 2R(z,X)  \\[.3cm]
R(z,X) =& (1/2)(\,X,X(z))^2 - 2\det X  \\[.3cm]
\omega =& dx\wedge dy/y^2.
\end{array} 
$${\bf 2.} If we want to see how this is related to our Schwartz form $\varphi _{KM},$ we have to realize it in the 
alternative {\it new} coordinates $w_1,w_2$ from (\ref{newcoord}) 
related to the old coordinates by 
$$
z_1 = \bar w_1, \,\, z_2 = - 1/\bar w.
$$
Then the special cycle $T_1 \simeq \mathbb D_1$ is given by $w_1 = w_2 = w.$ From above we have 
$$
\varphi (\tau ,z,M) = \varphi _0(v,z,M)e^{2\pi i\det M\tau } = e^{\pi i(M,M)_{\tau ,z}},
$$
$$
(M,M)_{\tau ,z} = u(M,M) + iv(M,M)_z,\,\, (M,M)_z = 2R(z,M)+(M,M).
$$
with 
$$
(M,M) = 2(\tilde d^2 + \det X)
$$
and introducing the new coordinates and putting  $w_1=w_2=w$ we get as in (\ref{RW})
$$
R(z,M) \longmapsto \frac{|aw_2 + b - cw_1w_2 - dw_1|^2}{2v_1v_2} =: \tilde R((w_1,w_2),M)
$$
and hence
$$
2\tilde R((w,w),M) + 2\det M = 2\tilde d^2 + 2R(w,X) + 2\det X,
$$
i.e. if we specialize the Schwartz function of the $\Or(2,2)-$case to $T_1,$ as predicted by Theorem 2.1 and Theorem 5.1 in [Fu],  we get  up to 
the theta function generating factor the Schwartz function of the $\Or(1,2)-$case. 
$$
\varphi _0(\tau ,(\bar w,-(1/\bar w)),M) = e^{2\pi i\tau \tilde d^2} \varphi _0(\tau ,w,X).
$$
{\bf 3.} In order to determine the behaviour of the coefficients $c_j$ of the Schwartz form, again we look at the map  
$$
\tilde G = \SL(2,\R)^2 \longrightarrow \Or(2,2)
$$
induced by the action $M \mapsto \gamma _1M\gamma _2^{-1},$ namely 
$\gamma := (\gamma _1,\gamma_2) \longmapsto A(\gamma)$ as in (\ref{A-1}) for
$\gamma = (g_{w_1},g_{w_2})$ with
$$
g_{w_j} = \begin{pmatrix} v_j^{1/2} & u_jv_j^{-1/2}\\ &v_j^{-1/2} \end{pmatrix} 
$$
one has with $\tilde t =(v_1v_2)^{1/2},\,\tilde s = (v_1/v_2)^{1/2}$
$$
A(w_1,w_2)^{-1} = A(\gamma )^{-1}
 = \begin{pmatrix} \tilde s^{-1}&0&-u_1\tilde s^{-1}&0\\
u_2\tilde t^{-1}&\tilde t^{-1}&u_1u_2\tilde t^{-1}&-u_1\tilde t^{-1}\\
0&0&\tilde t&0\\
0&0&u_2\tilde s&\tilde s\\
\end{pmatrix}. 
$$
We get 
$$
{\bf a'} = A(w_1,w_2)^{-1}{\bf a} 
$$
with
$$
\begin{array}{rl}
a' =& \tilde s^{-1}a - \tilde s^{-1}u_1c,\\
b' =& \tilde t^{-1}(u_2a + b -u_1u_2c - u_1d),\\
c' =& \tilde tc,\\
d' =& \tilde s(u_2c + d).
\end{array}
$$
A simple calculation specializing $w_1=w_2=w$ recovers the result already obtained above
$$
(M,M)_{(w,w)} = \ta{\bf a'}\cdot {\bf a'} \equiv 2\tilde d^2 + (X,X)_. 
$$

{\bf 4.} But to get control over the behaviour of the Schwartz form, 
one has to take the description obtained at the end of the last section, namely (\ref{phiz})
$$
\begin{array}{rl}
\varphi _{KM} =& ((c_1+c_3)(1/8i)\Omega +(c_1-c_3)(1/8i)\Omega ^+ + c_2(1/4)\Omega ^-)\varphi _0\\[.3cm]
=&: (c_{11}\frac{dz_1\wedge d\bar z_1}{y_1^2} + c_{22}\frac{dz_2\wedge d\bar z_2}{y_2^2} + 
c_{12}\frac{dz_1\wedge d\bar z_2}{y_1y_2} + c_{21}\frac{dz_2\wedge d\bar z_1}{y_1y_2})\varphi _0 
\end{array}
$$
with
$$
\begin{array}{rl}
c_{11} =& (1/8i)((a+d)^2 + (c-b)^2 -(1/\pi ))\\[.1cm]
c_{22} =& (1/8i)((a+d)^2 + (c-b)^2 -(1/\pi ))\\ [.1cm]
c_{12} =& (1/8i)((a+d)^2 - (c-b)^2 + 2i(a+d)(c-b))\\ [.1cm]
c_{21} =& (1/8i)((a+d)^2 - (c-b)^2 - 2i(a+d)(c-b)).\\
\end{array}
$$
\\[.3cm]
and analyze the change of the differentials. From 
$$
z_1 = \bar w_1,\,\,z_2 = -1/\bar w_2
$$
we get
$$
y_1 = -v_1,\,\, y_2 = -v_2/|w_2|^2, \,\,dz_1 = d\bar w_1,\,\,dz_2 = d\bar w /\bar w_2^2
$$
and hence

\begin{align}\label{Omegaw}
\Omega _{11} =& \frac{dz_1\wedge d\bar z_1}{y_1^2} =  - \frac{dw_1\wedge d\bar w_1}{v_1^2}\nonumber\\[.3cm]
\Omega _{22} =& \frac{dz_2\wedge d\bar z_2}{y_2^2} =  - \frac{dw_2\wedge d\bar w_2}{v_2^2}\nonumber\\[.3cm]
\Omega _{12} =& \frac{dz_1\wedge d\bar z_2}{y_1y_2} =  - \frac{dw_2\wedge d\bar w_1}{v_1v_2}\frac{|w_2|^2}{w_2^2}\nonumber\\[.3cm]
\Omega _{21} =& \frac{dz_2\wedge d\bar z_1}{y_1y_2} =  - \frac{dw_1\wedge d\bar w_2}{v_1v_2}\frac{|w_2|^2}{\bar w_2^2}.\nne
\end{align}

For  $w_1=w_2=w$ at $w = i$ one has
$$
\begin{array}{rl}
\Omega _{11} \,=& \Omega _{22} \equiv - {dw\wedge d\bar w}\\[.3cm]
\Omega _{12} \,=& \Omega _{21} \equiv \,\,dw\wedge d\bar w.\\[.3cm]
\end{array}
$$
and hence
$$
\begin{array}{rll}
\varphi _{KM}(v,(i,i),M)e^{2\pi i\tau \det M} \equiv & - (1/(8i))&(c_{11} + c_{22} - c_{12} - c_{21}) \varphi (\tau ,(i,i),M)\, dw\wedge d\bar w \\[.3cm]
\equiv & - (1/(2i))&((c-b)^2-(1/(2\pi )) e^{2\pi i\tau \tilde d^2} \varphi _0(\tau ,i,X) \,dw\wedge d\bar w .\\
\end{array}
$$
Thus we have shown that  for the special value $w_1=w_2= i$ the Schwartz form of the $\Or(2,2)-$case 
specializes up to the factor $e^{-2\pi \tilde d^2}$ to the 
Schwartz form of the $\Or(1,2)-$case and this is sufficient to have (\ref{phiw})
$$
\varphi _{KM}(v,(\bar w,-(1/\bar w),M) e^{2\pi i\tau \det M} = e^{2\pi i\tau \tilde d^2}\varphi _{(1,2)}(\tau ,w,X).
$$
\hfill$\Box$

{\bf Kudla's Green function and the Schwartz form}\\

With $- \operatorname{Ei}(-t) := \int\limits_1^\infty \exp(-tr)\frac{dr}{r}$  Kudla in \cite{Ku1} p.330  defines
$$
\xi (z,M) := - \operatorname{Ei}\left( - 2 \pi R(z_1,z_2,M)  \right) \exp(-\pi(M,M))
 $$
and states (as his Proposition 4.10)\\
\nn {\bf Proposition.} \label{PropKUDLA}\emph{One has as currents on $\mathbb D$
\begin{align}
\label{KUDLA}\nne
 {dd}^c \xi (z,M) + \exp(-\pi(M,M))\delta _{D_M} = [\varphi _{KM}(z,M)]
 \end{align}
 with} $D_M = \{ z \,: R(z,M) = 0 \}.$ \\

 {\bf Proof :} As the relation above is fundamental and Kudla only states the existence of a proof, we compute away from the singularities  
\begin{align*}
- \operatorname{dd}^c  \operatorname{Ei}\left(
- 2 \pi R(z_1,z_2,M)  \right)  &=
\frac{1}{2\pi i}\partial \bar \partial   \operatorname{Ei}\left(
 - 2 \pi R(z_1,z_2,M)  \right) \\
&= \frac{1}{2 \pi i} \partial \left( e^{-2\pi R} 
\bar \partial \log(R) \right)\\
&= \frac{1}{2 \pi i} \left( -2 \pi \partial R \wedge \bar\partial \log (R)
+ \partial\bar\partial \log(R) \right) e^{-2\pi R}. 
\end{align*}
Using
$$
\partial y^{-1} = -(1/2i) y^{-2} dz,\quad \bar \partial y^{-1} = (1/2i) y^{-2} d\bar z
$$
the second term is computed easily
$$
\begin{array}{rcl}
\partial\bar\partial \log(R) = - \partial\bar\partial \log(y_1y_2) &=& (i/2)(\frac{dx_1\wedge dy_1}{y_1^2} + 
\frac{dx_2\wedge  dy_2}{y_2^2})\\
&=& (1/4)(\frac{dz_1\wedge d\bar z_1}{y_1^2} + \frac{dz_2\wedge  d\bar z_2}{y_2^2}) = (1/4)\Omega .
\end{array}
$$
The first term is much more involved. We use again
$$
A := (a - bz_2 - cz_1 + dz_1z_2)
$$
such that $ R = A\bar A/(2y_1y_2)$ and moreover
$$
dB := y_1^{-1}dz_1 + y_2^{-1}dz_2
$$
such that
$$
\partial R = (i/2)R \,dB + (1/(2y_1y_2)) \bar A\,\partial A
$$
and
$$
\bar \partial \log R = -(i/2) d\bar B +(1/\bar A) \bar \partial \bar A
$$
Hence we get
$$
\begin{array}{rcl}
\partial R \wedge \bar\partial \log(R)&=&
((i/2)R dB + (1/(2y_1y_2)) \bar A\partial A)\wedge (-(i/2) d\bar B +(1/\bar A) \bar \partial \bar A)\\[0.4 cm]
&=& (1/(2y_1y_2))[(1/4) A\bar A dB \wedge d\bar B +(i/2) (AdB\wedge \bar \partial \bar A - \bar A\partial A\wedge d\bar B)
+\partial A\wedge \bar \partial \bar A ]\\
\end{array}
$$
and
$$
\begin{array}{rcl}
dd^c \xi 
&=& (i/(2y_1y_2))e^{-2(\pi R +\det M)}\\[0.4 cm]
&{}&\,\times [(1/4)A\bar A(dz_1\wedge d\bar z_1/y_1^2 + dz_2\wedge d\bar z_2/y_2^2 +
(dz_1\wedge d\bar z_2 +dz_2\wedge d\bar z_1)/y_1y_2)\\[0.4 cm]
&{}& \,\,\,\,+(i/2)(A((\bar z_2d-c)/y_1 dz_1\wedge d\bar z_1 +(\bar z_1d-b)/y_2 dz_2\wedge d\bar z_2\\
& {}&\quad \quad \quad \,\,\,\,\,\,+(\bar z_1d-b)/y_1 dz_1\wedge d\bar z_2 +(\bar z_2d-c)/y_2 dz_2\wedge d\bar z_1))\\
&{}& \quad \quad\quad\,-(\bar A(( z_2d-c)/y_1 dz_1\wedge d\bar z_1 +( z_1d-b)/y_2 dz_2\wedge d\bar z_2 \\
&{}&\quad\quad \quad\quad+( z_2d-b)/y_2 dz_1\wedge d\bar z_2 +(z_1d-c)/y_1 dz_2\wedge d\bar z_1))\\[0.4 cm]
&{}& \,\,\,+(\mid z_2d-c\mid ^2 dz_1\wedge d\bar z_1 + \mid z_1d-b\mid ^2 dz_2\wedge d\bar z_2 \\
&{}&\quad +(z_2d-c)(\bar z_1d-b) dz_1\wedge d\bar z_2 +(z_1d-b)(\bar z_2d-c) dz_2\wedge d\bar z_1) \\[0.4 cm]
&{}&\,\,\,-(1/4\pi )y_1y_2(dz_1\wedge d\bar z_1/y_1^2 + dz_2\wedge d\bar z_2/y_2^2 )]
\end{array}
$$
Again we write this as
$$
\begin{array}{rl}
\varphi 
=&: (\check c_{11}\frac{dz_1\wedge d\bar z_1}{y_1^2} + \check c_{22}\frac{dz_2\wedge d\bar z_2}{y_2^2} + 
\check c_{12}\frac{dz_1\wedge d\bar z_2}{y_1y_2} + \check c_{21}\frac{dz_2\wedge d\bar z_1}{y_1y_2})\varphi '_0 
\end{array}
$$

Though it is rather cumbersome, let us look at the coefficient $\check c_{11}$ of $\varphi '_0 dz_1\wedge d\bar z_1/y_1^2.$
We get
$$
i(1/(2y_1y_2)[A\bar A/4 + (i/2)(y_1A((\bar z_2d-c) - y_1\bar A(z_2d-c)) + y_1^2\mid z_2d - c \mid ^2) - (2y_1y_2)/(8\pi v)]
$$
Some easy calculation (one only has to verify this for $z_1 = z_2 = i$) leads to
$$
\check c_{11} = (i/8)((2R(z,M) + 2(M,M)) - 1/\pi) .
$$
and as well
$$
\check c_{22} = (i/8)((2R(z,M) + 2(M,M)) - 1/\pi) .
$$
For the coefficient $\check c_{12}$ of $(1/(y_1y_2)) \varphi '_0 dz_1\wedge d\bar z_2$ we get
$$
\check c_{12} = (i/2)(A\bar A/(4y_1y_2) + (i/2)(A((\bar z_2d-b)/y_1 -\bar A( z_2d-c)/y_2) + ( z_2d-c)(\bar z_2d-b))
$$
which with a similar computation for $z_1 = z_2 = i$ comes out as
$$
\check c_{12} = (i/8)((a+d)^2 - (c-b)^2 + 2i(a+d)(c-b))
$$
We leave the treatment of the last coefficients to the reader and 
see that one has the equality of $dd^c \xi $ and $-\varphi _{KM}$ with $\varphi _{KM}$ from (\ref{phiz}). Here, one has that the Schwartz form 
of Kudla's proposition is the negative of our $\varphi _{KM}$ built from the prescription of \cite{KM1} (The minus sign would disappear 
if we change to our alternative coordinates (\ref{newcoord}).\hfill$\Box$\\

{\bf Extension to the symplectic variable}\\

In the main part of our text and also later on here, we discuss generating functions 
and their modularity. To prepare the ground for this we describe the standard way to introduce a symplectic variable.\\

The general procedure to introduce in a situation near to our a symplectic variable goes like this: 
One starts with a Schwartz function $\varphi \in \Ss(V)$ with $\dim V = m, \sig V = (p,q)$ and variables 
written as $x, {\bf a},$ or $M.$ For $\tilde G = \Or(V) \supset \tilde P_0 \simeq \mathbb D = G/K$ with $A(z) \in \tilde P_0$ for $z \in \mathbb H^2$ 
(in our case $\mathbb D \simeq \mathbb H^2$), by the usual procedure, one defines
$$
\varphi (z,x) := \varphi (A(z)^{-1} x).
$$
resp. in our notation with coordinates as in (\ref{a}) $ \varphi (z,M) := \varphi (A(z)^{-1} {\bf a}).$ Moreover, using the 
mechanism of the Weil representation $\omega $ of $G' = \SL(2,,\R)$ on $\Ss(V),$ for $\tau = u+iv \in \mathbb H$ 
one defines 
$$
\varphi (\tau ,x) := j({g'}_\tau ,i)^\ell \omega ({g'}_\tau)\varphi (x)
$$
if $\varphi $ has weight $\ell$ for $K' = \SO(2)$. 
Hence, in our case we come to
\begin{align}
\label{taux}\nne
\varphi (\tau ,x) = v^{-\ell /2 +m/4} \varphi (\sqrt v x) e^{\pi i (x,x)u}.
\end{align}
We can use Kudla's statement (from \cite{Ku1} p.329): for $p=n, q=2,$
$$
\varphi _\infty (x,z) = \varphi _\infty (z,M)= e^{-\pi (x,x)_z} = e^{-2\pi R(z,x)} e^{- \pi (x,x)}
$$ 
has weight $\ell = (n/2)-1.$ \\

\nn {\bf Remark.} As already indicated after the definition \ref{DefKugr}, Kudla 
prefers to consider the Green function weighted by an exponential factor
\begin{align}
\label{xi}\nne
\xi (z,M) := - \operatorname{Ei}\left( - 2 \pi R(z_1,z_2,M)  \right) \exp(-\pi(M,M)).
\end{align}
Following Kudla, we treat $\xi$ as if it were a Schwartz function.
In our context we have $\ell =0$ and, hence, our 
$$
\xi (z,M) =  \int_1^\infty e^{-2\pi R(z,M)r}dr/r \,e^{-2\pi \det M}
$$
extends to 
\begin{align}
\label{tauzM}\nne
 \xi (\tau ,z,M) = ( \int_1^\infty e^{-2\pi vR(z,M)r}dr/r ) \,e^{2\pi i \tau \det M }.
\end{align}
An inspection of the proof of the Proposition in the last subsection shows 
that equation (\ref{KUDLA}) 
generalizes to a here central statement.\\

\nn {\bf Corollary.} \label{Cor:dd^cxitau} \emph{Outside of the singularities one has}
\begin{align}
\label{kudla}\nne
dd^c\xi (\tau ,z,M) = - \varphi _{KM}(\tau ,z,M).
\end{align}
Similarly, using [FM] Prop.4.20 with $m = 2, \ell = 2,$ we could extend the {\it boundary function} $\check \xi $ from (\ref{boundary function})
which will reappear instantly by the considerations of the next section and get
\begin{align}
\label{tauts}\nne
\check \xi (\tau ,(t,s),(b,c)) =  (t/\sqrt v)(B(v,s;b,c) - I(v,s;b.c))e^{2\pi i(-bc)\,\tau }.
\end{align}

\nn {\bf Remark.} In our alternative coordinates (\ref{newcoord}) one has the Kudla Green function 
$\tilde{\xi}(w,M) = \int_1^\infty e^{-2\pi \tilde{R}(w,M)r}dr/r$ with $\tilde{R}(w,M)$ as in (\ref{RW}) 
and, adapting the proof of Proposition \ref{PropKUDLA} and using the formulae (\ref{Omegaw}), this time (ouside 
of the singularities) we get  
$$
dd^c \tilde{\xi}(w,M) = \tilde{\varphi }_{KM}(w,M).
$$

\newpage

\section[Funke-Millson's restriction to the boundary]{Funke-Millson's restriction \\to (a neighbourhood of) the boundary}

We want to know what happens to the Green function and the Schwartz form at or near to the boundary of the orthogonal space $\mathbb D.$ 
In sections 2 and  3. we already discussed the approach from an ancient draft by Ulf K\"uhn, but now we want to check 
what comes out by transferring to our case the general theory from the paper \cite{FM1} by Funke and Millson.\\

{\bf A mixed model and the restriction of $\varphi _{KM}$} \\

We take Kudla-Millson's Schwartz form

$$
\begin{array}{rl}
\varphi _{KM} = \varphi _{(2,2)}  = & 2(x_1^2-(1/4\pi))\varphi _0(x)\omega _{13}\wedge \omega _{14}\\
                               &+  2x_1x_2 \varphi _0(x)(\omega _{13}\wedge \omega _{24}+\omega _{23}\wedge \omega _{14})\\
                               &+  2 (x_2^2-(1/4\pi))\varphi _0(x)\omega _{23}\wedge \omega _{24}\\[0.3 cm]
              =& c_1 \varphi _0(x) \Omega _1 + c_2 \varphi _0(x) \Omega _2 +c_3 \varphi _0(x) \Omega _3,\\[0.3 cm]
                  \varphi _0(x) =& \exp (-\pi (x_1^2+x_2^2+x_3^2+x_4^2)) ,            
\end{array}
$$
via the Weil representation, make it live on the space $\mathbb D$ and want to understand what happens 
if one goes near to the boundary of $\mathbb D.$ We follow Funke-Millson and introduce the 
{\it mixed model} of the Weil representation: We take the totally isotropic subspace 
$$
E = \,\,<u_1 = (1/\sqrt 2)(e_1 + e_4)>\,\, \subset  V = \,\,<e_1,e_2,e_3,e_4>,
$$ 
its dual 
$$
E' = \,\,<u'_1 = (1/\sqrt 2)(e_1 - e_4)>,
$$ 
and
$$
W = \,<u_2 = (1/\sqrt 2)(e_2 + e_3), u'_2 = (1/\sqrt 2)(e_2 - e_3)>.
$$
Then one has
$$
V \simeq E \oplus W \oplus E' = \tilde V
$$
where we denote the coordinates by $y = (y_1,y_2,y_3,y_4)$ 
and write with a slight abuse $\varphi (y)$ for the function corresponding to 
$\varphi (x) \in \mathcal S(V).$ Now, we perform the partial Fourier transformation
$\varphi \longmapsto \hat \varphi $ given by
$$
\hat \varphi (\zeta ,y_2,y_3,y_4) := \int_{-\infty }^\infty \varphi(y) e^{-2\pi iy_1\zeta }dy_1.
$$  
\nn {\bf Remark.} The {\it mixed model} of the Weil representation has as its 
representation space 
the space of these functions. The formulae for the action of $G' = \SL(2,\R)$ and 
$G = {\rm O}\,(\tilde V)$ are assembled in \cite{FM1} Lemma 4.1. We can recover these quickly. 
Using (\ref{a'})
$$
y' = A'(z)^{-1}y = \begin{pmatrix}
 t^{-1}(y_1 + x_2y_2 - x_1y_3 + x_1x_2y_4)\\ 
 s^{-1}(y_2 + x_1y_4)\\
s(y_3 - x_2y_4) \\
ty_4
\end{pmatrix},
$$
 we have 
$$
\begin{array}{rl}
n(x_1,x_2)a(t,s)\varphi (y) = \varphi (A(z)^{-1}y) = \varphi (y') =: \varphi '(y)
\end{array}
$$
and get
$$
\hat {\varphi '}(\zeta ,y_2,y_3,y_4) = te^{2\pi i\zeta (x_2y_2-x_1y_3+x_1x_2y_4)}\,\hat \varphi (t\zeta , s^{-1}(y_2 + x_1y_4),s(y_3 - x_2y_4),ty_4).
$$

We want to apply this to the Schwartz form and, as preparation, state the following formulae which easily follow 
by application of $\partial _\zeta $ to the first relation.\\

\nn {\bf Lemma.}\label{Fourierlem} \emph{We have}
$$
\begin{array}{rclrcl}
h(u) :&=& e^{-\pi u^2} = \varphi _0(u) \quad\Longrightarrow& \hat h(\zeta ) &=& \varphi _0(\zeta ),\\
h(u) :&=& u \varphi _0(u) \quad& \hat h(\zeta ) &=& -i\zeta \varphi _0(\zeta ),\\
h(u) :&=& u^2 \varphi _0(u) \quad& \hat h(\zeta ) &=& -(\zeta ^2 -(1/2\pi))\varphi _0(\zeta ),\\
h(u) :&=& (u^2 -(1/2\pi))\varphi _0(u) \quad& \hat h(\zeta ) &=& -\zeta ^2 \varphi _0(\zeta ),\\
h(u) :&=& (u^2 -(1/4\pi))\varphi _0(u) \quad& \hat h(\zeta ) &=& -(\zeta ^2 - (1/4\pi))\varphi _0(\zeta ),\\
h(u) :&=& \varphi(t^{-1}u) \quad& \hat h(\zeta ) &=& t \varphi (t\zeta ).\\
\end{array}
$$
At first, we take the Schwartz function
$$
\varphi _0(y) = e^{-\pi(y_1^2+y_2^2+y_3^2+y_4^2)}.
$$
We put $\varphi '_0(y) := \varphi _0(y')$ and get with $Y := x_2y_2-x_1y_3+x_1x_2y_4$
$$
\widehat {\varphi '_0}(\zeta ,y_2,y_3,y_4) = te^{-\pi t^2\zeta ^2}e^{-\pi (s^{-2}(y_2+x_1y_4)^2+s^2(y_3-x_2y_4)^2+t^2y_4^2)}e^{2\pi i\zeta Y}
$$

Moreover, one has 

$$
\begin{array}{rl}
\widehat {\varphi ' _{KM}}(\zeta ,y_2,y_3,y_4) =& (t^3(-i\zeta +y_4)^2 \Omega _1\\
& + t^2(-i\zeta +y_4)(s^{-1}(y_2+x_1y_4) + s(y_3 -x_2y_4)) \Omega _2\\
& + t((s^{-1}(y_2+x_1y_4) + s(y_3 -x_2y_4))^2 -(1/2\pi)) \Omega _3)\\
&\times e^{-\pi t^2\zeta ^2}e^{-\pi (s^{-2}(y_2+x_1y_4)^2+s^2(y_3-x_2y_4)^2+t^2y_4^2)}e^{2\pi i\zeta Y}.
\end{array}
$$
Here in each term the $\varphi _0-$factor goes to zero for $t \longmapsto \infty $ 
if one does not have $\zeta = y_4 = 0.$ In this case we are left with
$$
\begin{array}{rl}
\widehat {\varphi ' _{KM}}(0,y_2,y_3,0) =
& t((s^{-1}y_2 +sy_3)^2 -(1/2\pi))e^{-\pi (s^{-2}y_2^2+s^2y_3^2)} \,\omega _{23}\wedge \omega _{24}\\[.3cm]
=& ((s^{-1}y_2 +sy_3)^2 -(1/2\pi))e^{-\pi (s^{-2}y_2^2+s^2y_3^2)}\, ds/s \wedge ( dx_1/s - sdx_2 )
\end{array}
$$

This is (up to a factor) Funke-Millson's form $r_P \varphi _{KM}$ and we shall call it 
{\it the restricted Schwartz form} and also write it in our ${\bf a}-$coordinates as
$$
\check \varphi _{KM}(z;b,c) = \,t((b/s-cs)^2 - (1/(2\pi ))))e^{-\pi ((b/s)^2+(cs)^2)} \omega _{23}\wedge \omega _{24}.
$$
If, as done above several times, we introduce the symplectic variable $\tau = u+iv,$ 
we get\\

\nn {\bf Definition.} \label{DefresSchw} For $t\longrightarrow \infty $ one has as {\it restricted Schwartz form}
\begin{align}
\label{checkphi}\nne
\check \varphi _{KM}(\tau ,z;b,c) = t(v(b/s-cs)^2 - (1/(2\pi ))))e^{-\pi (v((b/s)-(cs))^2)}e^{-2\pi i\tau bc} \Omega _3.
\end{align}

This naming explains as via Poisson summation we get\\

\nn {\bf Proposition.} \label{PropresSchw} \emph{For $t\longrightarrow \infty $ one has }
\begin{align}
\label{phicheckphi}\nne
\sum _M \varphi _{KM}(\tau ,z,M) = \sum _M \hat{\varphi} _{KM}(\tau,z,M) \longrightarrow \sum_{b,c} \check \varphi _{KM}(\tau ,z;b,c).
\end{align}
The smoothness of the partition function $\rho $ immediately leads to\\

\nn{\bf Corollary.} \label{Corphirho} \emph{$\varphi _\rho (v,z,m)$ is a smooth function on $X.$ }\\

\vspace{.3cm}
{\bf Mixed model and restriction of the Green function}\\

 We try the same reasoning for the Green function and take (here already with the silent variable $v$ and slight abuse of notation)
$$
\xi (v,z,M) = \xi (v,z,a,b,c,d) = - \,{\rm Ei}\,(- 2\pi vR) =  \int_1^\infty e^{- 2\pi vR(z,M) r} dr/r
$$
with
$$
\begin{array}{rl}
R =&  ((1/(2y_1y_2))\mid a - bz_2 - cz_1 + dz_1z_2 \mid ^2\\[.3cm]
 =&  ((1/(2y_1y_2)) ((a-bx_2 -cx_1 + d(x_1x_2 - y_1y_2))^2 + (-by_2-cy_1 + d(x_1y_2+x_2y_1))^2)
\end{array} 
$$
For $t = \sqrt{y_1y_2} \longrightarrow \infty $ the integrand goes to zero if not $d = 0.$ 
Hence, assume $d = 0.$ With ${\bf x} :=bx_2+cx_1, {\bf y} := |by_2+cy_1|$ we get
$$
\xi (v,z;a,b,c,0) =   \int_1^\infty e^{- (\pi vr/t^2)((a-{\bf x})^2 +  {\bf y} ^2)} dr/r.
$$
We perform the Fourier transform concerning the variable $a$ and get
$$
\begin{array}{rl}
\hat \xi (v,z;\zeta ,b,c,0) =&  \int_{-\infty} ^\infty \int_1^\infty e^{- (\pi vr/t^2)((a-{\bf x})^2 +  {\bf y} ^2) - i2\pi a\zeta } (dr/r)da\\[.3cm]
=&  \,(t/\sqrt v) e^{-2i\pi {\bf x}\zeta } \int_1^\infty e^{- \pi (t^2\zeta ^2/(rv) + {\bf y} ^2 vr/t^2)} r^{-3/2}dr.
\end{array}
$$
Here we meet again the expression we already discussed in our boundary considerations, namely in calculating (\ref{hbc1}). 
For $\zeta \not=0$ the integrand goes to zero and for $\zeta =0$ we get
$$
\hat \xi (v,z;0,b,c,0) = (t/\sqrt v)\int_1^\infty e^{- \pi v( {\bf y} ^2 r/t^2)} r^{-3/2}dr.
$$
Using the Poisson formula we recover the first term of our old boundary function (\ref{boundary function})
$$
\check \xi (v,z;b,c) = (t/\sqrt v)(B(v\alpha (s) - I(v\alpha (s))),
$$ 
which makes sense for all $b,c\in \ZZ.$ For $-bc\geq 0$ by multiplication with $e^{-2\pi v\det M}$ we get
$$
\begin{array}{rl}
\check \xi (z;0,b,c,0)\,\cdot e^{2\pi vbc} =& (t/\sqrt v)  \int_1^\infty e^{- \pi v({\bf y}^2 u/t^2)} u^{-3/2}du \,\cdot e^{2\pi vbc}\\[.3cm]
=& \,(t/\sqrt v)  \int_1^\infty e^{- \pi  v(b/s + cs)^2 u} u^{-3/2}du \,\cdot e^{2\pi vbc}\\[.3cm]
=& \,(t/\sqrt v)B(\pi v(bs^{-1}+cs)^2)\,\cdot e^{2\pi vbc}
\end{array}
$$

For the moment, we shall call this {\it the restricted Green function} and denote it by $\check \xi '(v,z;b,c).$ \\

\nn {\bf Definition.} In its symplectic version one has the {\it  restricted Green function}
\begin{align}
\label{checkxi'}\nne
\check \xi '(\tau ,z;b,c) = (t/\sqrt v)B(\pi v((b/s)+(cs))^2)\,\cdot e^{2\pi i\tau (-bc)}.
\end{align}

{\bf $dd^c$ of the  restricted Green function}\\

 We relate the restrictions of our two objects.
\nn {\bf Proposition.} \label{PropComDiagr} \emph{If we associate $\check \xi '$ to the Green function $\xi $ 
and $\check \varphi _{KM}$ to the Schwartz form 
$\varphi _{KM},$ we get a commuting diagram: the equation (\ref{kudla})
$$
dd^c \xi = - \varphi _{KM}
$$
translates to }
\begin{align}
\label{ddccheckxi'}\nne
dd^c \check \xi '=  - \check \varphi _{KM}.
\end{align}
{\bf Proof}: In the section on More Derivatives we obtained the relation (\ref{ddtFs})
$$
\begin{array}{rl}
d d^c (tF(s))  =& (1/8\pi ) (F - sF' -s^2F'')(ds/s \wedge sd\mu _1 + ds/s \wedge d\mu _2/s)\\[.3cm]
=& (1/4\pi) t(F - sF' - s^2F'') \omega _{23}\wedge \omega _{24}.
\end{array}
$$
We put
$$
\alpha (s) := \pi (b/s + cs)^2
$$
and 
$$
F(s) := B(\alpha (s)), \, B(\alpha ) :=  \int^\infty _1 e^{-\alpha r} r^{-3/2}dr.
$$
and moreover as an auxiliary function
$$
F_0(s) := B_0(\alpha (s)),\,\,\,B_0(\alpha ) := \int^\infty _1 e^{-\alpha r} r^{-1/2}dr.
$$
By partial integration get for $\alpha >0$
$$
B(\alpha ) =  2(e^{-\alpha } - \alpha B_0(\alpha ))
$$
One has
$$
F'(s) = - \alpha '(s) \int^\infty _1 e^{-\alpha r} r^{-1/2}dr = -  \alpha '(s) B_0(\alpha (s)) = - \alpha '(s) F_0(s).
$$
and again by partial integration
$$
F'_0(s) = -\alpha '(s) \int^\infty _1 e^{-\alpha r} r^{1/2}dr = - (\alpha '/\alpha )(e^{-\alpha} + (1/2)F_0).
$$
From here we get
$$
F''(s) = - \alpha ''F_0 - \alpha 'F'_0 = ((\alpha ')^2/\alpha )e^{-\alpha } - ((\alpha ')^2/(2\alpha )) -\alpha '')F_0.
$$
and
$$
F - sF' - s^2F'' = 
 (2 - s^2(\alpha ')^2/\alpha )e^{-\alpha } - (2\alpha - s\alpha ' - s^2(\alpha '' + (\alpha ')^2/(2\alpha ))F_0.
$$
From
$$
\alpha = \pi (b/s + cs)^2 =: (\pi /2)(u+v)^2
$$
we get 
$$
s\alpha ' = 2\pi (v^2 - u^2)
$$ 
and
$$
s^2\alpha '' = 2\pi (v^2 + 3u^2).
$$
Hence, we have
$$
s^2(\alpha ')^2/\alpha = 4\pi (v-u)^2
$$
and 
$$
F - sF' - s^2F'' = - ( - 2 + 4\pi(s/b - cs)^2) e^{- \pi (b/s + cs)^2 },
$$
i.e.
$$
\begin{array}{rl} 
d d^c (tF(s)) =& (1/4\pi) t(F - sF' - s^2F'') \omega _{23}\wedge \omega _{24}\\[.3cm]

=&  - t( (b/s - cs)^2 - (1/(2\pi)))e^{- \pi (b/s + cs)^2 }\omega _{23}\wedge \omega _{24}\\[.3cm]
 =&  - ((b/s - cs)^2 - (1/(2\pi)))e^{- \pi (b/s + cs)^2 }(1/2)(ds/s \wedge sd\mu _1 + ds/s \wedge d\mu _2/s)\\[.3cm]
\end{array}
$$
By multiplication with $e^{2\pi bc}$ we get $ - \check \varphi $ and the formula in the Proposition.\hfill$\Box$\\

\nn{\bf Lemma.} \emph{One has}
$$
dd^c (tI(v,s;(b,c))) = 0.
$$
{\bf Proof.} As we know from {\bf 2.} and the main part of our paper, 
for $bc<0$ the restricted Green funcion $\check \xi $ from above still has to be modified 
(to get a smooth function with weight 2 for the Weil representation of $\SL(2,\R)$). As above we take
$$
\begin{array}{rll} 
I(u,v) :=& \min (\mid u \mid , \mid v \mid )\,\,&{\rm if}\,\,uv > 0\\
:=&\,\,0\,\,&{\rm if}\,\,uv \leq  0.
\end{array}
$$
and put
$$
I(s,u,v) := I(u/s,vs)
$$

i.e. one has inside the cone $u^2 - v^2 > 0$
$$
\begin{array}{rl}
I(s,u,v) =&  vs\,\,\,\,{\rm if}\,\,u > 0,\,v > 0, \\[.3cm]
:=& - vs \,\,\,\,{\rm if}\,\,u < 0,\,v < 0, \\
\end{array}
$$
and outside the cone, i.e. for $u^2 - v^2 < 0$
$$
\begin{array}{rl}
I(s,u,v) =& u/s \,\,\,\,{\rm if}\,\,u > 0, v > 0,\\[.3cm]
 :=& - u/s \,\,\,\,{\rm if}\,\,u < 0, v < 0,\\
\end{array}
$$
We see that for $u>v>0$ we have 
$$
\partial _sI = v , \,\, \partial ^2_sI = 0
$$
and for $v>u>0$
$$
\partial _sI = - u/s^2,\,\, \partial ^2_sI = 2u/s^3
$$
Hence one has 
$$
I(s) -sI'(s) - s^2 I''(s) = 0
$$
in both cases and as well, with a similar reasoning,  for $u<0,v<0.$\\

Using again the formulae from 
the subsection on derivatives, we come to
$$
\begin{array}{rl}
d(tI) =& tI \omega _{14} - tI'\omega _{23}\\
d'(tI) =& i(tI \omega _{13} + tI'\omega _{24})\\
dd^c(tI) =& 0.
\end{array}
$$\hfill$\Box$\\
Hence, we can add any multiple of $tIe^{2\pi bc}$ to the function $\,tB(\alpha (s))e^{2\pi bc} $ discussed above and still retain the 
relation $dd^c \check \xi =   \check \varphi .$ To get a smooth function (as function in $b$ and $c$), 
we have to proceed as we did in the main part of the paper.\\

{\bf A variant of the restricction of $\varphi _{KM}$}\\

There is another way to treat the restriction 
to the boundary of our Schwartz form, namely one also can follow 
the way we used in the proof of Theorem \ref{ThnearD} in section 2. Again we start with the Schwartz form
\begin{align}
\label{Sch}\nne
\varphi _{KM}(z,M) = \sum _{j=1}^3 c'_j\Omega _j \varphi _0(z,M)
\end{align}
with

\begin{align}
\label{c'}
\varphi _0(z,M) =& e^{- 2\pi v(R(z,M)+\det M)}, \,\,R(z,M) = (1/2t^2)|a-z_2b-z_1c+z_1z_2d|^2\nonumber\\
c'_1 =& v(a' + d')^2 - 1/(2\pi )\nonumber\\
c'_2 =& v(a' + d')(c' - b')\nne\\
c'_3 =& v(c' - b')^2 - 1/(2\pi )\nonumber\\
a'+d' =& t^{-1}(a-x_2b-cx_1+dx_1x_2+dt^2)\nonumber\\
c'-b' =& s(c-x_2d) - s^{-1}(b-x_1d)\nonumber\\
=& t^{-1}(y_1c-y_2b-d(y_1x_2-y_2x_1))\nonumber
\end{align}
For $t \longrightarrow \infty $ one has
$$
\varphi _0(z,M) \longrightarrow 0
$$
if not $d\not=0.$ Hence we look at $d=0$ and put
$$
{\bf x} := x_2b+x_1c,\,{\bf y} := y_2b+y_1c,\,{\bf y'} := y_1c - y_2b.
$$
We get
$$
\varphi _{KM}(z,a,b,c,0) = \sum _{j=1}^3 c'_j\Omega _j \varphi _0(z,a,b,c,0)
$$
with
$$
\varphi _0(z,a,b,c,0) = e^{- \pi v t^{-2}((a-{\bf x})^2 + {\bf y}^2)}e^{2\pi vbc}
$$
and
$$
a'+d' = t^{-1}(a-{\bf x}),\,c'-b' = t^{-1} {\bf y'}.
$$
Now, for $b$ and $c$ not both zero, the sum
$$
\phi _{bc}(z) := \sum_a c'_j\Omega _j\varphi _0(z,a,b,c,0)
$$
is a periodic function in ${\bf x} $ and has a Fourier expansion
$$
\phi _{bc}(z) = \sum_n a_{bc}(n)e(n{\bf x})
$$
with coefficients
$$
a_{bc}(n) = \int _{-\infty }^\infty \sum _{j=1}^3 \,c'_j\Omega _j \varphi _0(z,0,b,c,0))e(-n{\bf x}))d{\bf x}.
$$
Here we can use our formulae for the Fourier transforms from the beginning of this section and with ${\bf x}= ut/\sqrt v$ obtain
$$
\begin{array}{rl}
a_{1,bc}(\zeta ) =& \int _{-\infty }^\infty\,(vt^{-2}{\bf x}^2-1/(2\pi ))e^{-(\pi /t^2)v({\bf x}^2 
+ {\bf y}^2) - 2\pi i{\bf x}\zeta }e^{2\pi vbc}d{\bf x}\\[.3cm]
=& - (t/\sqrt v)\zeta ^2e^{-\pi t^2\zeta ^2/v}e^{-\pi v {\bf y}^2/t^2}e^{2\pi vbc},\\[.3cm]
a_{2,bc}(\zeta ) =& \int _{-\infty }^\infty\,(vt^{-2}{\bf x}{\bf y'}e^{-(\pi /t^2)v({\bf x}^2 
+ {\bf y}^2) - 2\pi i{\bf x}\zeta }e^{2\pi vbc}d{\bf x}\\[.3cm]
=& - (\sqrt v /t)i\zeta {\bf y'}e^{-\pi t^2\zeta ^2/v}e^{-\pi v {\bf y}^2/t^2}e^{2\pi vbc},\\[.3cm]
a_{3,bc}(\zeta ) =& \int _{-\infty }^\infty\,(vt^{-2}{\bf y'}^2-1/(2\pi ))e^{-(\pi /t^2)v({\bf x}^2 
+ {\bf y}^2) - 2\pi i{\bf x}\zeta }e^{2\pi vbc}d{\bf x}\\[.3cm]
=&  (vt^{-2}{\bf y'}^2-1/(2\pi ))(t/\sqrt v)e^{-\pi t^2\zeta ^2/v}e^{-\pi v {\bf y}^2/t^2}e^{2\pi vbc}.
\end{array}
$$
We see that, for $t\longrightarrow \infty $ and $j=1,2,3$ $a_{j,bc}(\zeta ) \longrightarrow 0$ if $\zeta \not=0$ and for $\zeta = 0$ 
we only are left with
$$
\begin{array}{rl}
a_{3,bc}(0) =&  (vt^{-2}{\bf y'}^2-1/(2\pi ))(t/\sqrt v)e^{-\pi v {\bf y}^2/t^2}e^{2\pi vbc}\\
=&  (vt^{-2}{\bf y'}^2-1/(2\pi ))(t/\sqrt v)e^{-\pi v ((b/s)^2+(cs)^2)}
\end{array}
$$
Hence, for $t\longrightarrow \infty $ one has
$$
\begin{array}{rl}
\sum_{M\in L} &\varphi _{KM}(z,M) \longrightarrow \\[.3cm]
&\sum_{b,c\in\mathbb Z} (v(b/s-cs)^2-1/(2\pi ))(t/\sqrt v)e^{-\pi v ((b/s)^2 + (cs)^2)}\Omega _3\\[.3cm]
 &+ \sum_a\,((v(a/t)^2-1/(2\pi )) e^{-\pi v(a/t)^2}\Omega _1 - 1/(2\pi )e^{-\pi v(a/t)^2}\Omega _3)\\[.3cm]
=& \sum_{b,c\in\mathbb Z} (v(b/s-cs)^2-1/(2\pi ))(1/\sqrt v)e^{-\pi v ((b/s)^2 + (cs)^2)}ds/s\wedge (dx_1/s-sdx_2)\\[.2cm]
 &+ \sum_a\,((v(a/t)^2-1/(2\pi )) e^{-\pi v(a/t)^2}\Omega _1 - 1/(2\pi )e^{-\pi v(a/t)^2}\Omega _3)\end{array}
$$
where in the sum $b$ and $c$ should not both be zero.\\

If one applies Poisson to the sum over $a,$ one has
$$
\sum_a\,((v(a/t)^2-1/(2\pi )) e^{-\pi v(a/t)^2} = \sum_a - (t/\sqrt v)^3a^2 e^{-\pi (t^2/v)a^2}
$$
and
$$
\sum_a\,(-1/(2\pi )) e^{-\pi v(a/t)^2} = \sum_a - (1/2\pi )(t/\sqrt v)e^{-\pi (t^2/v)a^2}.
$$
And, if again we use the argument that for $t \rightarrow \infty $ and $a \not= 0$ things go to zero, 
we get the same result as at the beginning of this section in (\ref{checkphi}) I.e., we have (again considering the symplectic variable)
\begin{align}
\label{phicheck}\nne
\check{\varphi } (\tau ,z) = \sum_{b,c\in\mathbb Z} (v(b/s-cs)^2-1/(2\pi ))(t/\sqrt v)e^{-\pi v ((b/s)^2 + (cs)^2)}e^{-2\pi iubc}\Omega _3 . 
\end{align}
In the next section we will show that this form is $\SL(2,\mathbb Z)-$modular of weight 2.\\

\section{Toy thetas and their modularity}\label{Toy}

At last, we collect and expand some classic material concerning theta functions adapted to our situation.\\

{\bf The Zero value of Jacobi's Theta and its Maass derivative }\\

As for instance in \cite{Mu} one has the classic theta series
$$
\theta (\tau ) = \sum_{n\in\mathbb Z} e^{\pi in^2\tau }
$$
with
$$
\theta (\tau ) = \theta (\tau + 2), \,\, \theta (\tau ) = (i/\tau )^{1/2} \theta (-1/\tau ),
$$
i.e., a modular form of weight $1/2$ for a certain subgroup of $\SL(2,\ZZ),$ the theta group. 
In \cite{Ma} Maass introduced certain {\it weight raising} resp. {\it lowering} differential operators. 
As a special case for weight $1/2$ here we take
$$
X_+ := 2iv\partial _\tau +1/2
$$
and set
$$
\begin{array}{rl}
 \theta _{KM}(\tau ) :=& X_+ \theta (\tau ) \\[.3cm]
=& \sum_n(-2\pi vn^2 + 1/2) e^{\pi in^2\tau } = -2\pi \sum_n(vn^2 - (1/(4\pi )))  e^{\pi in^2\tau } 
\end{array}
$$
We call this $\theta _{KM}$ because one gets a factor near to those in the definition of the Kudla-Millson 
Schwartz form.\\

\nn {\bf Proposition.} \label{Propthetakm}\emph{One has the relations}
\begin{align}
\label{thetakm}\nne
\theta _{KM}(\tau ) = \theta _{KM}(\tau +2), \,\, 
\theta _{KM}(\tau )= - (i/\tau )^{3/2}(i/\bar \tau )^{-1}\theta _{KM}(-1/\tau ).
\end{align}

{\bf Proof.} 
The first relation is clear and the second relation can be seen by the 
Poisson summation formula or by applying $X_+$ to the second
functional equation of $\theta :$ We have
$$
\theta _{KM}(-1/\tau ) = \sum_n (n^2v/\mid \tau \mid^2 - (1/(4\pi)))e^{\pi in^2(-1/\tau )}
$$
and
$$
\begin{array}{rl}
X_+ \,((\tau /i)^{-1/2} &\sum_n e^{\pi in^2(-1/\tau )})\\[.3cm]
=& i^{1/2}\sum_n (4iv(-1/2)\tau ^{-3/2} + \tau ^{-1/2} + 4iv\tau ^{-1/2}\pi in^2/\tau ^2)e^{\pi in^2(-1/\tau )}\\[.3cm]
=& (\tau /i)^{-1/2} \sum_n ( -2iv\tau ^{-1} + 1 - 4\pi vn^2/\tau ^2) e^{\pi in^2(-1/\tau )}\\[.3cm]
=& (\tau /i)^{-1/2} \sum_n ( \bar \tau /\tau  - 4\pi vn^2/\tau ^2) e^{\pi in^2(-1/\tau )}\\[.3cm]
=& 4\pi(\tau /i)^{-1/2} \bar \tau /\tau \sum_n ( 1/(4\pi) -  vn^2/\mid \tau \mid^2) e^{\pi in^2(-1/\tau )}\\[.3cm]
=& - 4\pi (\tau /i)^{-1/2} \bar \tau /\tau \,\,\theta (-1/\tau ).
\end{array}
$$
\hfill$\Box$\\
{\bf A signature (1,1) Siegel theta.}\\

To get a nearer approach to the theta forms growing out of the Kudla-Millson Schwartz forms, 
we go to Siegel's work on Theta series for indefinite quadraic forms \cite{S1} and \cite{S2} and look at the special signature (1,1) Siegel theta
\begin{align}
\label{Siegel theta}\nne
\theta _{(1,1)}(\tau ) := \sum_{b,c}\,E(b,c;\tau ),\,\,E(b,c;\tau ) := e^{-\pi v(b^2+c^2)}e^{-2\pi iubc}
\end{align}
By specializing Siegel's original procedure we get\\

\nn {\bf Proposition.} \label{Proptheta1}\emph{One has the relations
\begin{align}
\label{theta1}\nne
\theta _{(1,1)}(\tau ) = \theta _{(1,1)}(\tau +1),\,\,\theta _{(1,1)}(\tau ) = (1/\sqrt {\bar \tau })(1/\sqrt {\tau })\theta _{(1,1)}(-1/\tau ),
 \end{align}
i.e., $\theta _{(1,1)}(\tau )$ is modular of weight} (1/2,1/2).\\

{\bf Proof.}  
$\theta _{(1,1)}(\tau )$ clearly stays invariant under $\tau \mapsto 1+\tau.$ 
For $\tau \mapsto  -(1/\tau )$ we use the Poisson formula 
$$
\theta _{(1,1)}(\tau ) = \sum_{b,c}\,E(b,c;\tau ) = \sum_{\zeta ,\eta }\,\hat{E}(\zeta ,\eta ;\tau )
$$
and, with $b=(1/\sqrt 2)(x+y),\,c=(1/\sqrt 2)(x-y),$ compute
$$
\begin{array}{rl}
\hat{E}(\zeta ,\eta ;\tau ) &= \int\int  e^{-\pi v(b^2+c^2)}e^{-2\pi iubc} e^{-2\pi i(b\zeta +c\eta )}dbdc\\[.3cm]
&= \int\int  e^{-\pi v(x^2+y^2)}e^{-\pi iu(x^2-y^2} e^{-2\pi i((x+y)\zeta +(x-y)\eta )/\sqrt 2}dxdy\\[.3cm]
&= \int  e^{-\pi x^2(v+iu)}e^{-2\pi ix(\zeta +\eta )/\sqrt 2}dx \int e^{-\pi y^2(v-iu)}e^{-2\pi iy(\zeta -\eta )/\sqrt 2}dy\\[.3cm]
&= \int  e^{-\pi x^2i\bar{\tau } }e^{-2\pi ix(\zeta +\eta )/\sqrt 2}dx \int e^{-\pi y^2(-i\tau )}e^{-2\pi iy(\zeta -\eta )/\sqrt 2}dy\\[.3cm]
\end{array}
$$
moreover, with $w:=x\sqrt {i\bar\tau }$ resp.  $\tilde{w}:=y\sqrt {-i\tau },$
$$
\begin{array}{rl}
\hat{E}(\zeta ,\eta ;\tau ) 
=& \int  e^{-\pi w^2 }e^{-2\pi w(\zeta +\eta )(1/\sqrt 2)(1/\sqrt {i\bar \tau })}dw(1/\sqrt {i\bar \tau })\\
 &\,\times \,\int e^{-\pi \tilde{w} }e^{-2\pi i\tilde w(\zeta -\eta )(1/\sqrt 2)(i/\sqrt {i\tau })}d\tilde{w} (i/\sqrt {i\tau })\\[.3cm]
 =& (1/\sqrt {i\bar \tau })(i/\sqrt {i\tau })e^{-\pi (\zeta +\eta )^2/(2i\bar \tau )}e^{-\pi (\zeta -\eta )^2/(2i \tau )}\\[.3cm]
 =& (1/\sqrt {\bar \tau })(1/\sqrt {\tau })e^{-\pi ((\zeta +\eta )^2 (-i)(u+iv) - (\zeta -\eta )^2(-i)(u-iv))/(2|\tau|^2) }\\[.3cm]
 =& (1/\sqrt {\bar \tau })(1/\sqrt {\tau })e^{-\pi iu2\zeta \eta  + v(\zeta ^2+\eta ^2)^2 )/|\tau|^2 }\\[.3cm]
 =& \,\,(1/\sqrt {\bar \tau })(1/\sqrt {\tau })E(\zeta ,\eta ; -1/\tau ).
\end{array}
$$
\hfill$\Box$\\
{\bf The Maass derivative of the signature (1,1) Siegel theta.}\\

Again we take the Maass operator
$$
X_+ = 2iv\partial _\tau  + 1/2 = iv(\partial _u - i\partial _v) + 1/2
$$
and set 
\begin{align}
\label{theta2}\nne
\theta _{(1,1)}^{KM}(\tau ) := (X_+\theta _{(1,1)}(\tau ) = -\pi \sum_{b,c}\,(v(b-c)^2 - (1/2\pi ))E(b,c;\tau ) .
\end{align}
Here we get\\

\nn{\bf Proposition.} \label{Proptheta22}\emph{One has the relations
\begin{align}
\label{theta22}\nne
\theta _{(1,1)}^{KM}(\tau )= \theta _{(1,1)}^{KM}(\tau +1),\,\,
\theta _{(1,1)}^{KM}(\tau ) = \tau ^{-3/2}\bar{\tau }^{1/2}\theta _{(1,1)}^{KM}(-1/\tau ),
\end{align}
i.e., $\theta _{(1,1)}^{KM}(\tau )$ is of weight }(3/2,-1/2).\\

{\bf Proof.} Application of $X_+$ to the transformation formula (\ref{theta1}) leads to 
$$
\begin{array}{rl}
X_+ \theta _{(1,1)}(\tau ) &= X_+((1/\sqrt {\bar \tau })(1/\sqrt {\tau })\theta _{(1,1)}(-1/\tau ))\\[.3cm]
&= (2iv(-1/2)\tau ^{-3/2}\bar{\tau }^{-1/2} + (1/2)\tau ^{-1/2}\bar{\tau }^{-1/2})\theta _{(1,1)}(-1/\tau ) \\
&\,\,\quad \,\,\,\quad\quad \,\,\,\quad\quad \,\,\,\quad\quad \,\,\,\quad\quad \,\,\,\quad\quad \,
+ 2iv \tau ^{-1/2}\bar{\tau }^{-1/2}\partial _\tau \theta _{(1,1)}(-1/\tau )\tau ^{-2} \\[.3cm]
&= (-iv/\tau +(1/2)\tau /\tau )\tau ^{-1/2}\bar{\tau }^{-1/2})\theta _{(1,1)}(-1/\tau ) \\
&\,\,\quad \,\,\,\quad\quad \,\,\,\quad\quad \,\,\,\quad\quad \,\,\,\quad\quad \,\,\,\quad\quad \,
- \pi v\tau ^{-1/2}\bar{\tau }^{-1/2})\sum (b-c)^2E(b,c;-1/\tau )\tau ^{-2}\\[.3cm]
&= - \pi\tau ^{-3/2}\bar{\tau }^{1/2}\sum_{b,c}\, (v/|\tau |^2(b-c)^2 - (1/2\pi ))E(b,c;-1/\tau )\\[.3cm]
&= \tau ^{-3/2}\bar{\tau }^{1/2}X_+\theta _{(1,1)}(-1/\tau ).
\end{array}
$$
\hfill$\Box$\\

{\bf The Theta coming from our restricted Schwartz form}\\

Our restricted Schwartz form $\check{\varphi }_{KM}(v,t,s;b,c)$ is a two-form 
and for $t=s=1$ has as its coefficient
$$ 
(1/\sqrt v )(v(b-c)^2 - (1/2\pi ))e^{-\pi v(b+c)^2}.
$$
Now, we set
\begin{align}
\label{theta3}\nne
\check{\theta } _{(1,1)}^{KM}(\tau ) :=& (1/\sqrt v )\theta _{(1,1)}^{KM}(\tau )\nonumber\\
=& (1/\sqrt v )\sum_{b,c}(v(b-c)^2 - (1/2\pi ))e^{-\pi v(b^2+c^2)}e^{-2\pi iubc}.
\end{align}

and get\\

\nn {\bf Corollary.\label{Cortheta33}} \emph{One has the relation
\begin{align}
\label{theta33}\nne
\check{\theta } _{(1,1)}^{KM}(\tau ) = \check{\theta } _{(1,1)}^{KM}(\tau +1) = \check{\theta } _{(1,1)}^{KM}(-1/\tau ),
\end{align}
i.e., $\check{\theta } _{(1,1)}^{KM}(\tau )$ is modular of weight 2, and 
hence, as well,}
$$
\theta _{\check{\varphi }_{KM}}(\tau ,z) := \sum_{b,c} \check{\varphi }_{KM}(v,t,s;b,c)e^{2\pi ibc\tau }.
$$ 
{\bf Proof.} For $f(\tau ) = 1/\sqrt v$ one has $f(-1/\tau )=|\tau |/\sqrt v.$ Hence the 
relations in the corollary follow from (\ref{theta22}). And one has
$$
\theta _{\check{\varphi }_{KM}}(\tau ,(i,i)) = \check{\theta } _{(1,1)}^{KM}(\tau )
$$
where the modularity in $\tau $ translates from $(i,i)$ to any $z.$
\hfill$\Box$\\

{\bf A signature (2,2) Siegel Theta.} \\

As this works so nicely, let's look at the thetas coming from the 
four-dimensional space. Certainly, the same way as above, we can treat the signature (2,2) Siegel theta
\begin{align}
\label{Siegel Theta}\nne
\Theta (\tau ,z) := \sum_{M\in M_2(\mathbb Z)} E(M;\tau ,z), \,\,E(M;\tau ,z) = e^{\pi i(u(M,M)+iv(M,M)_z)},
\end{align}
resp. 
$$
\Theta (\tau ) := \Theta (\tau ,(i,i)) = \sum_{M\in M_2(\mathbb Z)}E(M;\tau ,(i,i))
$$
with $x_1=(a+d)/\sqrt 2, x_2=(c-b)/\sqrt 2, x_3=-(c+b)/\sqrt 2, x_4=(a-d)/\sqrt 2$
$$
\begin{array}{rl}
E(M;\tau ) := E(M;\tau ,(i,i)) =& e^{\pi i(u(M,M)+iv(M,M)_{(i,i)})}\\[.3cm]
=& e^{2\pi iu\det M - \pi v(a^2+b^2+c^2+d^2)}\\[.3cm] 
=& e^{\pi iu(x_1^2+x_2^2-x_3^2-x_4^2) - \pi v(x_1^2+x_2^2+x_3^2+x_4^2)}\\[.3cm] 
=& e^{\pi i\tau x_1^2}e^{\pi i\tau x_2^2}e^{-\pi i\bar\tau x_3^2}e^{-\pi i\bar\tau x_4^2}\\[.3cm] 
\end{array}
$$
Convergence again is obvious and we get\\

\nn {\bf Proposition.} \label{PropSiegel}\emph{As predicted by Siegel,  we see that $\Theta $ as function of $\tau $ is modular of weight }(1,1).\\

{\bf Proof.} One has the Fourier transform
$$
\begin{array}{rl}
\hat{E} (\zeta ,\eta ,\xi ,\varrho ;\tau ) :=& \int E(a,b,c,d;\tau )e^{-2\pi i(a\zeta +b\eta +c\xi +d\varrho )}dadbdcdd\\[.3cm]
=& \int e^{\pi i\tau x_1^2 - 2\pi ix_1(\zeta +\varrho )/\sqrt 2}dx_1
\int e^{\pi i\tau x_2^2 - 2\pi ix_2(\xi -\eta  )/\sqrt 2}dx_2\\
&\times \int e^{-\pi i\bar\tau x_3^2 + 2\pi ix_3(\xi +\eta )/\sqrt 2}dx_3
\int e^{-\pi i\bar\tau x_4^2 - 2\pi ix_4(\zeta -\varrho )/\sqrt 2}dx_4\\[.3cm]
=& (1/|\tau |^2) e^{-\pi i/2(((\zeta +\varrho )^2+(\xi -\eta )^2)/\tau - ((\xi +\eta )^2+(\zeta -\varrho )^2)/\bar\tau )}\\[.3cm]
=& (1/|\tau |^2) e^{-2\pi iu(\eta \xi -\zeta \varrho ) - \pi v(\zeta ^2+\eta ^2+\xi ^2+\varrho ^2)/|\tau |^2}\\[.3cm]
=& (1/|\tau |^2){E} (\zeta ,\eta ,\xi ,\varrho ;-1/\tau ) .
\end{array}
$$
and by Poisson summation we come to the hypothesis.
\hfill$\Box$\\

{\bf The Theta coming from our $\varphi _{KM}$.}\\

The Siegel theta $\Theta $ we just dicussed is a kind of kernel to the (1,1)-form theta series 
growing out of our Kudla-Millson Schwartz form
$$
\varphi _{KM}(v,z,M) =   (c'_{11}\Omega + c'_{12}\Omega _{12} +c'_{21}\Omega _{21})\varphi _0(v,z,M)
$$
namely
\begin{align}
\label{ThetaKM}\nne
\Theta _{\varphi _{KM}}(\tau ,z) = \sum_M \varphi _{KM}(v,z,M)e^{2\pi i\det M \tau }.
\end{align}
Here we get as a special instance of the general Kudla-Millson results\\

\nn{\bf Proposition.} \label{PropSiegelKudla}\emph{ $\Theta _{\varphi _{KM}}(\tau ,z)$ is as function in $\tau $ modular for $\SL(2,\mathbb Z)$ of weight }2.\\

{\bf Proof.} This time the Maass operator concerning $\Theta $ is
$$
X_+ = 2iv\partial _\tau  + 1 = iv(\partial _u - i\partial _v) + 1.
$$
If we apply it to our $\Theta (\tau,(i,i)) =: \Theta (\tau)$ we get
$$ 
X_+\Theta (\tau ) = \sum_{M\in M_2(\ZZ)} - \pi(v(a^2+b^2+c^2+d^2 + 2\det M) - 1/\pi ) E(M;\tau ),
$$
i.e., the first terms $c_{11}\varphi _0 = c_{22}\varphi _0$ of our Schwarz form $\varphi _{KM}$ as in (\ref{phiz}). 
And if we apply $X_+$ to the transformation formula
$$
E(M;\tau ) = (1/|\tau |^2)E(M;-1/\tau ) =:(1/|\tau |^2)E^\ast =: F
$$
we get
$$
\begin{array}{rl}
X_+F =& (1/|\tau |^2)(1 - 2iv/\tau )E^\ast + (1/|\tau |^2)((2\pi v\det M -\pi v(a^2+b^2+c^2+d^2)/|\tau |^2 \\
&- (1/|\tau |^2)(2iv/\tau )(-2\pi iu\det M - \pi v((a^2+b^2+c^2+d^2))E^\ast \\[.3cm]
=& (1/\tau ^2)((v/|\tau |^2)(a^2+b^2+c^2+d^2 + 2\det M) - 1/\pi)E(M;-1/\tau ),
\end{array}
$$
i.e., as it should, $X_+\Theta (\tau )$ is modular of weight 2.\\

We don't know whether there is another differential operator producing the other terms from our Schwarz form $\varphi _{KM}.$ 
Hence we have to go the hard way: One has from (\ref{phiz})
$$
c_{12}\varphi _0 = (1/8i)((a+d)^2+(b-c)^2 + 2i(a+d)(c-b))vE(M,\tau )
$$
and, if one takes the other coordinates as above in the discussion of the Siegel (2,2) Theta,
$$ 
c_{12}\varphi _0 = (1/4i)v(x_1+ix_2)^2 e^{\pi i\tau x_1^2}e^{\pi i\tau x_2^2}e^{-\pi i\bar\tau x_3^2}e^{-\pi i\bar\tau x_4^2}.
$$
Here again one can separate the variables and apply the already often used formulae for the Fourier transforms. 
We won't reproduce this but only confirm that one also has modularity of weight 2.\hfill$\Box$ \\

\section{Application to Hirzebruch-Zagier} 

 In the formulation and derivation of the next result we follow a procedure outlined and used in a similar situation in the thesis  of Funke [Fu]. It gives the alternative proof of Corollary \ref{CorthetaTt} announced in Remark \ref{RemthetaTt}.\\
 
\nn {\bf Proposition.}\label{PropthetaTt} \emph{The following identity of generating series
\begin{align}
\lefteqn {- \int_{T(1) }\theta _{\varphi _{KM}} (\tau ,L) =} \nonumber\\& & \sum _{N \geq 0} H_1(N) q^N + 
\sum _{b,c\in\mathbb Z}(1/(8\pi \sqrt v)) \int_1^\infty  e^{-\pi v(b+c)^2r}r^{-3/2}dr)q^{-bc}.\nne
\label{thetaTt}
\end{align}
is an identity of modular forms of weight $2$ for $\SL_2(\ZZ)$.}\\

{\bf Proof.} {\bf 1.} We want to evaluate the integral
$$
\int_{T(1)} \sum_{M\in L}(1/2)\varphi _{KM}(v,z,M)q^{\det M} = \sum_{m\in \mathbb Z} \int_{T(1) } \varphi _{KM}(v,z,m) q^m, \,\,q = e^{2\pi i\tau }.
$$
One has (see (\ref{phireal})
\begin{align*}
\varphi _{KM}(v,(z_1,z_2),M) q^{\det M}&=\\- (1/4)\varphi _\infty (\tau ,(z_1,z_2),M)&\big((v(a'+d')^2 + v(c'-b)^2 - (1/\pi ))(\frac{dx_1\wedge dy_1}{y_1^2}+\frac{dx_2\wedge dy_2}{y_2^2})\nonumber   \\
\!+&(v(a'+d')^2 - v(c'-b')^2)\frac{dx_1\wedge dy_2-dy_1\wedge dx_2}{y_1y_2}\\
 -&2v(a'+d')(c'-b')\frac{dx_1\wedge dx_2+dy_1\wedge dy_2}{y_1y_2}\big)
\end{align*}
with
$$
\varphi _\infty (\tau ,(z_1,z_2),M) = e^{-2\pi (vR((z_1,z_2),M) + \tau \det M}
$$
and (see (\ref{a'+d'}))
\begin{align*}
a'+d' = \,t^{-1}(a-x_2b-x_1c+x_1x_2d +d t^2),\,\,
 c'-b' = \,t^{-1}(y_1(c-x_2d) - y_2(b-x_1d))
\end{align*}

For $z_1=z_2=z$ we get
\begin{align}
\label{phizz}
\varphi _{KM}(v,(z,z),M) q^{\det M}&=- \varphi _\infty (\tau ,(z,z),M)((v(a'+d')^2 + v(c'-b')^2 - (1/2\pi ))\frac{dx\wedge dy}{y^2}\nonumber\\
&= -\varphi _\infty (\tau ,(z,z),M)((1/y^2)v(a-(b+c)x+d|z|^2)^2 -1/2\pi )\frac{dx\wedge dy}{y^2}.\nne
\end{align}
Moreover one has
$$
R((z,z),M) + \det M = (1/2v^2)(a-(b+c)x+d|z|^2)^2 - (ad-(b+c)^2/4) + (c-b)^2/4.
$$
{\bf 2.} Now we compare this to the signature (1,2)-case. The $\Or(1,2)-$case as discussed for instance in [Fu] and \cite{BF} uses a vector space 
$$
V_0 = \{ X \in M_2(\R), {\rm tr} X = 0 \}
$$ 
with quadratic form $q(X) = \det X = -x_1^2 - x_2x_3$  for 
$$
X := \begin{pmatrix} x_1&x_2\\x_3&-x_1 \end{pmatrix}
$$
and coordinates $z =x+iy \in \mathbb H$  stemming from the action of $\SL_2(\R)$ on $V$ 
given by $X \mapsto \gamma X \gamma ^{-1}.$ Hence the associated Schwartz form is (\cite{BF} (3.7) or [Fu] p.296)
\begin{align}
\label{phiw}
\varphi (\tau ,z,X) :=& \varphi _{(1,2)}(\tau ,z,X) \nonumber\\
=& (v(X,X(z)^2 -(1/2\pi )) e^{\pi i(X,X)_{\tau ,z}} \omega \nne
\end{align}
where
$$
\begin{array}{rl}
(X,X(z)) =& -(1/y)(x_3z\bar z-2\tilde x_1x-x_2)\\[.3cm]
(X,X)_{\tau ,z} =& u(X,X) + iv(X,X)_z \\[.3cm]
(X,X)_z =& 2\det X + 2R(z,X)  \\[.3cm]
R(z,X) =& (1/2)(\,X,X(z))^2 - 2\det X  \\[.3cm]
\omega =& dx\wedge dy/y^2.
\end{array} 
$$
To relate this to the formula (\ref{phizz}) above, we put 
$$
M' = \begin{pmatrix} a'&b'\\c'&d' \end{pmatrix} := M \begin{pmatrix} &-1\\1& \end{pmatrix} = X + (b-c)/2E_2
$$
with
$$
x_1 = (b+c)/2, x_2 = -a, x_3 = d, \tilde c = (b-c)/2
$$
and get
$$
\varphi _{KM}(v,(z,z),M) q^{\det M'} = q^{\tilde c}\varphi _{(1,2)}(\tau ,z,X).
$$
and, hence,
\begin{align*}
\sum_{M\in L}\varphi _{KM}(v,(z,z),M) q^{\det M} &= \sum_{M'\in L}\varphi _{KM}(v,(z,z),M') q^{\det M}\\
&= \sum_{M'\in L}q^{\tilde c}\varphi _{(1,2)}(\tau ,z,X).
\end{align*}
Here we can use calculations from part 3 and 5 of Funke's thesis [Fu]. 

For each $M' \in L = M_2(\Z)$ we can write
$$
M' = (\hat c+(j/2))E_2 + \begin{pmatrix} \hat b+(j/2)&-a\\d&-(\hat b+(j/2)) \end{pmatrix}
$$
with $j=0$ or $j=1$ and $\hat b,\hat c,a,d \in \mathbb Z.$
Our lattice $L = M_2(\mathbb Z)$ decomposes as
$$
L = L_0 + L_1 = \sum_{j=0,1} (\{(\hat c+(j/2))E_2,\,\hat d\in\mathbb Z\} + K_j )
$$
with 
$$
K_j := (j/2)X_1 + K\,,\,\, K := \{ \begin{pmatrix} \hat b&-a\\d&-\hat b \end{pmatrix}, \,\hat b,a,d\in\mathbb Z\},\,
X_1 = \begin{pmatrix}1&\\&-1 \end{pmatrix}.
$$
{\bf 3.} Now  here we have a special example of [Fu] Theorem 5.1 for $p=2$ 
and we get

\begin{align}\label{thetaT}
I_{\varphi _{KM}} (\tau ,L,T(1) ) &= \int_{T(1) }\theta _{\varphi _{KM}} (\tau ,L)\nonumber\\[.3cm]
&= (1/2)\sum_{j=0,1}\vartheta _j(\tau ) I_{\varphi _{U^\bot }}(\tau ,L_j),\,\, \vartheta _j(\tau ) = \sum q^{(m+(j/2))^2}.\nne
\end{align}

Here $\varphi _{U^\bot } = \varphi _{(1,2)}$ is the Schwartz form for the three-dimensional vector space with quadratic form 
$q(X) = -\tilde b^2 +ad$ of signature (1,2). As Funke kindly remarked, here we can take over the calculation in the proof of 
(the Hirzebruch-Zagier-) Theorem 5.4 in [Fu] for the special case $D=1.$ From [Fu] (3.41) one has
$$
(1/2)I_{\varphi _{U^\bot }}(\tau ,K_j) = \sum_{N=0}^\infty H(4N-j)q^{N-j/4} + (1/(8\pi \sqrt v))\sum_{n\in\mathbb Z}\beta (4\pi(n+(j/2))^2v)q^{-(n+(j/2))^2}
$$
with $\beta (u) := \int_1^\infty e^{-u t}t^{-3/2}dt = B(u)$ and, for $N>0,$ $H(N)$ the  class number of binary positive definite integral quadratic forms 
of discriminant $-N$ and  
$$
H(0) = (1/2){\rm vol}\, (\,\SL(2,\Z)\setminus \mathbb H\,) = -1/12.
$$
 
{\bf 4.} Using this result, in (\ref{thetaT}) we have as coefficients of $q^N, N>0$ the sum over 
$$
j=0,1,\,\, (m+(j/2))^2 + (n-(j/4)) = N \,\,{\rm or}\,\,4n-j = 4N - (2m+j)^2
$$
i.e.
$$
\sum_{j=0,1}\sum_{(2m+j)2\leq 4N} H(4N-(2m+j)^2) = \sum_{s2\leq 4N} H(4N-s^2) = H_1(N)
$$
which by inspection of \cite{Hi} p.2 is the same as the number $H_1(N)$ counting the intersections of $T(1) $ and $T(N)$ in the interior of $X.$ 
and hence the holomorphic term
$$
P(\tau ) = \sum _N H_1(N) q^N.
$$

{\bf 5.} Moreover, in (\ref{thetaT}) we have as coefficients of $q^r, r\in \mathbb Z$ the sum over  
$$
j=0,1\,,\,\,(m+(j/2))^2 - (n+(j/2))^2 = r
$$
i.e.
$$
1/(4\pi \sqrt v)(\sum_{m^2-n^2=r} \beta (4\pi n^2v) + \sum_{(m+(1/2))^2-(n+(1/2))^2=r} \beta (4\pi (n+(1/2))^2v))
$$
Putting here $m=\hat b-\hat c, \,n=\hat b+\hat c$ where $\hat b = b,\, \hat c=c$ and  $\hat b = (1/2)+b,\, \hat c= (1/2)+c$ with integer $b,c$ 
we get the coefficient
$$
1/(8\pi \sqrt v)(\sum_{b,c\in\mathbb Z,\,-bc=r} \beta (\pi (b+c)^2v)
$$

and hence the claim follows.
\hfill $\Box$





\addcontentsline{toc}{section}{\bf Bibliography}
\begin{thebibliography}{AKuKueneu1}
\bibitem[AS] {AS}Abramovitz, M., Stegun, I.A.: {\it Handbook of mathematical functions.} US National bureau of standards 1972.
\bibitem[BJ] {BJ}Borel, A., Ji, L.: {\it Compactifications of Symmetric and Locally Symmetric Spaces.} Birkh\"auser, Boston 2006.
\bibitem[BF] {BF}Bruinier, J.H., Funke, J.: {\it Traces of CM values of modular functions.} J. reine angew. Math. {\bf 594} (2006) 1-33.
\bibitem[BBK] {BBK}Bruinier, J.H., Burgos Gil, J.I., K\"uhn, U.: {\it Borcherds Products and Arithmetic Intersection Theory on Hilbert Modular Surfaces.}Duke Math. J. {\bf 139} (2007) 1-88.
\bibitem[BKK] {BKK} Burgos Gil, J.I., Kramer, J., K\"uhn, U.: {\it Cohomogical arithmetic Chow rings.} 
J. Inst. Math. Jussieu {\bf 6} (2007) 1-172.
\bibitem[Fu]{Fu}Funke, J,: {\it Heegner Divisors and Nonholomorphic Modular Forms} 
Compositio Math. {\bf 133} (2002) 289-321. 
\bibitem[FM1]{FM1}Funke, J., Millson, J.: {\it Boundary behavior of special cohomology classes arising from the Weil 
representation.} Preprint 2009.
\bibitem[FM2]{FM2}Funke, J., Millson, J.: {\it The geometric theta correspondence for Hilbert modular surfaces.} Preprint 2011.
\bibitem[FM3]{FM3}Funke, J., Millson, J.: {\it On a theorem by Hirzebruch and Zagier.} forthcoming.
\bibitem[GR]{GR} G\"ortz, U., Rapoport, M.(ed.): {\it ARGOS Seminar on Intersections of Modular Correspondences.} Ast\'erisque {\bf 312} 2007.
\bibitem[Hi] {Hi}Hirzebruch,F.: {\it Kurven auf den Hilbertschen Modulflaechen und Klassenzahlrelationen.}
pp.75-93 in LNM {\bf 412} Springer, Berlin 1974.
\bibitem[HZ] {}Hirzebruch,F., Zagier, D.: {\it Intersection numbers of curves on Hilbert Modular surfaces and modular forms of nebentypus.}
Inventiones math. {\bf 36} (1978) 57-113.
\bibitem[Iw] {Iw}Iwaniec, H.: {\it Spectral Methods of Automorphic Forms.} GSM 53, AMS 2002.
\bibitem[Kn] {Kn}Knapp, A.W.: {\it Elliptic Curves.} Princeton University Press 1992.
\bibitem[Ku1]{Ku1}Kudla, S.: {\it Integrals of Borcherds forms.} Compositio math. {\bf 137} (2003) 239-349. 
\bibitem[Ku2]{Ku2}Kudla, S.: {\it Central Derivatives of Eisenstein series and Height pairings.} Ann. math. {\bf 146} (2003) 545-646.
\bibitem[Ku3]{Ku3}Kudla, S.: {\it Special Cycles and derivatives of Eisenstein series.} p. 243-270 in {\it Heegner Points and 
Rankin L-Series.} Math. Sci. Res. Inst. Publ. {\bf 49} Cambridge Univ. Press, Cambridge, 2004.
\bibitem[KM1]{KM1}Kudla, S., Millson, J.: {\it The theta correspondence and harmonic forms I.} Math. Ann. {\bf 247} (1986) 353-378.
\bibitem[KM2]{KM2}Kudla, S., Millson, J.: {\it The theta correspondence and harmonic forms II.} Math. Ann. {\bf 277} (1987) 267-314.
\bibitem[KM3]{KM3}Kudla, S., Millson, J.: {\it Intersection numbers of cycles on locally symmetric spaces and
Fourier coefficients of holomorphic modular forms in several complex variables.} IHES Publ. Math.. {\bf 71} (1990) 121-172.
\bibitem[Ma]{Ma}Maass, H.: {\it Die Differentialgleichungen in der Theorie der elliptischen Modulfunktionen.} Math. Ann. {\bf 125} (1953) 235-263.
\bibitem[Mo]{Mo}Moore, C.N.: {\it On the Summability of the Double Fourier's Series of Discontinuous Functions.} Math. Ann. {\bf 74} (1913) 555-572. 
\bibitem[Mu]{Mu}Mumford. D.: {\it Tata Lectures on Theta I.}\, PM {\bf 74} Birkh\"auser, Boston 1983. 

\bibitem[Og] {Og}Ogg, A.: {\it Modular forms and Dirichlet Series.} Benjamin, New York 1969.
\bibitem[Ran] {Ran}Rankin, R.A.: {\it Modular Forms and Functions.} 
Cambridge Univerity Press, 1977 
\bibitem[S1] {S1}Siegel, C. L.: {\it \"Uber die Zetafunktionen indefiniter quadratischer Formen II.} Math. Z. {\bf 44} (1939) 398-426.
\bibitem[S2] {S2}Siegel, C. L.: {\it Indefinite quadratische Formen und Funktionentheorie I.} Math. Ann. {\bf 124} (1951) 17-54.
\bibitem[SABK] {SABK}Soul\'e, C, Abramovich, D., Burnol, J.-F., Kramer, J.: {\it Lectures on Arakelov Geometry.} Cambridge University Press 1992.
\bibitem[Sou] {Sou}Soul\'e, C.: {\it Hermitian vector bundles on arithmetic varieties.} Proc. Symp. Pure Math. {\bf 62.I} (1997) 383-419.
\end{thebibliography}
\end{document}